\documentclass[12pt]{amsart}

\usepackage{amscd,verbatim}
\usepackage{amssymb}
\usepackage[all]{xy}
\usepackage[colorlinks,linkcolor=blue,citecolor=blue,urlcolor=red]{hyperref}
\usepackage{bm}

\newcommand{\A}{\mathbf{A}}

\newcommand{\G}{\mathbf{G}}

\renewcommand{\P}{\mathbf{P}}

\newcommand{\Z}{\mathbb{Z}}

\newcommand{\sK}{\mathcal{K}}

\newcommand{\sO}{\mathcal{O}}

\newcommand{\sU}{\mathcal{U}}
\newcommand{\sV}{\mathcal{V}}
\newcommand{\sW}{\mathcal{W}}
\newcommand{\sX}{\mathcal{X}}
\newcommand{\sY}{\mathcal{Y}}

\newcommand{\bZ}{\mathbb{Z}}
\newcommand{\fm}{\mathfrak{m}}

\newcommand{\Xb}{{\overline{X}}}

\newcommand{\Yb}{{\overline{Y}}}

\newcommand{\Cor}{\operatorname{\mathbf{Cor}}}

\newcommand{\HI}{\operatorname{\mathbf{HI}}}

\newcommand{\Ext}{\operatorname{Ext}}
\newcommand{\ul}[1]{{\underline{#1}}}

\newcommand{\PST}{{\operatorname{\mathbf{PST}}}}
\newcommand{\PS}{{\operatorname{\mathbf{PS}}}}

\newcommand{\NST}{\operatorname{\mathbf{NST}}}

\newcommand{\Hom}{\operatorname{Hom}}
\newcommand{\uHom}{\operatorname{\underline{Hom}}}

\newcommand{\Ker}{\operatorname{Ker}}

\newcommand{\Coker}{\operatorname{Coker}}

\newcommand{\Pic}{\operatorname{Pic}}

\newcommand{\Spec}{\operatorname{Spec}}

\newcommand{\Comp}{\operatorname{Comp}}
\newcommand{\Sm}{\operatorname{\mathbf{Sm}}}
\newcommand{\tSm}{\widetilde{\Sm}}
\newcommand{\Sch}{\operatorname{\mathbf{Sch}}}

\newcommand{\Ab}{\operatorname{\mathbf{Ab}}}

\newcommand{\tr}{{\operatorname{tr}}}

\newcommand{\dlog}{{\operatorname{dlog}}}

\newcommand{\fin}{{\operatorname{fin}}}

\newcommand{\red}{{\operatorname{red}}}

\newcommand{\Nis}{{\operatorname{Nis}}}
\newcommand{\et}{{\operatorname{\acute{e}t}}}

\newcommand{\inj}{\hookrightarrow}

\newcommand{\Res}{\operatorname{Res}}

\newcommand{\ch}{{\operatorname{ch}}}

\newcommand{\Frac}{{\operatorname{Frac}}}

\renewcommand{\lim}{\operatornamewithlimits{\varprojlim}}
\newcommand{\colim}{\operatornamewithlimits{\varinjlim}}

\newcommand{\ol}{\overline}

\renewcommand{\phi}{\varphi}
\renewcommand{\epsilon}{\varepsilon}
\renewcommand{\div}{\operatorname{div}}

\newcommand{\MNST}{\operatorname{\mathbf{MNST}}}

\newcommand{\MCor}{\operatorname{\mathbf{MCor}}}

\newcommand{\MPST}{\operatorname{\mathbf{MPST}}}
\newcommand{\CI}{\operatorname{\mathbf{CI}}}

\newcommand{\bcube}{{\ol{\square}}}
\newcommand{\cube}{\square}

\newcommand{\Zb}{{\overline{Z}}}

\def\rmapo#1{\overset{#1}{\longrightarrow}}

\newcommand{\ulMPST}{\operatorname{\mathbf{\underline{M}PST}}}
\newcommand{\ulMNST}{\operatorname{\mathbf{\underline{M}NST}}}
\newcommand{\ulMCor}{\operatorname{\mathbf{\underline{M}Cor}}}
\newcommand{\ulomega}{\underline{\omega}}
\newcommand{\ulomegaCI}{\underline{\omega}^{\CI}}

\newcommand{\MV}{\operatorname{MV}}

\def\bZ{\mathbb{Z}}
\def\Ztr{\bZ_\tr}

\def\isom{\overset{\simeq}{\longrightarrow}}

\newcommand{\Vb}{\overline{V}}

\def\Vb{\overline{V}}
\def\Wb{\overline{W}}

\newcounter{spec}
{\end{list}}%
\def\otMPST{\otimes_{\MPST}}
\def\otuMPST{\otimes_{\ulMPST}}

\def\hcubesp{h_0^{\bcube,sp}}
\def\hcube{h_0^{\bcube}}
\def\hcubespNis{h_{0,\Nis}^{\bcube,sp}}

\def\aNisCI{\aNis^{\CI}}
\newcommand{\eq}[2]{\begin{equation}\label{#1}#2 \end{equation}}

\newtheorem{lemma}{Lemma}[section]
\newtheorem{thm}[lemma]{Theorem}

\newtheorem{prop}[lemma]{Proposition}

\newtheorem{cor}[lemma]{Corollary}

\theoremstyle{definition}

\theoremstyle{remark}
\newtheorem{qn}[lemma]{Question}

\newtheorem{remark}[lemma]{Remark}

\newtheorem{claim}[lemma]{Claim}

\numberwithin{equation}{section}

\def\aNis{a_{\Nis}}
\def\ulaNis{\underline{a}_{\Nis}}

\def\qaq{\quad\text{ and }\quad}

\def\Comp{\Comp^{\fin}}

\def\ulMNST{\operatorname{\mathbf{\ul{M}NST}}}

\def\MNST{\operatorname{\mathbf{MNST}}}

\def\RSC{\operatorname{\mathbf{RSC}}}

\def\Comp{\operatorname{\mathbf{Comp}}}
\def\uli{\ul{i}}

\def\qfor{\text{ for }\;\;}
\def\CIt{\operatorname{\mathbf{CI}}^\tau}
\def\CItsp{\operatorname{\mathbf{CI}}^{\tau,sp}}
\def\CIsp{\operatorname{\mathbf{CI}}^{sp}}
\def\CIspNis{\operatorname{\mathbf{CI}}^{sp}_{\Nis}}
\def\CItspNis{\operatorname{\mathbf{CI}}^{\tau,sp}_{\Nis}}
\def\CItNis{\operatorname{\mathbf{CI}}^{\tau}_{\Nis}}
\def\otCIsp{\otimes_{\CI}^{sp}}
\def\otCINissp{\otimes_{\CI}^{\Nis,sp}}
\def\otCIspNis{\otimes_{\CI}^{\Nis,sp}}

\begin{document}
\def\ihF#1{F^{#1}}
\def\ihFA{\ihF {\A^1}}

\def\istm{\iota_{st,m}}
\def\im{\iota_m}
\def\est{\epsilon_{st}}
\def\tL{\tilde{L}}
\def\tX{\tilde{X}}
\def\tY{\tilde{Y}}
\def\omegaCI{\omega^{\CI}}
\def\qwith{\;\text{ with} }
\def\aVNis{a^V_\Nis}
\def\ulMCorls{\ulMCor_{ls}}

\def\Zinf{Z_\infty}
\def\Einf{E_\infty}
\def\Xinf{X_\infty}
\def\Yinf{Y_\infty}
\def\Pinf{P_\infty}

\def\Lot{{\cubegm\otimes\cubegm}}

\def\Ln#1{\Lambda_n^{#1}}
\def\tLn#1{\widetilde{\Lambda_n^{#1}}}
\def\tild#1{\widetilde{#1}}
\def\otuCINis{\otimes_{\underline{\CI}_\Nis}}
\def\otCI{\otimes_{\CI}}
\def\otCINis{\otimes_{\CI}^{\Nis}}
\def\tF{\widetilde{F}}
\def\tG{\widetilde{G}}
\def\cubegm{\bcube^{(1)}}
\def\cubegma{\bcube^{(a)}}
\def\cubegmb{\bcube^{(b)}}
\def\cubegmred{\bcube^{(1)}_{red}}
\def\cubegmreda{\bcube^{(a)}_{red}}
\def\cubegmredb{\bcube^{(b)}_{red}}

\def\LT{\bcube^{(1)}_{T}}
\def\LU{\bcube^{(1)}_{U}}
\def\LV{\bcube^{(1)}_{V}}
\def\LW{\bcube^{(1)}_{W}}
\def\LTred{\bcube^{(1)}_{T,red}}
\def\Lred{\bcube^{(1)}_{red}}
\def\LTred{\bcube^{(1)}_{T,red}}
\def\LUred{\bcube^{(1)}_{U,red}}
\def\LVred{\bcube^{(1)}_{V,red}}
\def\LWred{\bcube^{(1)}_{W,red}}
\def\PP{\P}
\def\AA{\A}

\def\LL{\bcube^{(2)}}
\def\LLred#1{\bcube^{(2)}_{#1,red}}
\def\LLredd{\bcube^{(2)}_{red}}
\def\Lredd#1{\bcube_{#1,red}}

\def\Lnredd#1{\bcube^{(#1)}_{red}}

\def\LLT{\bcube^{(2)}_T}
\def\LLTred{\bcube^{(2)}_{T,red}}

\def\LLU{\bcube^{(2)}_U}
\def\LLUred{\bcube^{(2)}_{U,red}}

\def\LLS{\bcube^{(2)}_S}
\def\LLSred{\bcube^{(2)}_{S,red}}
\def\tMCor{\Hom_{\MPST}}
\def\otHINis{\otimes_{\HI}^{\Nis}}

\title{Cancellation theorems for reciprocity sheaves}
\author[A. Merici]{Alberto Merici}
\address{Insutitut f\"ur Mathematik\\
	Universit\"at Z\"urich\\
	Winterthurerstrasse 190, CH-8057 Z\"urich\\ Switzerland}
\email{alberto.merici@gmail.com}
\author[S. Saito]{Shuji Saito}
\address{Graduate School of Mathematical Sciences\\
University of Tokyo\\
3-8-1 Komaba, Tokyo 153-8941\\ Japan}
\email{sshuji@msb.biglobe.ne.jp}
\thanks{The first author is supported by the Swiss National Science Foundation (SNF), project 200020\_178729}
\thanks{The second author is supported by JSPS KAKENHI Grant (15H03606).}

\begin{abstract}
    We prove cancellation theorems for reciprocity sheaves and cube-invariant modulus sheaves with transfers of Kahn--Miyazaki--Saito--Yamazaki.
It generalizes a cancellation theorem for $\A^1$-invariant sheaves with transfers,
which was proved by Voevodsky. As an application, we get some new formulas for internal hom's of the sheaves  $\Omega^i$ of absolute K\"ahler differentials.
\end{abstract}

\subjclass[2020]{19E15 (14F42, 19D45, 19F15)}
\maketitle
\tableofcontents
\setcounter{section}{-1} 

\section{Introduction}

We fix once and for all a perfect field $k$. 
Let $\Sm$ be the category of separated smooth schemes of finite type over $k$.
Let $\Cor$ be the category of finite correspondences:
$\Cor$ has the same objects as $\Sm$ and morphisms in $\Cor$ are finite correspondences. 
Let $\PST$ be the category of additive presheaves of abelian groups on $\Cor$, called presheaves with transfers. 
Let $\NST\subset \PST$ be the full subcategory of Nisnevich sheaves, i.e.
those objects $F\in \PST$ whose restrictions $F_X$ to the small \'etale site $X_{\et}$ over $X$ are Nisnevich sheaves for all $X\in \Sm$. 
By a fundamental result of Voevodsky, the inclusion $\NST\to \PST$ has an exact left adjoint $a^V_\Nis$ such that for any $F\in \PST$ and $X\in \Sm$, 
$(a^V_\Nis F)_X$ is the Nisnevich sheafication of $F_X$ as a presheaf on $X_\Nis$. 
In Voevodsky's theory of motives, a fundamental role is played by $\A^1$-invariant objects $F\in \NST$, namely such $F$ that $F(X)\to F(X\times\A^1)$ induced by the projection $X\times\A^1\to X$ are isomorphisms for all $X\in \Sm$.
The $\A^1$-invariant objects form a full abelian subcategory $\HI_\Nis\subset\NST$ that carries a symmetric monoidal structure $\otHINis$ such that
\[F\otHINis G=h_0^{\A^1,\Nis}\aVNis(F\otimes_{\PST} G) \;\text{ for } F,G\in \HI_\Nis,\] 
where $\otimes_{\PST}$ is the symmetric monoidal structure on $\PST$ induced formally from that on $\Cor$ and $h_0^{\A^1,\Nis}$ is a left adjoint to the inclusion functor $\HI_\Nis\to \NST$, which sends an object of $\NST$ to its maximal $\A^1$-invariant quotient in $\NST$. 
For integers $n>0$, the twists of $F\in \HI_\Nis$ are then defined as
\[ F(1) = F\otHINis \G_m,\quad F(n) := F(n-1)\otHINis \G_m;\]
where $\G_m\in \NST$ is given by $X\to \Gamma(X,\sO^\times)$ for $X\in \Sm$.

Noting that $-\otHINis \G_m$ is an endo-functor on $\HI_\Nis$, we get a natural map:
\begin{equation}\label{eq;cancelHI}
 \iota_{F,G}: \Hom_{\PST}(F,G) \to \Hom_{\PST}(F(1),G(1))\;\text{ for } F,G\in \HI_\Nis.
\end{equation}
One key ingredient in Voevodsky's theory is the Cancellation theorem \cite[Cor, 4.10]{voecancel}, which implies the following theorem:

\begin{thm}\label{voecancel}
    For $F,G\in \HI_\Nis$, $\iota_{F,G}$ is an isomorphism.
\end{thm}

\medbreak

The purpose of this paper is to generalize the above theorem to reciprocity sheaves. 
The category $\RSC_{\Nis}$ of reciprocity sheaves was introduced in \cite{ksyI} and \cite{ksyII} as a full subcategory of $\NST$ that contains $\HI_{\Nis}$ as well as interesting non-$\A^1$-invariant objects such as the additive group scheme $\G_a$, the sheaf of absolute K\"ahler differentials $\Omega^i$ and the de Rham-Witt sheaves $W_n\Omega^i$. 
In \cite{rsy}, a lax monoidal structure $(\_\;,\_)_{\RSC_\Nis}$ on $\RSC_\Nis$ 
is defined in such a way that 
\[(F,G)_{\RSC_\Nis}=F\otHINis G\;\text{ for } F,G\in \HI_\Nis.\]
It allows us to define the twists for $F\in \RSC_\Nis$ recursively as\[
F\langle 1 \rangle := (F,\G_m)_{\RSC_\Nis},\quad F\langle n\rangle := (F\langle n-1\rangle,\G_m)_{\RSC_\Nis}.\] 
Some examples of twists were computed in \cite{rsy}:
if $F\in \HI_{\Nis}$, then $F\langle n \rangle=F(n)$, in particular 
$\Z\langle n\rangle\cong \mathcal{K}_n^M$ (the Milnor $K$-sheaf), 
and $\G_a\langle n\rangle\cong \Omega^n$ if $\ch(k)=0$.

By the fact that $(-,\G_m)_{\RSC_\Nis}$ is an endo-functor on $\RSC_\Nis$, we get a natural map (cf. \eqref{eq;cancelRSC}) :
\begin{equation}\label{eq;cancelRSCintro}
 \iota_{F,G}: \Hom_{\PST}(F,G) \to \Hom_{\PST}(F\langle 1 \rangle,G\langle 1 \rangle)
\;\text{ for } F,G\in \RSC_\Nis,
\end{equation}
which coincides with \eqref{eq;cancelHI} if $F,G\in \HI_\Nis$.
We will also get a natural map in $\NST$:
\begin{equation}\label{iotaFGmintro}
\lambda_F: F \to \uHom_{\PST}(\sK_n^M,F\langle n\rangle)\;\text{ for } F\in \RSC_\Nis,
\end{equation}
using the functoriality of $(-,\G_m)_{\RSC_\Nis}$, where $\uHom_{\PST}$ denotes the internal hom in $\PST$.

The main result of this paper is the following:

\begin{thm}[Theorems \ref{thm;closeenoughtocancel} and \ref{cor;cancelmilnorK}]\label{thm;main1}
    The maps $\iota_{F,G}$ and $\lambda_F$ are isomorphisms.
\end{thm}

As an application of the above theorem, we prove the following.

\begin{cor}\label{cor;cancelOmegaintro}(Theorem \ref{pr;cancelOmega})
Assume $\ch(k)=0$. For integers $m,n \geq 0$, there are natural isomorphisms in $\NST$:
\[\uHom_{\PST}(\Omega^n,\Omega^m)\cong \Omega^{m-n}\oplus \Omega^{m-n-1},\]
\[\uHom_{\PST}(\mathcal{K}_n^M,\Omega^m)\cong \Omega^{m-n},\]
where $\Omega^i=0$ for $i<0$ by convention.
\end{cor}

Let $\PS$ be the category of additive presheaves of abelian groups on $\Sm$ (without transfers).
Note that $\PST$ is viewed as a subcategory of $\PS$. 
By a lemma due to Kay R\"ulling (see Lemma \ref{lem;MapPSvsPST}), we have a natural isomorphism in $\PS$:
\begin{equation}\label{eq;homuhom}
\uHom_{\PST}(G,\Omega^m)\cong \uHom_{\PS}(G,\Omega^m)
\;\text{ for any } G\in \PST,
\end{equation}
where $\uHom_{\PS}$ is the internal hom in $\PS$. 
Thanks to \eqref{eq;homuhom}, the isomorphisms of Corollary \ref{cor;cancelOmegaintro} and its explicit descriptions \eqref{uPhi} and \eqref{uPsi} imply 
\[  \Hom_{\PS}(\Omega^n,\Omega^m) =
\{\omega_1 \wedge (-) + \omega_2\wedge d(-)\;|\; 
\omega_1\in \Omega^{m-n}_k,\; \omega_2\in \Omega^{m-n-1}_k\},\]
\[  \Hom_{\PS}(\mathcal{K}_n^M,\Omega^m) =\{\omega \wedge \dlog (-) \;|\; 
\omega\in \Omega^{m-n}_k\},\]
where $\dlog: \mathcal{K}_n^M \to\Omega^m$ is the map
$\{x_1,\dots,x_n\} \to \dlog x_1\wedge\cdots\wedge \dlog x_n$.
It would be an interesting question if there is a direct proof of these formulas which does not use the machinery of modulus sheaves with transfers explained below.
\bigskip

Reciprocity sheaves are closely related to \emph{modulus sheaves with transfers}
introduced in \cite{kmsyI} and \cite{kmsyII}: Voevodsky's category $\Cor$ of finite correspondences is enlarged to a new category $\ulMCor$ of \emph{modulus pairs}: Its objects are pairs $\sX=(X,D)$ where $X$ is a separated scheme of finite type over $k$ and $D$ is an effective Cartier divisor on $X$ such that $\sX^\circ:=X-|D|\in \Sm$ ($\sX^\circ$ is called the interior of $\sX$).
The morphisms are finite correspondences on interiors satisfying some admissibility and properness conditions. 
Let $\MCor\subset \ulMCor$ be the full subcategory of such objects $(X,D)$
that $X$ is proper over $k$. 
There is a symmetric monoidal structure $\otimes$ on $\ulMCor$, which also induces that on $\MCor$ by restriction (cf. \S\ref{recollection}\eqref{tensor}) .

We then define $\ulMPST$ (resp. $\MPST$) as the category of additive presheaves of abelian groups on $\ulMCor$ (resp. $\MCor$). 
We have a functor 
\[\ulomega:\ulMCor \to \Cor\;;\; (\Xb,\Xinf) \to \Xb - |\Xinf|,\]
and two adjunctions
\begin{equation*}\label{eq;tauomegaadjunction-intro}
\MPST\begin{smallmatrix} \tau^*\\ \longleftarrow\\ \tau_!\\ \longrightarrow\\
\end{smallmatrix}\ulMPST, 
\quad
\ulMPST
\begin{smallmatrix} \ulomega^*\\ \longleftarrow\\ \ulomega_!\\ \longrightarrow\\
\end{smallmatrix}\PST,
\end{equation*}
where $\ulomega^*$ is induced by $\ulomega$ and $\ulomega_!$ is its left Kan extension, and $\tau^*$ is induced by the inclusion $\tau:\MCor\to \ulMCor$ and $\tau_!$ is its left Kan extension, which turned out to be exact and fully faithful.

For $F\in \ulMPST$ and $\sX=(X,D)\in \ulMCor$ write $F_{\sX}$ for the presheaf
on the small \'etale site $X_{\et}$ over $X$ given by $U\to F(\sX_U)$ for $U\to X$ \'etale, where $\sX_U=(U,D\times_{X}U)\in \ulMCor$.
We say $F$ is a Nisnevich sheaf if so is $F_{\sX}$ for all $\sX\in \ulMCor$.
We write $\ulMNST\subset \ulMPST$ for the full subcategory of Nisnevich sheaves. 

The replacement of the $\A^1$-invariance in this new framework is 
the $\bcube$-invariance, where $\bcube:=(\P^1,\infty)\in \MCor$:
Let $\CI\subset \MPST$ be the full subcategory of those objects
$F$ that $F(\sX) \to F(\sX\otimes \bcube)$ induced by the projection 
$\sX\otimes\bcube \to \sX$ are isomorphisms for all 
$\sX\in \MCor$. Let $\CIt\subset \ulMPST$ be the essential image of $\CI$ under $\tau_!$ and $\CItsp\subset \CIt$ be the full subcategory of \emph{semipure} objects $F$, namely such objects that the natural map
$F(X,D)\to F(X-D,\emptyset)$ are injective for all $(X,D)\in \ulMCor$.
We also define $\CItsp_\Nis=\CItsp\cap \ulMNST$ as a full subcategory of $\ulMNST$.
A symmetric monoidal structure $\otCIsp$ (resp. $\otCINissp$)
on $\CItsp$ (resp. on $\CItsp_\Nis$) can be defined in the same spirit as $\otHINis$ (see \S\ref{sec3}).

The relationship between reciprocity (pre)sheaves and $\bcube$-invariant modulus (pre)sheaves with transfers is encoded in 
\[ \RSC =\ulomega_!(\CItsp)\qaq \RSC_\Nis =\ulomega_!(\CItsp_\Nis).\]
There is a pair of adjoint functors 
\[
\CItsp\begin{smallmatrix} \ulomega^{\CI}\\ \longleftarrow\\ \ulomega_!\\ \longrightarrow\\ \end{smallmatrix}\RSC\qaq
\CItsp_\Nis\begin{smallmatrix} \ulomega^{\CI}\\ \longleftarrow\\ \ulomega_!\\ \longrightarrow\\ \end{smallmatrix}\RSC_\Nis\]
such that $\ulomega^{\CI} F =\ulomega^* F$ for $F\in \HI$.
Moreover, the lax monoidal structure on $\RSC_{\Nis}$ is induced by the symmetric monoidal structure on 
$\CItsp_{\Nis}$ via the formula:
\[ (F,G)_{\RSC_\Nis}:= \ulomega_!(\ulomegaCI F\otCIspNis \ulomegaCI G)\qfor F,G\in \RSC_\Nis.\] 
The endo-functor $-\otCINissp \ulomega^*\G_m$ on $\CItsp_\Nis$ induces a natural map for $F\in \CItsp_\Nis$:
\eq{iotaFNisintro}{\iota_{F}:  F\to \uHom_{\ulMPST}(\ulomega^*\G_m, F \otCINissp\ulomega^*\G_m),}
where $\uHom_{\ulMPST}$ denotes the internal hom in $\ulMPST$. 
Now Theorem \ref{thm;main1} will be a consequence of the following result:

\begin{thm}[Cor \ref{cor3;weakcancel}]\label{thm;main2}
  For $F\in \RSC_\Nis$ and $\tF=\ulomega^{\CI}F \in \CItsp_\Nis$, the map
$\iota_{\tF}$ is an isomorphism.
\end{thm}
\bigskip

We give an outline of the content of the paper:
\begin{itemize}
    \item In section \ref{recollection} we first review basic definitions and results of the theory of modulus (pre)sheaves with transfers and reciprocity sheaves from \cite{kmsyI}, \cite{kmsyII} and \cite{shuji}. 
   We also prove some technical lemmas which will be used in 
the later sections.
   
    \item In section \ref{sec2} we define the contraction functors $\gamma$ on $\CItsp$ and $\CItsp_\Nis$, which generalize Voevodsky's contraction functors on $\HI$ and $\HI_\Nis$ (cf. \cite[Lecture 23]{mvw}) to the setting of modulus (pre)sheaves with transfers. We prove some technical lemmas which will be used in the later sections.
    
    \item In section \ref{sec3} we define the symmetric monoidal structure $\otCIsp$ (resp. $\otCINissp$) on $\CItsp$ (resp. on $\CItsp_\Nis$) using results from section \ref{recollection}. 
The endo-functor $-\otCIsp \ulomega^*\G_m$ on $\CItsp$ induces a natural map for $F\in \CItsp$:
\eq{iotaFintro}{\iota_{F}:  F\to \uHom_{\ulMPST}(\ulomega^*\G_m, F \otCIsp\ulomega^*\G_m).}
We state the main Theorem \ref{thm;weakcancel}: $\iota_F$ is an isomorphism.
Theorem \ref{thm;main2} is deduced from it by using results from sections \ref{sec2}.

The last half of the section is devoted to the proof of the split-injectivity of the map
$\iota_F$ \eqref{iotaFintro}. In order to construct a section of $\iota_F$, we follow the same strategy as \cite{voecancel} by generalizing the techniques used in loc. cite. 
    
    \item In section \ref{sec4} we finish the proof of Theorem \ref{thm;weakcancel} by showing the surjectivity of $\iota_F$. We again follow the same strategy as \cite{voecancel} by generalizing the results of \cite[Section 2.7]{gramot}: here a technical problem is that for $(X,D)\in \ulMCor$, the diagonal map $X\to X\times X$ does not induce a map $(X,D)\to (X,D)\otimes(X,D)$ in $\ulMCor$ but only induces a map
$(X,2D)\to (X,D)\otimes(X,D)$, where $2D\hookrightarrow X$ is the thickening of 
$D\hookrightarrow X$ defined by the square of the ideal sheaf.
This is the main reason why we need to work with $\CItsp$ instead of $\CIt$ 
employing much more intricate arguments than those in \cite{voecancel} and 
\cite[Section 2.7]{gramot}, for which we need the technical results in \S\ref{recollection} and \S\ref{sec2}.

    \item In section \ref{implicationRSC} we deduce Theorem \ref{thm;main1} from Theorem \ref{thm;main2}.

    \item In section \ref{intHomOmega} we deduce Corollary \ref{cor;cancelOmegaintro} from Theorem \ref{thm;main1}. 
\end{itemize}

\subsection*{Acknowledgements} The authors would like to thank Kay R\"ulling for letting them include his proofs of Lemmas \ref{lem;MapPSvsPST} and \ref{lem;Kay} and for pointing out a mistake in the first version of this paper and also for many valuable comments.
The first author would like to thank his PhD supervisor Joseph Ayoub for suggesting the study of modulus sheaves and for many discussions that led to the formulation of the results in this paper. He would like to thank Lorenzo Mantovani and Federico Binda for many helpful discussions.
The second author would like to thank Joseph Ayoub for the invitation to the university of Z\"urich where his collaboration with the first author started.

\subsection*{Conventions}
In the whole paper we fix a perfect base field $k$. 
Let $\tSm$ be the category of $k$-schemes $X$ which are essentially smooth over $k$, i.e. $X$ is a limit $\lim_{i \in I} X_i$ over a filtered set $I$, 
where $X_i $ is smooth over $k$ and all transition maps are \'etale.
Note $\Spec K\in \tSm$ for a function field $K$ over $k$ thanks to the assumption that $k$ is perfect.
We frequently allow $F\in \PST$ to take values 
on objects of $\widetilde{\Sm}$ by
$F(X) := \colim_{i \in I} F(X_i)$ for $X$ as above.
 
\section{Recollection on modulus sheaves with transfers}\label{recollection}

In this section we recall the definitions and basic properties of modulus sheaves with transfers from \cite{kmsyI} and \cite{shuji} (see also \cite{ksyII} for a more detailed summary).
\begin{enumerate}
    \item Denote by $\Sch$ the category of separated schemes of finite type over $k$ and by $\Sm$ the full subcategory of smooth schemes. For $X, Y \in \Sm$, an integral closed subscheme of $X\times Y$ that is finite and surjective over a connected component of $X$ is called a \emph{prime correspondence from $X$ to $Y$}. 
The category $\Cor$ of finite correspondences has the same objects as $\Sm$, and for $X, Y \in \Sm$, $\Cor(X,Y)$ is the free abelian group on the set of all prime correspondences from $X$ to $Y$ (see \cite{mvw}). We consider $\Sm$ as a subcategory of $\Cor$ by regarding a morphism in $\Sm$ as its graph in $\Cor$. 

Let $\PST=Fun(\Cor,\Ab)$ be the category of additive presheaves of abelian groups on $\Cor$ whose objects are called \emph{presheaves with transfers}. 
Let $\NST\subseteq \PST$ be the category of Nisnevich sheaves with transfers and let 
\eq{aVNis}{\aNis^V:\PST\to \NST}
be Voevodsky's Nisnevich sheafification functor, which is an exact left adjoint to the inclusion $\NST\to \PST$.
Let $\HI\subseteq \PST$ be the category of $\A^1$-invariant presheaves and put  $\HI_\Nis=\HI\cap \NST\subseteq \NST$. 
The product $\times$ on $\Sm$ yields a symmetric monoidal structure on $\Cor$,
which induces a symmetric monoidal structure on $\PST$ in the usual way.

\item We recall the definition of the category $\ulMCor$ from \cite[Definition 1.3.1]{kmsyI}. A pair $\sX = (X,D_X)$ of $X \in \Sch$ and an effective Cartier divisor $D$ on $X$ is called a \emph{modulus pair} if $X - |D_X| \in \Sm$.  
Let $\sX=(X,D_X)$, $\sY=(Y,D_Y)$ be modulus pairs and $\Gamma \in \Cor(X-D_X,Y-D_Y)$ be a prime correspondence. Let $\overline{\Gamma} \subseteq X \times Y$ be the closure of $\Gamma$, and let $\overline{\Gamma}^N\to X\times Y$ be the normalization. We say $\Gamma$ is \emph{admissible} (resp. \emph{left proper}) if $(D_X)_{\overline{\Gamma}^N}\geq (D_Y)_{\overline{\Gamma}^N}$
(resp. if $\overline{\Gamma}$ is proper over $X$). Let $\ulMCor(\sX, \sY)$ be the
subgroup of $\Cor(X-D_X,Y-D_Y)$ generated by all admissible left proper prime correspondences. The category $\ulMCor$ has modulus pairs as objects and $\ulMCor(\sX, \sY)$ as the group of morphisms from $\sX$ to $\sY$.

\item\label{ulMCorls}
Let $\ulMCorls\subset \ulMCor$ be the full subcategory of $(X,D)\in \ulMCor$ with $X\in \Sm$ and $|D|$ a simple normal crossing divisor on $X$. As observed in \cite[Remark 1.14]{shuji}, after assuming resolution of singularities, we can assume $\ulMCor\cong \ulMCorls$, as for every object $(X,D)\in \ulMCor$ there exists a proper birational map $p:X'\to X$ that is an isomorphism on $X-|D|$ and such that $|p^*D|$ is a simple normal crossing divisor. Hence the modulus correspondence $(X',D')\to (X,D)$ induced by the graph of $p$ is invertible in $\ulMCor$.

\item
There is a canonical pair of adjoint functors $\lambda\dashv \ulomega$:\[
\lambda: \Cor\to \ulMCor\quad X\mapsto (X,\emptyset),\]
\[
\ulomega: \ulMCor\to\Cor\quad (X,D)\mapsto X-|D|,
\]

\item
There is a full subcategory $\mathbf{MCor}\subset \ulMCor$ consisting of 
\emph{proper modulus pairs}, where
a modulus pair $(X,D)$ is \emph{proper} if $X$ is proper.
Let $\tau:\MCor\hookrightarrow \ulMCor$ be the inclusion functor and $\omega=\ulomega\tau$.

\item\label{sXn}
For all $n>0$ there is an endofunctor $(\_)^{(n)}$ on $\ulMCor$ preserving $\MCor$,
such that $(X,D)^{(n)}=(X,nD)$ where $nD$ is the $n$-th thickening of $D$.

\item  
We have two categories of \emph{modulus presheaves with trasnfers}:
\[\MPST=Fun(\MCor,\mathbf{Ab})\text{ and }\ulMPST=Fun(\ulMCor,\mathbf{Ab}).\]
Let $\Ztr(\sX)=\ulMCor(-,\sX) \in \ulMPST$ be the representable presheaf for
$\sX \in \MCor$. In this paper we frequently write $\sX$ for $\Ztr(\sX)$ for simplicity.

\item\label{ulomega} The adjunction $\lambda\dashv \ulomega$ induce a string of $4$ adjoint functors $(\lambda_!=\ulomega^!,\lambda^*=\ulomega_!,\lambda_*=\ulomega^*,\ulomega_*)$ (cf. \cite[Pr. 2.3.1]{kmsyI}):
\[\ulMPST\begin{smallmatrix}\ulomega^!\\\longleftarrow\\\ulomega_!\\\longrightarrow\\\ulomega^*\\\longleftarrow\\\ulomega_*\\ \longrightarrow\end{smallmatrix}\PST \]
where $\ulomega_!,\ulomega_*$ are localisations and $\ulomega^!$ and $\ulomega^*$ are fully faithful.

\item\label{omega} The functor $\omega$ yields a string of $3$ adjoint functors 
$(\omega_!,\omega^*,\omega_*)$ (cf. \cite[Pr. 2.2.1]{kmsyI}):
\[\MPST\begin{smallmatrix}
\omega_!\\\longrightarrow\\\omega^*\\\longleftarrow\\\omega_*\\\longrightarrow\end{smallmatrix}\PST \]
where $\omega_!,\omega_*$ are localisations and $\omega^*$ are fully faithful.

\item\label{tau}
The functor $\tau$ yields a string of $3$ adjoint functors $(\tau_!,\tau^*,\tau_*)$:
\[\MPST\begin{smallmatrix}\tau_!\\\longrightarrow\\\tau^*\\\longleftarrow\\\tau_*\\ \longrightarrow\end{smallmatrix}\ulMPST \]
where $\tau_!,\tau_*$ are fully faithful and $\tau^*$ is a localisation; $\tau_!$ 
has a pro-left adjoint $\tau^!$, hence is exact (cf. \cite[Pr. 2.4.1]{kmsyI}). We will denote by $\ulMPST^{\tau}$ the essential image of $\tau_!$ in $\ulMPST$. Moreover, we have (cf. \cite[Lem. 2.4.2]{kmsyI})
\eq{tauomegaulomega}{\omega_!=\ulomega_!\tau_! , \; \omega^*=\tau^*\ulomega^*,\;
\tau_!\omega^* =\ulomega^*.}

\item\label{ulMNST} For $F\in \ulMPST$ and $\sX=(X,D)\in \ulMCor$, write $F_{\sX}$ for the presheaf
on the small \'etale site $X_{\et}$ over $X$ given by $U\to F(\sX_U)$ for $U\to X$ \'etale, where $\sX_U=(U,D_{|U})\in \ulMCor$. We say $F$ is a Nisnevich sheaf if so is $F_{\sX}$ for all $\sX\in \ulMCor$ (see \cite[Section 3]{kmsyI}).
We write $\ulMNST\subset \ulMPST$ for the full subcategory of Nisnevich sheaves.
Let $\MNST\subset \MPST$ be the full subcategory of such objects $F$ that $\tau_! F\in \ulMNST$.
By \cite[Prop. 3.5.3]{kmsyI} and \cite[Theorem 2]{kmsyII}, the inclusion functors
\[\uli_\Nis: \ulMNST\to \ulMPST\qaq i_\Nis: \MNST \to \MPST\]
admit exact left adjoints $\ulaNis$ and $\aNis$ respectively and there are natural isomorphisms
\begin{equation}\label{tauaNis}
\tau_!\aNis \simeq \ulaNis \tau_!\qaq \aNis\tau^* \simeq \tau^*\ulaNis,
\end{equation}
and the adjunction from \eqref{tau} induces an adjunction
\[\MNST\begin{smallmatrix}\tau_!\\\longrightarrow\\\tau^*\\ \longleftarrow\end{smallmatrix}\ulMNST.\]
The functor $\ulaNis$ has the following description: For $F\in \ulMPST$ and $\sY\in \ulMCor$, let $F_{\sY,\Nis}$ be the usual Nisnevich sheafification of $F_{\sY}$. 
Then, for $(X,D)\in \ulMCor$ we have
\eq{ulaNisformular}{
\ulaNis F (X,D) = \colim_{f:Y\to X} F_{(Y,f^*D),\Nis}(Y),}
where the colimit is taken over all proper maps $f:Y\to X$ that induce isomorphisms $Y-|f^*D|\xrightarrow{\sim}X-|D|$.
\def\Sh{\mathrm{Sh}}

\item\label{b*leftexact}
For $X\in \Sch$, let $\Sh(X_\Nis,\Ab)$ be the abelian category of additive sheaves on $X_\Nis$. 
By definition of $\ulMNST$, we have an additive functor for $\sX=(X,D)\in \ulMCor$,
\[  \ulMNST \to  \Sh(X_\Nis,\Ab)\;;\; F \to F_\sX.\]
The functor is not exact in general but it is left exact by \eqref{ulaNisformular}.

\item\label{NST} By \cite[Pr. 6.2.1]{kmsyII}, the functors $\ulomega^*$ and $\ulomega_!$ respect $\ulMNST$ and $\NST$, and induce a pair of adjoint functors
\[\ulMNST\begin{smallmatrix}\\\ulomega_!\\\longrightarrow\\\ulomega^*\\\longleftarrow\end{smallmatrix}\NST ,\]
which are both exact. Moreover, we have
\[\ulomega_!\ulaNis=\aNis^V\ulomega_!\qaq
\ulaNis\ulomega^*=\ulomega^* \aNis^V.
\]


\item\label{semipure} We say that $F\in \ulMPST$ (resp. $\MPST$) is \emph{semi-pure} if the unit map
\[ u: F\to \ulomega^*\ulomega_!F \quad (\text{resp. } u: F\to \omega^*\omega_!F)\]
is injective. For $F\in \ulMPST$ (resp. $F\in \MPST$),  let $F^{sp}\in \ulMPST$ (resp. $F^{sp}\in \MPST$) be the image of
$F\to \ulomega^*\ulomega_! F$ (resp. $F\to \omega^*\omega_! F$) (called the semi-purification of $F$). One easily sees that the association $F\to F^{sp}$ gives a left adjoint to the inclusion  of the full subcategories of semipure objects into $\ulMPST$ and $\MPST$. For $F\in \MPST$  we have
\eq{tausp}{ \tau_!(F^{sp}) \simeq (\tau_! F)^{sp}.}
This follows from the fact that $\tau_!$ is exact and commutes with 
$\omega^*\omega_!$ and $\ulomega^*\ulomega_!$
since $\tau_!\omega^*=\ulomega^*$ and $\ulomega_!\tau_!=\tau_!$ (cf. \eqref{tau}). 
In particular $F\in \MPST$ is semiupre if and only if so is $\tau_!F\in \ulMPST$.
For $F\in \ulMPST$  we have
\eq{ulaNissp}{ \ulaNis(F^{sp}) \simeq (\ulaNis F)^{sp},}
where the $(\_)^{sp}$ on the right is defined for $F\in \ulMNST$ in the same way as above. This follows from the fact that $\ulaNis$ is exact and commutes with 
$\omega^*\omega_!$ and $\ulomega^*\ulomega_!$ (cf. \eqref{NST}).

\item\label{CI} 
Let $\bcube:=(\P^1,\infty)\in \MCor$.
We say $F \in \MPST$ is $\bcube$-invariant if $p^* : F(\sX) \to F(\sX \otimes \bcube)$ is an isomorphism for any $\sX \in \MCor$, where $p : \sX \otimes \bcube \to \sX$ is the projection.
Let $\CI$ be the full subcategory of $\MPST$ consisting of all $\bcube$-invariant objects. 

Recall from \cite[Theorem 2.1.8]{ksyII} that $\CI$ is a Serre subcategory of $\MPST$,
and that the inclusion functor $i^\bcube : \CI \to \MPST$
has a left adjoint $h_0^\bcube$ and a right adjoint $h^0_\bcube$ given 
for $F \in \MPST$ and $\sX \in \MCor$ by
\begin{align*} 
&h_0^\bcube(F)(\sX)
=\Coker(i_0^* - i_1^* : F(\sX \otimes \bcube) \to F(\sX)),
\\
\notag
&h^0_\bcube(F)(\sX)=\Hom(h_0^\bcube(\sX), F),
\end{align*}
where for $a\in k$ the section $i_a: \sX\to \sX\otimes \bcube$ is induced by the map $k[t]\to k[t]/(t-a)\cong k$.

For $\sX\in \MCor$, we write $h_0^\bcube(\sX)=h_0^\bcube(\Ztr(\sX))\in \CI$.

\item\label{CItspNis} 
Let $\CIt=\tau_!\CI\subset \ulMPST$ be the essential image of $\CI$ under $\tau_!$.
In this paper, for $F\in \CI$, we let $F$ denote also $\tau_!F\in \CIt$ by abuse of notation.
Let $\CIsp\subset\CI$ (resp. $\CItsp\subset\CIt$) be the full subcategory of semipure objects. By \eqref{tausp}, we have
\eq{FCItsp}{F^{sp}\in \CIt \qfor F\in \CIt,}
and $\tau_! $ and $\tau^*$ induce an equivalence of categories
\eq{CIspCItsp}{ \tau_! : \CIsp \simeq \CItsp: \tau^*}
with natural isomorphisms $\tau^*\tau_!\simeq id$ and $\tau_!\tau^*\simeq id$.

We also consider the full subcategories
\[\CIspNis =\CIsp\cap \MNST\subset \MNST,\]
\[\CItNis =\CIt\cap \ulMNST \subset \ulMNST.\]
\[\CItspNis =\CItsp\cap \ulMNST \subset \ulMNST.\]
By \cite[Th. 0.4]{shuji}, we have
\eq{ulaNisCItsp}{\ulaNis(\CItsp)\subset \CItspNis.}
By \cite[Th. 2 (1)]{kmsyII}, $\tau_! $ and $\tau^*$ induce an equivalence of categories
\eq{CIspCItspNis}{ \tau_! : \CIspNis \simeq \CItspNis: \tau^*}
with natural isomorphisms $\tau^*\tau_!\simeq id$ and $\tau_!\tau^*\simeq id$.


\item\label{RSC} 
We write  $\RSC\subseteq \PST$ for the essential image of $\CI$ under
$\omega_!$ (which is the same as the essential image of $\CItsp$ under
$\ulomega_!$ since $\omega_!=\ulomega_!\tau_!$ and $\ulomega_!F=\ulomega_!F^{sp}$).
Put $\RSC_{\Nis}=\RSC\cap \NST$.
The objects of $\RSC$ (resp. $\RSC_\Nis$) are called reciprocity presheaves 
(resp. sheaves). 
We have $\HI\subseteq \RSC$ and it contains also smooth commutative group schemes (which may have non-trivial unipotent part), and the sheaf $\Omega^i$ of K\"ahler differentials, and the de Rham-Witt sheaves $W\Omega^i$ (see \cite{ksyI} and \cite{ksyII}).

\item\label{omegaCI}
\def\hF{\hat{F}}
By \cite[Prop. 2.3.7]{ksyII} we have a pair of adjoint functors:
\begin{equation}\label{omegaCIadjoint0}
\CI\begin{smallmatrix}\omega_!\\\longrightarrow\\\omega^{\CI}\\\longleftarrow
\end{smallmatrix}\RSC, 
\end{equation}
where $\omega^{\CI}=h^0_\bcube\omega^*$ and it is fully faithful.
It induces a pair of adjoint functors:
\begin{equation}\label{omegaCIadjoint}
\CIt\begin{smallmatrix}\ulomega_!\\\longrightarrow\\\ulomega^{\CI}\\\longleftarrow
\end{smallmatrix}\RSC, 
\end{equation}
where $\ulomega^{\CI}=\tau_!h^0_\bcube\omega^*$ and it is fully faithful.
Indeed, let $F=\tau_!\hF$ for $\hF\in \CI$ and $G \in \RSC$. 
In view of \eqref{CI} and the exactness and full faithfulness of $\tau_!$, we have
\begin{multline*}
    \Hom_{\CIt}(F,\tau_!h^0_\bcube\omega^* G) \simeq  
\Hom_{\CI}(\hF,h^0_\bcube\omega^*G) \simeq \\
    \Hom_{\MPST}(\hF,\omega^*G) \simeq
\Hom_{\ulMPST}(\tau_! \hF,\ulomega^*G) \simeq \Hom_{\RSC}(\ulomega_! F,G).
\end{multline*}

By \cite[Theorem 2.4.1(2)]{ksyII}, \eqref{omegaCIadjoint} induce a pair of adjoint functors:
\begin{equation}\label{omegaCIadjointNis}
\CItsp_\Nis\begin{smallmatrix}\ulomega_!\\\longrightarrow\\\ulomega^{\CI}\\\longleftarrow \end{smallmatrix}\RSC_\Nis, 
\end{equation}
If $F\in \CIt$, the adjunction induces a canonical map\[
F\to \ulomega^{\CI}\ulomega_!F \]
which is injective if $F\in \CItsp$.

\item\label{tensor} 
$\ulMCor$ is equipped with a symmetric monoidal structure given by\[
    (X,D_X)\otimes (Y,D_Y):=(X\times Y,D_X\times Y + X\times D_Y),
\]
and $\mathbf{MCor}$ is clearly a $\otimes$-subcategory. Notice that the product is not a categorical product since the diagonal map is not admissible. It is admissible as a correspondence \[
    (X,D_X)^{(n)}\to (X,D_X)\otimes (X,D_X)\qquad \text{for }n\geq 2
\]
The symmetric monoidal structure $\otimes$ on $\ulMCor$ (resp. $\MCor$) induces a symmetric monoidal structure on $\ulMPST$ (resp. $\MPST$) in the usual way, and
$\tau_! $, $\omega_!$ and $\ulomega_!$ from \eqref{tau}, \eqref{omega} and \eqref{ulomega} are all monoidal  (see \cite[\S3]{rsy}).

\end{enumerate}
\bigskip

We end this section with some lemmas that will be needed in the rest of the paper.

The proof of the following Lemma is due to Kay R\"ulling. We thank him for letting us include it in our paper.
\begin{lemma}\label{lem;MapPSvsPST}
Let $p$ be the exponential characteristic of the base field $k$.
Let $F\in \PST$ such that 
\begin{enumerate}
\item\label{1} for all dominant \'etale maps $U\to X$ in $\Sm$ 
the pullback $F(X)\to F(U)$ is injective,
\item\label{2} $F$ has no $p$-torsion.
\end{enumerate}
Then, for any $G\in \PST$, the natural map
\[ \uHom_{\PST}(G,F) \to \uHom_{\PS}(G,F)\]
is an isomorphism.
\end{lemma}
\begin{proof}(Kay R\"ulling)
First we prove $\Hom_{\PST}(G,F)=\Hom_{\PS}(G,F)$, i.e. any morphism $\varphi: G\to F$ of presheaves on $\Sm$ is also a morphism in $\PST$.
We have to show $\varphi(f^*a)=f^*\varphi(a)$ in $F(X)$, 
for $a\in G(Y)$ and $f\in  \Cor(X,Y)$ a prime correspondence.
By (1) we can reduce to the case $X=\Spec K$, with $K$ a function field over $k$.
In this case we can write $f^*= h_*g^*$, where  $h: \Spec L\to \Spec K$ is induced by a finite field extension $L/K$
and $g: \Spec L\to Y$ is a morphism. 
Since $\varphi$ is a morphism of presheaves on $\Sm$, 
we are reduced to show 
\begin{equation}\label{(*)}
	h_*\varphi(a)=\varphi(h_*a),\quad a\in G(L).
	\end{equation}
It suffices to consider the following two cases:

{\em 1st case: L/K is finite separable.}
Let $E/K$ be a finite Galois extension containing $L/K$ and denote by $j:\Spec E\to \Spec K$ the induced morphism and 
by $\sigma_i : \Spec E\to \Spec L$ the morphism induced by all $K$-embeddings of $L$ into $E$.
Since $G\in \PST$ we obtain in $G(E)$
\[j^*h_*a= (h^t\circ j)^*a=\sum_i \sigma_i^*(a).\]
Thus 
\[j^*\varphi(h_*a)= \varphi(j^*h_*a)=\varphi(\sum_i \sigma_i^*(a))=\sum_i\sigma_i^*\varphi(a)= j^*h_*\varphi(a).\]
Since $j^*: F(L)\to F(E)$ is injective by (1) this shows \eqref{(*)} in this case.

{\em 2nd case: L/K is purely inseparable of degree $p$.}
In this case we have $h^*h_*=(h^t\circ h): G(L)\to G(L)$ is multiplication by $p$ as well as
$h_* h^*: G(K)\to G(K)$.
Thus 
\[h^*\varphi(h_*a)= \varphi(h^*h_*a)= p\varphi(a)= h^*h_*\varphi(a);\]
applying $h_*$ yields
\[p\varphi(h_*a)=p h_*\varphi(a);\]
thus \eqref{(*)} follows from (2).
\medskip

Next we prove the analogous statement for internal hom's. 
Indeed, note that for $X\in \Sm$, $\uHom_{\PST}(\Z_{\tr}(X), F)\in \PST$ also satisfies (1) and (2) above and that we have  
\begin{equation}\label{(**)}\uHom_{\PST}(\Z_{\tr}(X), F)= F(X\times -)= \uHom_{\PS}(h_X, F) \quad \text{in }\PS,
	\end{equation}
where $h_X= \Z(\Hom_{\Sm}(-, X))$.
Thus for $G\in \PST$ 
\begin{align*}
\uHom_{\PST}(G, F)(X)& = \Hom_{\PST}(\Z_\tr(X), \uHom_{\PST}(G,F))\\
                                  & = \Hom_{\PST}(G\otimes^{\PST} \Z_\tr(X), F)\\
                                  & = \Hom_{\PST}(G, \uHom_{\PST}(\Z_\tr(X), F))\\
                                  &= \Hom_{\PS}(G, \uHom_{\PST}(\Z_\tr(X), F)), & & \text{by } \eqref{(*)}\\
                                  &=\Hom_{\PS}(G, \uHom_{\PS}(h_X, F)), & & \text{by } \eqref{(**)}\\
                                  &= \Hom_{\PS}(G\otimes^{\PS} h_X, F)\\
                                  &=\Hom_{\PS}(h_X, \uHom_{\PS}(G,F))\\
                                  &=\uHom_{\PS}(G,F)(X).
\end{align*}
This completes the proof of Lemma \ref{lem;MapPSvsPST}.
\end{proof}

\begin{lemma}\label{lem0;omegauHom}
For $F\in \PST$ and $X\in \Sm$, we have a natural isomorphism
\[ \ulomega^*\uHom_{\PST}(\Ztr(X), F) \simeq 
\uHom_{\ulMPST}(\Ztr(X,\emptyset),\ulomega^* F).\]
\end{lemma}
\begin{proof}
For $\sY=(Y,E)\in\ulMCor$ with $V=Y-|E|$, we have natural isomorphisms
\begin{multline*}
\ulomega^*\uHom_{\PST}(\Ztr(X), F)(\sY) \simeq\uHom_{\PST}(\Ztr(X), F)(V)
\simeq  \Hom_{\PST}(X\times V,F) \\
\simeq \Hom_{\ulMPST}((X,\emptyset)\otimes \sY,\ulomega^* F)
\simeq \uHom_{\ulMPST}(\Ztr(X,\emptyset),\ulomega^* F)(\sY).
\end{multline*}
This proves the lemma.
\end{proof}

\begin{lemma}\label{lem;omegauHom}
For $F\in \ulMPST$ and $X\in \Sm$, we have a natural isomorphism
\[ \ulomega_!\uHom_{\ulMPST}(\Ztr(X,\emptyset), F) \simeq 
\uHom_{\PST}(\Ztr(X),\ulomega_! F).\]
\end{lemma}
\begin{proof}
For $Y\in \Sm$, we have natural isomorphisms
\begin{multline*}
\ulomega_!\uHom_{\ulMPST}(\Ztr(X,\emptyset), F)(Y) \simeq
\uHom_{\ulMPST}(\Ztr(X,\emptyset), F)(Y,\emptyset) \\
\simeq \Hom_{\ulMPST}(\Ztr(X\times Y,\emptyset), F) 
\simeq \Hom_{\PST}(X\times Y,\ulomega_! F)\\
\simeq\uHom_{\PST}(\Ztr(X),\ulomega_! F)(Y).
\end{multline*}
This proves the lemma.
\end{proof}

\begin{lemma}\label{lem;RSCexactness}
A complex in $C^\bullet$ in $\NST$ such that $C^n\in \RSC$ for all $n\in \Z$ is exact if and only if $C^\bullet(K)$ is exact as a complex of abelian groups for any function field $K$ .
\end{lemma}
\begin{proof}
The cohomology sheaves $H^n_{\Nis}(C^\bullet)$ are in $\RSC_\Nis$ by \cite[Th.0.1]{shuji}.
Hence for all $X\in \Sm$, by \cite[Th. 0.2]{shuji} there is an injective map $(H^n_{\Nis}C^\bullet)(X)\hookrightarrow (H^n_{\Nis}C^\bullet)(k(X))$. Hence the lemma follows from the fact that
$(H^n_{\Nis}C^\bullet)(k(X))\cong H^n(C^\bullet(k(X)))$, since $k(X)$ is henselian local.
\end{proof}

\begin{lemma}\label{lem;homrec}
For $G\in \RSC$ and $F\in \PST$ such that $F$ is a quotient of a representable sheaf, $\uHom_{\PST}(F,G)\in \RSC$.
\end{lemma}
\begin{proof} 
First assume $F=\Ztr(X)$ with $X\in \Sm$. 
Put $\tG=\ulomega^{\CI}G\in \CIt$ (cf. \eqref{omegaCI}). 
The adjunction \eqref{omegaCIadjoint} implies $\ulomega_!\tG\simeq G$.
Lemma \ref{lem;omegauHom} implies a natural isomorphism
\[ \uHom_{\PST}(\Ztr(X),G)\simeq \ulomega_! \uHom_{\ulMPST}(\Ztr(X,\emptyset),\tG) .\]
Thus it suffices to show 
\[ \uHom_{\ulMPST}(\Ztr(X,\emptyset),\tG)  \in \CIt.\]
The $\bcube$-invariance follows directly from the one for $\tG$. 
The fact that it is in $\MPST^\tau$ follows from \cite[Lemma 1.27(2)]{shuji}.

Now assume there is a surjection $\Ztr(X)\to F$ in $\PST$,
where $X\in \Sm$. It induces an injection
\[\uHom_{\PST}(F,G)\hookrightarrow \uHom_{\PST}(\Ztr(X),G).
\]
Since $\uHom_{\PST}(\Ztr(X),G)\in \RSC$ as shown above and
$\RSC\subset \PST$ is closed under finite products and subobjects,
 we get $\uHom_{\PST}(F,G)\in \RSC$ as desired. This completes the proof.
\end{proof}

\begin{lemma}\label{P1cohomolgyvanishing}
Let $F\in\ulMNST$ be such that $F^{sp}\in \CItNis$ (cf. \eqref{CItspNis}).
For any function field $K$ over $k$, we have
\[ H^i(\P^1_K,F_{(\P^1_K,0+\infty)})=0\qfor i>0.\]
\end{lemma}
\begin{proof}
If $F$ is semi-pure, the assertion follows from \cite[Th. 9.1]{shuji}.
In general we use the exact sequence in $\ulMNST$:
\[ 0\to C \to F \to F^{sp}\to 0\]
to reduce to the above case noting $H^i(\P^1_K,C_{(\P^1_K,0+\infty)})=0$ for $i>0$ since
$C_{(\P^1_K,0+\infty)}$ is supported on $\{0,\infty\}$.
\end{proof}

\begin{lemma}\label{lem;CIsheafication}
For $F\in \CIt$ and a function field $K$ over $k$, we have
\[ \ulaNis F(K) \isom \ulaNis F(\bcube\otimes K) .\]
\end{lemma}
\begin{proof}
We consider the exact sequence in $\ulMPST$:
\[ 0\to C\to F \to F^{sp}\to 0\qwith \;\; \ulomega_! C=0.\]
Since $\ulaNis$ is exact, from this we get an exact sequence in $\ulMNST$:
\[ 0\to \ulaNis C\to \ulaNis F \to \ulaNis F^{sp}\to 0.\]
Since $C_{(\P^1_K,0+\infty)}$ is supported on $\{0_K,\infty_K\}$,
we have by \eqref{ulaNisformular}
\[(\ulaNis C)_{(\P^1_K,0+\infty)} = C_{(\P^1_K,0+\infty)}.\]
Hence the exact sequence gives rise to a commutative diagram
\[\xymatrix{
0\ar[r] & C(K) \ar[r]\ar[d]^{\simeq} & F(K) \ar[r]\ar[d] & F^{sp}(K) \ar[r]\ar[d]^{\simeq} & 0\\
0\ar[r] & C(\bcube\otimes K) \ar[r] & \ulaNis F(\bcube\otimes K) \ar[r] & \ulaNis F^{sp}(\bcube\otimes K) \\}\]
The left (resp. right ) vertical map is an isomorphism since $C\in \CIt$ (resp. 
thanks to \cite[Th. 10.1]{shuji}).
This completes the proof.
\end{proof}
\medbreak

\def\lamgm{\lambda_{\G_m}}
\def\lamga{\lambda_{\G_a}}
\def\cubegmm{\bcube_{\G_m}}
\def\cubega{\bcube_{\G_a}}

Let $\A^1_t=\Spec k[t]$ be the affine line with the coordinate $t$.
Consider the map in $\PST$:
\[ \lamgm: \Ztr(\A^1_t-\{0\}) \to \G_m\]
given by $t\in \G_m(\A^1_t-\{0\}) =k[t,t^{-1}]^{\times}$, and the map in $\PST$:
\[ \lamga: \Ztr(\A^1_t) \to \G_a\]
given by $t\in \G_a(\A^1_t) =k[t]$.
Note that $\lamgm$ and $\lamga$ factor through
\[ \Coker(\Z\rmapo{i_1}  \Ztr(\A^1_t-\{0\})) \qaq
\Coker(\Z\rmapo{i_0}  \Ztr(\A^1_t)),\]
with $i_1$ and $i_0$ induced by the points $1\in \A^1_t-\{0\}$ and $0\in \A^1_t$ respectively.

\begin{lemma}\label{lem;GmGa}
\begin{itemize}
\item[(1)]
The composite map
\[ \omega_!\Ztr(\P^1,0+\infty)\simeq \Ztr(\A^1_t-\{0\}) \rmapo{\lamgm} \G_m\]
induces an isomorphism
\begin{equation}\label{Gmh0cube}
\aVNis \omega_! h_0^\bcube(\cubegmm) \isom \G_m,
\end{equation}
where
$\cubegmm=\Coker(\Z\rmapo{i_1}  \Ztr(\P^1,0+\infty))\in \MPST $.
\item[(2)]
The composite map
\[ \omega_!\Ztr(\P^1,2\infty)\simeq \Ztr(\A^1_t) \rmapo{\lamga} \G_a\]
induces an isomorphism
\begin{equation}\label{Gah0cube}
\aVNis \omega_! h_0^\bcube(\cubega) \isom \G_a,
\end{equation}
where
$\cubega=\Coker(\Z\rmapo{i_0}  \Ztr(\P^1,2\infty)) \in \MPST$.
\end{itemize}
\end{lemma}
\begin{proof}
We prove only (2). The proof of (1) is similar.
By \cite[Lem. 1.36 and Th. 0.1]{shuji}, we have
$\aVNis \omega_! h_0^\bcube(\cubega)\in \RSC_\Nis$. 
Hence, by Lemma \ref{lem;RSCexactness}, it suffices to show that 
the map $\Ztr(\A^1)(K)\rmapo{\lamgm} \G_a(K)=K$ for a function field $K$ over $k$, induces an isomorphism
$\omega_! h_0^\bcube(\cubega)(K) \simeq K$.
We know that $\Ztr(\A^1_t)(K)$ is identified with the group of $0$-cycles on $\A^1_K=\A^1\otimes_k K$. Then, by \cite[Th. 3.2.1]{ksyII}, the kernel of 
$\Ztr(\A^1)(K)\to \omega_! h_0^\bcube(\cubega)(K)$ is generated by 
the class of $0\in \A^1_K$ and 
$\div_{\A^1_K}(f)$ for $f\in K(t)^\times$ such that 
$f\in 1+ \fm_\infty^2 \sO_{\P^1_K,\infty}$, where $\fm_\infty$ is the maximal ideal of the local ring $\sO_{\P^1_K,\infty}$ of $\P^1_K$ at $\infty$.
Now (2) follows by an elementary computation.
\end{proof}

\begin{lemma}\label{lem;uHomNis}
We have
\[\uHom_{\ulMPST}(G,F)\in \ulMNST\qfor 
G\in \ulMPST,\; F\in \ulMNST.\]
\end{lemma}
\begin{proof}
Put $H=\uHom_{\ulMPST}(G,F)$.
Let $\sX\in \ulMCor$ and 
\[\xymatrix{
\sW \ar[r]\ar[d] &\sV \ar[d]\\
\sU \ar[r] & \sX\\}\]
be a $\underline{\MV}^{\fin}$-square as defined in \cite[Def. 3.2.1]{kmsyI}.
By \cite[Def. 4.5.2 and Lem. 4.2.3]{kmsyI}, it suffices to show the exactness of
\[ 0\to H(\sX) \to H(\sU)\oplus H(\sV) \to H(\sW).\]
By the adjunction, we have
\[ H(\sX) = \Hom_{\ulMPST}(G,F^\sX)\qwith \;F^\sX=\uHom_{\ulMPST}(\Ztr(\sX),F)).\]
Hence it suffices to show the exactness of the following sequence in $\ulMPST$:
\[ 0\to F^\sX \to F^\sU \oplus F^\sV  \to F^\sW.\]
Taking $\sY\in \ulMCor$, this is reduced to showing the exactness of 
 \[ 0\to F(\sX\otimes\sY) \to F(\sU\otimes\sY)\oplus F(\sV\otimes\sY) \to
F(\sW\otimes\sY).\]
This follows from the fact that $\underline{\MV}^{\fin}$-squares are preserved by the product $\otimes$ in $\ulMCor$.
\end{proof}

\begin{prop}\label{prop;aNisasoociative} 
\begin{itemize}
\item[(i)]
For $F,G\in \ulMPST$, we have a natural isomorphism 
\[\ulaNis(F\otuMPST G) \simeq \ulaNis(\ulaNis F\otuMPST \ulaNis G)\]
induced by the natural maps $F\to \ulaNis F$ and $G\to \ulaNis G$.
\item[(ii)]
For $F,G\in \MPST$, we have a natural isomorphism
\[\aNis(F\otMPST G) \simeq \aNis(\aNis F\otMPST \aNis G)\]
induced by the natural maps $F\to \aNis F$ and $G\to \aNis G$.
\end{itemize}
\end{prop}
\begin{proof}
For $H\in \ulMNST$, we have isomorphisms
\[\begin{aligned}
\Hom_{\ulMNST}(\ulaNis(F\otuMPST G),H) 
&\simeq \Hom_{\ulMPST}(F\otuMPST G,H) \\
&\simeq \Hom_{\ulMPST}(F,\uHom_{\ulMPST}(G,H)) \\
&\overset{(*1)}{\simeq} \Hom_{\ulMPST}(\ulaNis F,\uHom_{\ulMPST}(G,H)) \\
&\simeq \Hom_{\ulMPST}(\ulaNis F\otuMPST G,H) \\
&\simeq \Hom_{\ulMPST}(G,\uHom_{\ulMPST}(\ulaNis F,H)) \\
&\overset{(*2)}{\simeq} \Hom_{\ulMPST}(\ulaNis G,\uHom_{\ulMPST}(\ulaNis F,H)) \\
&\simeq \Hom_{\ulMPST}(\ulaNis F\otuMPST \ulaNis G,H) \\
&\simeq \Hom_{\ulMNST}(\ulaNis(\ulaNis F\otuMPST \ulaNis G),H) \\
\end{aligned}\]
where $(*1)$ and $(*2)$ follow from the fact 
$\uHom_{\ulMPST}(A,H)\in \ulMNST$ for $A\in \ulMPST$ by Lemma \ref{lem;uHomNis}.
This proves (i).

For $F,G\in \MPST$, we have isomorphisms
\[\begin{aligned}
\tau_! \aNis(F\otMPST G) 
&\overset{(*1)}{\simeq} \ulaNis \tau_! (F\otMPST G)  \\
&\overset{(*2)}{\simeq} \ulaNis (\tau_! F\otuMPST \tau_!G)  \\
&\overset{(*3)}{\simeq} \ulaNis (\ulaNis\tau_! F\otuMPST \ulaNis\tau_!G)  \\
&\overset{(*4)}{\simeq} \ulaNis (\tau_!\aNis F\otuMPST \tau_!\aNis G)  \\
&\overset{(*5)}{\simeq} \ulaNis \tau_!(\aNis F\otuMPST \aNis G)  \\
&\overset{(*6)}{\simeq} \tau_!\aNis(\aNis F\otuMPST \aNis G)  \\
\end{aligned}\]
where $(*1)$, $(*4)$ and $(*6)$ follow from \eqref{tauaNis}, $(*2)$ and $(*5)$
follow from the monoidality of $\tau_!$ (\cite[\S3.8]{rsy}) and $(*3)$ follows from (i). Since $\tau_!$ is fully faithful, this implies (ii).
This completes the proof of the lemma. 
\end{proof}

\begin{lemma}\label{lem1;MPSTsp}
There are natural isomorphisms for $F,G\in \MPST$ 
\eq{MPSTsp1}{ (F\otMPST G)^{sp} \simeq (F^{sp}\otMPST G)^{sp}\simeq 
(F^{sp}\otMPST G^{sp})^{sp}.} 
\end{lemma}
\begin{proof}
We have an exact sequence in $\MPST$:
\[ 0\to C\to F \to F^{sp}\to 0 \qwith\;\; \omega_! C =0.\]
Since $(-)\otMPST G:\MPST\to \MPST$ is right exact, we get an exact sequence
\[  C\otMPST G \to F\otMPST G  \to F^{sp}\otMPST G \to 0. \]
We have $\omega_!(C\otMPST G )=0$ since $\omega_! : \MPST \to \PST$ 
is monoidal by \cite[\S3.6]{rsy}. Hence we get an isomorphism
$(F\otMPST G)^{sp}  \simeq (F^{sp}\otMPST G)^{sp}$.
This implies \eqref{MPSTsp1}.
\end{proof}

\begin{lemma}\label{lem1;hcubesp}
There are natural isomorphisms for $F,G,H\in \MPST$
\eq{h0cube2}{h_0^\bcube(F^{sp})^{sp}\simeq h_0^\bcube(F)^{sp},} 
\eq{h0cube3}{\hcube(F\otMPST G) \simeq h_0^\bcube(h_0^\bcube(F)\otMPST h_0^\bcube(G)).}
\end{lemma}
\begin{proof}
We have an exact sequence in $\MPST$:
\[ 0\to C\to F \to F^{sp}\to 0 \qwith\;\; \omega_! C =0.\]
From this we get an exact sequence in $\MPST$:
\[  \hcube(C) \to \hcube(F)  \to \hcube(F^{sp})\to 0 \]
since $\hcube:\MPST\to \MPST$ is right exact.
We have $\omega_!\hcube(C)=0$ since $\omega_! : \MPST \to \PST$ 
is exact and $\hcube(C)$ is a quotient of $C$. Hence we get an isomorphism
$\omega_!\hcube(F)  \simeq \omega_!\hcube(F^{sp})$.
This implies \eqref{h0cube2}.

For $H\in \CI$, we have isomorphisms
\[\begin{aligned}
\Hom_{\CI}(\hcube(F\otMPST G),H)
&\simeq \Hom_{\MPST}(F\otMPST G,H)\\
&\simeq \Hom_{\MPST}(F,\uHom_{\MPST}(G,H))\\
&\overset{(*)}{\simeq} \Hom_{\MPST}(\hcube(F),\uHom_{\MPST}(G,H))\\
&\simeq \Hom_{\MPST}(\hcube(F)\otMPST G,H))\\
&\simeq \Hom_{\CI}(\hcube(\hcube(F)\otMPST G),H))\\
\end{aligned}\]
where $(*)$ follows from the fact that $\uHom_{\MPST}(G,H)\in \CI$ for $H\in \CI$,
which follows easily from the definition. This shows  
\[\hcube(F\otMPST G)\simeq \hcube(\hcube(F)\otMPST G),\]
which implies \eqref{h0cube3}.
\end{proof}
\medbreak

From \eqref{ulaNisCItsp}, we have $\ulaNis(\CItsp)\subset \CItspNis$, which implies 
\[\aNis(\CIsp)\subset \CIspNis.\]
Indeed, for $F\in \CIsp$, we have $\tau_! \aNis F\simeq \ulaNis \tau_! F\in \CItspNis$
by \eqref{tauaNis}, which implies $\aNis F\in \CIspNis$ by definition (cf. \eqref{ulMNST} and \cite[Def. 3]{kmsyII}).
Thus we get an induced functor
\eq{aCI}{\aNisCI: \CIsp \to \CIspNis.}
By definition we have 
\eq{aNisaNisCI}{ \aNisCI(F) = \aNis j(F)\qfor F\in \CIsp,}
where $j: \CIsp \to \MPST$ is the inclusion.

\begin{lemma}\label{lem;aNsiCI}
$\aNisCI$ is a left adjoint to the inclusion $\CIspNis\to \CIsp$.
\end{lemma}
\begin{proof}
This follows easily from the fact that $\aNis$ is a left adjoint to the inclusion $\MNST\to \MPST$ and the inclusions $\CIsp\to \MPST$ and 
$\CIspNis\to \MNST$ are fully faithful. 
\end{proof}

\begin{lemma}\label{lem2;hcubesp}
Consider the functors
\[ \hcubesp : \MPST\to \CIsp\;:\; F \to \hcube(F)^{sp},\]
\[  \hcubespNis: \MPST\to \CIspNis\;:\; F\to \aNisCI\hcubesp(F).\]
\begin{itemize}
\item[(i)]
The functor $\hcubesp$ (resp. $\hcubespNis$) is a left adjoint to the inclusion
$\CIsp \to \MPST$ (resp. $\CIspNis \to \MPST$).
For $F\in \MPST$, we have natural isomorphisms
\[ \hcubesp(F) \simeq \hcubesp\hcubesp(F)\qaq
\hcubespNis(F) \simeq \hcubespNis\hcubespNis(F).\]
\item[(ii)]
For $F\in \MPST$, the natural map $F\to \aNis F$ induces isomorphisms 
\[\hcubespNis(F) \simeq  \hcubespNis( \aNis F).\]
\item[(iii)]
For $F\in \MPST$, we have natural isomorphisms
\[ \hcubesp(F\otMPST G) \simeq \hcubesp(\hcubesp(F)\otMPST \hcubesp(G)) ,\]
\[ \hcubespNis(F\otMPST G) \simeq \hcubespNis(\hcubespNis(F)\otMPST \hcubespNis(G)) .\]
\end{itemize}
\end{lemma}
\begin{proof}
The first statement of (i) follows from the left-adjointness of $\hcube$, $(-)^{sp}$ and $\aNis$.
The second statement of (i) is a formal consequence of the first since the inclusions are fully faithful.

To show (ii), consider the commutative diagram 
\[\xymatrix{
\CIspNis \ar[r]^{i_{\CI}} \ar[d]^{j_\Nis} &\CIsp\ar[d]^{j} \\
\MNST \ar[r]^{i}  & \MPST
}\]
where the functors are inclusions.
For $F\in \MPST$ and $G\in \CIspNis$, we have isomorphisms
\[\begin{aligned}
\Hom_{\CIspNis}(\hcubespNis i \aNis F, G)
&\overset{(*1)}{\simeq} \Hom_{\CIsp}(\hcubesp i \aNis F, i_{\CI}G)\\
&\overset{(*2)}{\simeq} \Hom_{\MPST}(i\aNis F, j i_{\CI}G)\\
&\simeq \Hom_{\MPST}(i\aNis F, i j_{\Nis}G) \\
&\overset{(*3)}{\simeq}\Hom_{\MNST}(\aNis F, j_{\Nis}G) \\
&\simeq \Hom_{\MPST}(F, i j_{\Nis}G) \\
&\simeq \Hom_{\MPST}(F, j i_{\CI}G) \\
&\overset{(*4)}{\simeq} \Hom_{\MPST}( \hcubesp F, i_{\CI}G) \\
&\overset{(*5)}{\simeq} \Hom_{\MPST}(\aNisCI \hcubesp F, G) \\
\end{aligned}\]
where $(*1)$ and $(*5)$ (resp. $(*2)$ and $(*4)$, resp. $(*3)$) follow from Lemma \ref{lem;aNsiCI} (resp. (i), resp. the full faithfulness of $i$). This proves (ii).

For $F,G\in \MPST$, we have natural isomorphisms
\[\begin{aligned}
\hcubesp(F\otMPST G) 
&\overset{\eqref{h0cube2}}{\simeq} \hcube((F\otMPST G)^{sp})^{sp} \\ 
&\overset{\eqref{MPSTsp1}}{\simeq} \hcube((F^{sp}\otMPST G^{sp})^{sp})^{sp} \\
&\overset{\eqref{h0cube2}}{\simeq} \hcube(F^{sp}\otMPST G^{sp}))^{sp} \\
&\overset{\eqref{h0cube3}}{\simeq} \hcube(\hcube(F^{sp})\otMPST \hcube(G^{sp}))^{sp} \\
&\overset{\eqref{h0cube2}}{\simeq} \hcube((\hcube(F^{sp})\otMPST \hcube(G^{sp}))^{sp})^{sp} \\
&\overset{\eqref{MPSTsp1}}{\simeq} \hcube((\hcube(F^{sp})^{sp}\otMPST \hcube(G^{sp})^{sp})^{sp})^{sp}\\ 
&\overset{\eqref{h0cube2}}{\simeq} \hcube((\hcubesp(F)\otMPST \hcubesp(G))^{sp})^{sp} \\
&\overset{\eqref{h0cube2}}{\simeq} \hcube(\hcubesp(F)\otMPST \hcubesp(G))^{sp} \\
& = \hcubesp(\hcubesp(F)\otMPST \hcubesp(G))\\
\end{aligned}\]
This proves the first isomorphism of (iii). From this we get natural isomorphisms
\[\begin{aligned}
\hcubespNis(F\otMPST G) 
& \simeq  \hcubespNis(\hcubesp(F)\otMPST \hcubesp(G)) \\
&\overset{(*1)}{\simeq}  \hcubespNis \aNis (\hcubesp(F)\otMPST \hcubesp(G)) \\  
&\overset{(*2)}{\simeq}  \hcubespNis \aNis(\hcubespNis(F)\otMPST \hcubespNis(G)) \\
&\overset{(*3)}{\simeq}  \hcubespNis(\hcubespNis(F)\otMPST \hcubespNis(G)) \\
\\
\end{aligned}\]
where $(*1)$ and $(*3)$ follow from (ii) and $(*2)$ follows from Proposition \ref{prop;aNisasoociative} in view of \eqref{aNisaNisCI}.
This completes the proof of the lemma.
\end{proof}

\section{Some lemmas on contractions}\label{sec2}

For an integer $a\geq 1$ put $\cubegma=(\P^1,a(0+\infty))\in \MCor$ and
\[ \cubegmreda= \Ker\big(\Ztr(\cubegma) \to \Z=\Ztr(\Spec k,\emptyset)\big).\]
The inclusion $\A^1-\{0\}\hookrightarrow \A^1$ induces a map $\bcube^{(a)}\to \bcube$ in $\MCor$ for all $a$.
Note that the composite map 
\begin{equation}\label{eq;cuberedgm}
\cubegmred \hookrightarrow \cubegm \to \cubegmm
\end{equation}
is an isomorphism, where $\cubegmm$ is from \eqref{Gmh0cube}.

For $F\in \ulMPST$, we write
\[ \gamma F=\Coker\big(\uHom_{\ulMPST}(\bcube,F)\to
\uHom_{\ulMPST}(\cubegm,F)\big)\in \ulMPST,\]
where the map is induced by $\cubegm\to \bcube$ in $\MCor$. 
If $F\in \CIt$, the projection $\bcube\to \Spec k$ induces an isomorphism
\[ F=\uHom_{\ulMPST}(\Spec k,F)\simeq \uHom_{\ulMPST}(\bcube,F).\]
Thus we get an isomorphism
\begin{equation}\label{eq;gammaFCI}
\gamma F\simeq \uHom_{\ulMPST}(\cubegmred,F)\overset{(*)}=\uHom_{\ulMPST}(h_0^{\bcube}(\cubegmred),F) \qfor F\in \CIt,
\end{equation}
where the equality $(*)$ follows from the adjunction from \eqref{CI}. Note $\gamma F\in \CItsp$ for $F\in \CItsp$.
We also define
\[ \gamma_\Nis F =\ulaNis \gamma F\in \ulMNST.\]
By \eqref{eq;gammaFCI} and Lemma \ref{lem;uHomNis}, we have  
\[\gamma_\Nis F = \gamma F\qfor F\in \CIt_\Nis,\]
We write for an integer $n\geq 1$ (cf, \S1\eqref{tensor})
\begin{equation}\label{eq;gammaFCIn2}
\gamma^n F\cong \uHom_{\ulMPST}((\cubegmred)^{\otimes_{\ulMPST}n},F) 
\cong \overbrace{\gamma\gamma\cdots \gamma}^{n\;\text{times}} F.
\end{equation}
Notice that it is  
\medbreak
 
The proof of the following Lemma is due to Kay R\"ulling. We thank him for letting us include it in our paper.

\begin{lemma}\label{lem;Kay}
The unit map
\begin{equation}\label{eq:Gm1}
\ulaNis h_0^{\bcube}(\cubegm)^{sp}\xrightarrow{\simeq} 
\ulomega^*\ulomega_! \ulaNis h_0^{\bcube}(\cubegm)
\cong \ulomega^*(\G_m\oplus \Z)
\end{equation}
is an isomorphism, where the second isomorphism in \eqref{eq:Gm1} holds by 
Lemma \ref{lem;GmGa} and \eqref{eq;cuberedgm}.
\end{lemma}
\begin{proof}\emph{(Kay R\"ulling)}
The unit map is injective by semipurity.
It remains to show the surjectivity. 
By definition of the sheafification functor, it suffices to show the surjectivity on $(\Spec R, (f))$, where $R$ is an integral local $k$-algebra and $f\in R\setminus\{0\}$, such that $R_f$ is regular.
Denote by
\[\psi: \Z_{tr}(\P^1,0+\infty) (R,f)\to R_f^\times\oplus\Z\]
the precomposition of $\eqref{eq:Gm1}$ evaluated at $(R,f)$ with the 
quotient map $\Z_{tr}(\P^1,0+\infty)(R,f)\to \ulaNis h_0^{\bcube}(\cubegm)^{sp}$.

We show that $\psi$ is surjective. To this end, observe that for $a\in R_f^\times$ we find $N\ge 0$ and  $b\in R$ such that
\begin{equation}\label{lem:Gm4}
    ab=f^{N}, \quad \text{and}\quad af^{N}\in R.
\end{equation}
Set $W:=V(t^{N}-a)\subset \Spec R_f[t, 1/t]$ and $K:=\text{Frac}(R)$.

The map $\Cor(K,\A^1-\{0\})\to \Pic(\P^1_K,0+\infty)\cong K^\times\oplus \Z$ which induces the second isomorphism of \eqref{eq:Gm1} sends a prime correspondence $V(a_0+a_1t+\ldots a_rt^r)$ to $((-1)^ra_0/a_r,r)$, hence we have:
\begin{equation}\label{lem:Gm3.5}
\psi(V(a_0+ a_1t+\ldots a_r t^r))= ((-1)^r a_0/a_r, r)
\end{equation}
provided that $V(a_0+ a_1t+\ldots a_r t^r)\in \ulMCor((R,f),(\P^1,0+\infty))$.

For any $a\in R_f^\times$, consider $h=t^{N}-a$ and let $h=\prod_i h_i$ be the decomposition into monic irreducible factors in $K[t,1/t]$
and denote by $W_i\subset \Spec R_f[t, 1/t]$ the closure of $V(h_i)$. (Note that $W_i=W_j$ for $i\neq j$ is allowed.)

The $W_i$ correspond to the components of $W$ which are dominant over $R_f$; since $W$ is finite and surjective over $R_f$, so are the $W_i$. We claim
\begin{equation}\label{lem:Gm5}
W_i\in \ulMCor((R,f),(\P^1,0+\infty))
\end{equation}
Indeed, let $I_i$ (resp. $J_i$) be the ideal of the closure of $W_i$ in $\Spec R[t]$ (resp. $\Spec R[z]$ with $z=1/t$).
By \eqref{lem:Gm4} 
\[bt^N- f^N\in I_i \quad \text{and} \quad f^N- f^Na z^N\in J_i.\]
Hence $(f/t)^N\in R[t]/I_i$ and $(f/z)^N\in R[z]/J_i$.
It follows that $f/t$  (resp. $f/z$) is integral over $R[t]/I_i$ (resp. $R[z]/J_i$);
thus \eqref{lem:Gm5} holds. 
We claim 
\[\psi(\sum_i W_i)=((-1)^{N+1}a, N).\]
Indeed, it suffices to show this after restriction to the generic point of $R$, 
in which case it follows directly from the definition of the $W_i$ and \eqref{lem:Gm3.5}.
Since $\psi(V(t\pm1))=(-(\pm1),1)$, this implies the surjectivity of $\psi$ and proves the lemma.

\end{proof}

\begin{cor}\label{cor;gammaGm}
\par
\begin{itemize}
\item[(1)]
There is a natural isomorphism
\[ \ulaNis h_0^{\bcube}(\cubegmred)^{sp} \simeq \ulomega^* \G_m.\]
\item[(2)]
For $F\in \CItsp_\Nis$, 
we have a natural isomorphism
\begin{equation}\label{eq;gammaGm}
\gamma F \simeq \uHom_{\ulMPST}(\ulomega^* \G_m, F).
\end{equation}\end{itemize}
\end{cor}
\begin{proof}
(1) is a direct consequence of Lemma \ref{lem;Kay}.
In view of \eqref{eq;gammaFCI}, (2) follows from (1) and the adjunctin of $\ulaNis$ and 
that from \S1\eqref{semipure}.
\end{proof}
\bigskip

\begin{lemma}\label{lem;gammaK}
Consider an exact sequence $0\to A\to B \to C\to 0$ in $\ulMNST$.
\begin{itemize}
\item[(1)]
Assume $A,B,C\in \CIt_{\Nis}$. Then the following sequence in $\NST$
\[0\to \ulomega_! \gamma A \to \ulomega_! \gamma B \to \ulomega_! \gamma C \to 0\]
 is exact.
\item[(2)]
Assume  $\ulomega_! A=0$ and $C\in \CItspNis$. Then the following sequence 
\[0\to \gamma A(K) \to \gamma B(K) \to \gamma C(K) \to 0\]
is exact for any function field $K$ over $k$.
\end{itemize}
\end{lemma}
\begin{proof}
First assume $A,B,C\in \CIt_{\Nis}$. Then all terms of the sequence in (1) are in $\RSC_\Nis$.
By Lemma \ref{lem;RSCexactness}, it suffices to show the exactness of
\[ 0\to \gamma A(K) \to \gamma B(K) \to \gamma C(K) \to 0\]
for a function field $K$ over $k$.

By \eqref{eq;gammaFCI}, we have $\gamma F(K)=\Hom({\cubegmred}_{,K},F)$ for all $F\in \CIt$ where ${\cubegmred}_{,K}=\cubegmred\otimes \Spec K$. 
Since ${\cubegmred}_{,K}$ is a direct summand of $\Ztr(\P^1_K,0+\infty)$, it is enough to show that \[
\Ext^1_{\ulMNST}(\Ztr(\P^1_K,0+\infty),A)=0.
\]
By using \cite[Th.2(2)]{kmsyI} we can compute
\[ \Ext^1_{\ulMNST}(\Ztr(\P^1_K,0+\infty),A) \simeq H^1_\Nis(\P_K^1,A_{(\P^1_K,0+\infty)}),\]
where we used the fact that any proper birational map $X\to \P^1_K$ is an isomorphism.
Thus the vanishing follows from Lemma \ref{P1cohomolgyvanishing}. This proves (1).

Next we assume $\ulomega_! A=0$ and $C\in \CItspNis$.
For a function field $K$ over $k$, we have a commutative diagram
\[\xymatrix@C=.7cm{
0\ar[r] & A(\P^1_K,\infty) \ar[r] \ar[d] &B(\P^1_K,\infty) \ar[r] \ar[d] 
&C(\P^1_K,\infty) \ar[r] \ar[d]^c & 0\\
0\ar[r] & A(\P^1_K,0+\infty) \ar[r]  &B(\P^1_K,0+\infty) \ar[r] &C(\P^1_K,0+\infty) \ar[r] & 0\\}\]
where the sequences are exact since for every effective Cartier divisor $D$ on $\P^1_K$,
\[ \Ext^1_{\ulMNST}(\Ztr(\P^1_K,D),A) \simeq H^1_\Nis(\P^1_K,A_{(\P^1_K,D)})=0,\]
by \cite[Th.2(2)]{kmsyI} and the fact that $A_{(\P^1_K,D)}$ is supported on the zero-dimensional scheme $|D|$ by the assumption. Finally, $\Ker(c)=0$ by \cite[Th. 3.1]{shuji}, hence the snake lemma gives the exact sequence of (2).
\end{proof}

\begin{prop}\label{lem;H1gamma}
\begin{itemize}
\item[(1)]
Take $F\in \CItspNis$ (cf. \S1\eqref{CItspNis}).  For $\sX=(X,D_X)\in \ulMCorls$ (cf. \S1\eqref{ulMCorls}), there exists a map functorial in $\sX$:
\begin{equation}\label{eq;H1gamma}
 \gamma F (\sX) \to H^1(\P^1\times X, F_{\P^1\otimes \sX}).
\end{equation}
Moreover, if $X$ is henselian local, it is an isomorphism.
\item[(2)]
Let $F\in \ulMNST$ be such that $F^{sp}\in \CItspNis$.
For $X\in \Sm$, there exists a map functorial in $X$:
\begin{equation}\label{eq2;H1gamma}
 \gamma F (X) \to H^1(\P^1\times X, F_{\P^1\times X}).
\end{equation}
Moreover, it is an isomorphism either if $F\in \CIt_\Nis$ and $X$ is henselian local,
or if $X=\Spec(K)$ is the spectrum of a function field over $k$ and the natural map $F(K) \to F(\bcube\otimes K)$
is an isomorphism.
\end{itemize}
\end{prop}
\begin{proof}
Let $L=(\P^1,0)$. We prove (1).
By \eqref{eq;gammaFCI} and \cite[Lem. 7.1]{shuji}, there exists an exact sequence of sheaves on $(\P^1\times X)_\Nis$:
\begin{equation}\label{eq;gysinseq}
0\to F_{\P^1\otimes \sX} \to  F_{L\otimes \sX} \to  i_*\gamma F_{\sX} \to 0,
\end{equation}
where $i: X \to \P^1\times X$ is induced by $0\in\P^1$.
Taking cohomology, we get the map \eqref{eq;H1gamma}.
If $X$ is henselian local, we have
\begin{equation}\label{eq;cubecohvashining}
H^1(\P^1\times X, F_{L\otimes \sX})\simeq H^1(X, F_{\sX})=0
\end{equation}
thanks to \cite[Th. 9.3]{shuji}. 
Note that the map $F(\sX)\to F(L\otimes \sX)$ induced by the projection
$L\otimes \sX \to \sX$ is an isomorphism by the $\bcube$-invariance of $F$.
Since the projection factors as $L\otimes \sX \to \P^1\otimes \sX \to \sX$, this implies
the map $F(\P^1\otimes \sX)\to F(L\otimes \sX)$ is surjective. This implies that the map \eqref{eq;H1gamma} is an isomorphism.

We prove (2).
Consider the exact sequence of sheaves on $(\P^1\times X)_\Nis$:
\begin{equation}\label{eq;Fexactsequence}
0\to F_{\P^1\times X} \to  F_{L\otimes X} \to  i_*\lambda_X F \to 0,
\end{equation}
where $\lambda_X F=i^* (F_{L\otimes X}/F_{\P^1\times X})$. 
The injectivity of the first map follows from \cite[Th. 3.1]{shuji} noting 
$F_{\P^1\times X}=F^{sp}_{\P^1\times X}$
\footnote{The point is that $X$ has the empty modulus.}
and $F^{sp}\in \CItspNis$ by the assumption.
Taking cohomology over an \'etale $U\to X$, we get a map natural in $U$:
\[ \lambda_X F(U) \to H^1(\P^1\times U, F_{\P^1\times U}).\]
To define the map \eqref{eq2;H1gamma}, it suffices to show the following.

\begin{claim}\label{claim;H1gamma}
There exists a natural map of sheaves on $X_\Nis$:
\[\phi_{F,X}: (\gamma_\Nis F)_X \to \lambda_X F.\]
It is an isomorphism if $F\in \CIt_\Nis$.
If $F\in \ulMNST$ and $F^{sp}\in \CIt_\Nis$, then
$\phi_{F,K}: (\gamma_\Nis F)_K = (\gamma F)_K\to \lambda_K F$ is an isomorphism for a function field $K$ over $k$.
\end{claim}

By definition, $\lambda_X F$ is the sheaf on $X_\Nis$ associated to the presheaf 
\begin{equation}\label{eq3;H1gamma}
\widetilde{\lambda_X F}: U \to \colim_{V} F(V,0_V)/F(V,\emptyset),
\end{equation}
where $V$ ranges over \'etale neighborhoods of $0_U=i(U)\subset \P^1\times U$. On the other hand, we have
\[(\gamma F)_X( U) = F(\P^1\times U,0+\infty)/F(\P^1\times U,\infty).\]
Since the colimit in \eqref{eq3;H1gamma} does not change when taken over \'etale neighborhood of
$0_U\subset \A^1\times U$, there is a natural map
\[(\gamma F)_X( U) \to 
F(\A^1\times U,0)/F(\A^1\times U,\emptyset) \to \widetilde{\lambda_X F}(U),\]
which induces the desired map $\phi_{F,X}$.

Next we show $\phi_{F,X}$ is an isomorphism if $F\in \CIt_\Nis$, or
if $F\in \ulMNST$ with $F^{sp}\in \CItspNis$ and $X=K$ is a function field over $k$. 
If $F$ is semi-pure, the assertion follows from \cite[Lem. 7.1]{shuji}.
In general we consider the exact sequence in $\ulMNST$:
\begin{equation}\label{eq;CFFsp}
 0\to C \to F \to F^{sp}\to 0 \qwith \;\;\ulomega_! C=0.
\end{equation}
It gives rise to a commutative diagram of sheaves on $(\P^1\times X)_\Nis$:
\begin{equation*}\label{eq;H1gammaCD}
\xymatrix{
0\ar[r] & C_{\P^1\times X} \ar[r]\ar[d] & F_{\P^1\times X} \ar[r]\ar[d] 
& F^{sp}_{\P^1\times X} \ar[r]\ar[d] & 0\\
0\ar[r] & C_{L\otimes X} \ar[r] & F_{L\otimes X} \ar[r]& F^{sp}_{L\otimes X} & \\}\end{equation*}
where the upper (resp. lower) sequence is exact by the exactness of 
$\ulomega_!:\ulMNST\to \NST$ from \S1\eqref{NST} (resp. by \eqref{b*leftexact}).
The right vertical map is injective by \cite[Th. 3.1]{shuji}.
This implies the exactness of the lower sequence of the following commutative diagram on $X_{\Nis}$:
\[\xymatrix{
0\ar[r] & (\gamma C)_{X} \ar[r]\ar[d]^{\phi_{C,X}} & (\gamma F)_{X} \ar[r]\ar[d]^{\phi_{F,X}}  & (\gamma F^{sp})_{X} \ar[r]\ar[d]^{\phi_{F^{sp},X}}  & 0\\
0\ar[r] & \lambda_X C \ar[r] & \lambda_X F \ar[r]& \lambda_XF^{sp} & \\}\]
The upper sequence is exact by Lemma \ref{lem;gammaK}.
Since we know that $\phi_{F^{sp},X}$ is an isomorphism, it suffices to show that 
$\phi_{C,X}$ is an isomorphism. Indeed, for an \'etale $U\to X$ with $U$ henselian local, we have
\begin{multline*}
 (\gamma C)_X(U) = C(\P^1\times U,0+\infty)/C(\P^1\times U,\infty)\\
\simeq \colim_{V} C(V,0_V)/C(V,\emptyset) =\widetilde{\lambda_X C}(U),
\end{multline*}
where $V$ are as in \eqref{eq3;H1gamma} and the isomorphism comes from the excision noting that $C_{(\P^1\times U,0+\infty)}$ (resp. $C_{(\P^1\times U,\infty)}$) is supported on $\{0_U,\infty_U\}$ (resp. $\infty_U$). 
This proves that $\phi_{C,X}$ is an isomorphism and completes the proof of the claim.
\medbreak

To show the second assertion of (2), we look at the cohomology exact sequence arising from \eqref{eq;Fexactsequence}. Note that
$F(\P^1\times X) \to F(L\otimes X)$ is surjective since
$F(X) \isom F(L\otimes X)$ by the assumption.
Hence it suffices to show 
$H^1(\P^1\times X, F_{L\otimes X})=0$. If $F$ is semi-pure, this follows from \eqref{eq;cubecohvashining}.
In general it is reduced to the above case using \eqref{eq;CFFsp} and noting
$H^1(\P^1\times X,C_{L\otimes X})=0$ since $C_{L\otimes X}$ is supported on $0\times X$.
This completes the proof of the lemma.
\end{proof}

\begin{cor}\label{cor;H1gamma}
Let $G\in \CIt$ and $K$ be a function field $K$ over $k$.
\begin{itemize}
\item[(1)]
There is a natural isomorphism
\[ \gamma \ulaNis G(K) \simeq H^1(\P^1_K,\ulaNis G).\]

\item[(2)]
The natural map 
\[\gamma \ulaNis G(K) \to \gamma \ulaNis G^{sp}(K)\] 
is an isomorphism.
\end{itemize}
\end{cor}
\begin{proof}
Letting $F=\ulaNis G$, we have $F^{sp}=\ulaNis G^{sp}\in \CItspNis$ by 
\S1\eqref{ulaNisCItsp}.
By Lemma \ref{lem;CIsheafication}, $F$ satisfies the second assumption of Proposition \ref{lem;H1gamma}(2).
 Hence (1) follows from Proposition \ref{lem;H1gamma}(2).
(2) follows from isomorphisms
\begin{multline*}
 \gamma \ulaNis G(K) \simeq H^1(\P^1_K,  \ulaNis G) \simeq 
H^1(\P^1_K,  \ulomega_!\ulaNis G)\simeq H^1(\P^1_K,  \aVNis \ulomega_! G)\\
\simeq H^1(\P^1_K,  \aVNis \ulomega_! G^{sp})\simeq H^1(\P^1_K,  \ulaNis G^{sp})\simeq \gamma \ulaNis G^{sp}(K),
\end{multline*} 
where the third (resp. last) isomorphism follows from \S1\eqref{NST} (resp. Proposition \ref{lem;H1gamma}).
\end{proof}

\begin{lemma}\label{lem;gammasheafication}
Let $F\in \CIt$.
\begin{itemize}
\item[(1)]
The natural map 
\[\gamma F (K) \to \gamma \ulaNis F(K)\] 
is an isomorphism for any function field $K$ over $k$.
\item[(2)]
The natural map 
$ \ulaNis\gamma F^{sp} \to  \gamma \ulaNis F^{sp}$ 
is injective. 
\item[(3)]
The natural map 
$\ulomega_! \ulaNis\gamma F^{sp} \to \ulomega_! \gamma \ulaNis F^{sp}$ 
is an isomorphism. 
\end{itemize}
\end{lemma}
\begin{proof}
Consider the exact sequence in $\ulMPST$:
\begin{equation}\label{eq2;CFFsp}
 0\to C \to F \to F^{sp}\to 0 \qwith \;\;\ulomega_! C=0.
\end{equation}
By \S1\eqref{FCItsp}, we have $C,F^{sp}\in \CIt$.
It gives rise to an exact sequence in $\ulMNST$:
\begin{equation*}
 0\to \ulaNis C \to  \ulaNis F \to  \ulaNis F^{sp}\to 0
\end{equation*}
and a commutative diagram
\[\xymatrix{
0\ar[r] & \gamma C(K) \ar[r]\ar[d] & \gamma F(K) \ar[r]\ar[d] &\gamma  F^{sp}(K) \ar[r]\ar[d] & 0\\
0\ar[r] & \gamma \ulaNis C(K) \ar[r] & \gamma  \ulaNis F(K) \ar[r] 
&\gamma   \ulaNis F^{sp}(K) \ar[r] & 0\\}\]
The upper sequence is exact thanks to \eqref{eq;gammaFCI}.
The lower sequence is exact by Lemma \ref{lem;gammaK}(2) noting $\ulaNis F^{sp}\in \CItspNis$ by \cite[Th. 10.1]{shuji} and
$\ulomega_! \ulaNis C=\aVNis\ulomega_! C=0$ (cf. \S1\eqref{NST}). 
Since $C_{(\P^1_K,0+\infty)}$ is supported on $\{0_K,\infty_K\}$,
we have by \S1\eqref{ulaNisformular}
\[(\ulaNis C)_{(\P^1_K,0+\infty)} = C_{(\P^1_K,0+\infty)},\]
where we used the fact that any proper birational map between normal schemes of dimension $1$ is an isomorphism.
Hence the left vertical map is an isomorphism. Hence we may assume that
$F$ is semi-pure. 
By \S1\eqref{ulaNisCItsp},
we have $\ulaNis F\in \CItspNis$.
By \cite[Lem. 5.9]{shuji}, we have natural isomorphisms
\[ \gamma F (K) \simeq F(\A_K^1,0)/F(\A^1_K,\emptyset),\]
\[ \gamma \ulaNis F (K) \simeq \ulaNis F(\A_K^1,0)/\ulaNis F(\A^1_K,\emptyset).\]
Hence (1) follows from \cite[Th. 4.1]{shuji}.

To show (2) and (3), first note that $F^{sp}\in \CItsp$ by the assumption and \S1\eqref{FCItsp} and hence $\gamma F^{sp} \in \CItsp$.
By \S1\eqref{ulaNisCItsp}, $ \ulaNis\gamma F^{sp}$ and $\gamma \ulaNis F^{sp}$ are in $\CItsp_\Nis$, and hence $ \ulomega_!\ulaNis\gamma F^{sp}$ and $\ulomega_!\gamma \ulaNis F^{sp}$ are in $\RSC_\Nis$. 
Hence (2) (resp. (3)) follows from (1) for $F=F^{sp}$ and \cite[Cor. 3.4]{shuji}
(resp. Lemma \ref{lem;RSCexactness}).

\end{proof}

\begin{lemma}\label{lem;gammaK2}
Consider a sequence $ A\to B \to C$ in $\CIt$ such that  
\[\ulomega_!\ulaNis A\to \ulomega_!\ulaNis B \to \ulomega_!\ulaNis C\to 0\]
is exact  in $\NST$. Then the following sequence 
\[ \gamma \ulaNis A(K) \to \gamma \ulaNis B(K) \to \gamma \ulaNis C(K) \to 0\]
 is exact for any function field $K$ over $k$.
\end{lemma}
\begin{proof}
In view of the right exactness of the functor
\[ H^1(\P_K, -) : \NST \to \Ab,\]
the lemma follows from Corollary \ref{cor;H1gamma}(1) by applying the above functor to the first exact sequence. 
\end{proof}

\begin{cor}
\label{cor;gammadoesntseeconductors}
Let $F\in \CItsp_\Nis$. Then for any function field $K$ we have an isomorphism $\gamma F(K)\cong \gamma \ulomega^{\CI}\ulomega_! F(K)$
\end{cor}
\begin{proof}
Let $q:\gamma(F)(K)\to \gamma(\ulomega^{\CI}\ulomega_! F)(K)$ be the map induced by the unit map $F\hookrightarrow \ulomega^{\CI}\ulomega_! F$ for the adjunction \eqref{omegaCIadjointNis}, which is injective since it factors the map $F\hookrightarrow \ulomega^*\ulomega_!F$. Notice that $q$ is injective by \eqref{eq;gammaFCI} and the fact that $\Hom_{\ulMPST}({\cubegmred}_{,K},\_)$ preserves injective maps, hence it is enough to show that it is surjective. Let $Q$ be the presheaf cokernel of $F\to \ulomega^{\CI}\omega_! F$, hence $Q\in \CIt$ and $\ulomega_!Q=0$. By Lemma \ref{lem;gammaK2} we have an exact sequence\[
\gamma F(K) \xrightarrow{q} \gamma \ulomega^{\CI}\omega_! F(K) \to \gamma \ulaNis Q(K) \to 0.
\]
By Corollary \ref{cor;H1gamma}(2) we have that \[
\gamma \ulaNis Q(K)\cong \gamma \ulaNis Q^{sp}(K)=0,
\]  
hence $q$ is surjective.
\end{proof}

\begin{prop}\label{prop;omegagamma}
For $F\in \CItsp_\Nis$, there is a natural isomorphism
\[ \ulomega_! \gamma F \simeq 
\ulomega_!\uHom_{\ulMPST}(\ulomega^*\G_m,F)
\simeq \uHom_{\PST}(\G_m,\ulomega_! F).\]
\end{prop}
\begin{proof}
The first isomorphism follows from \eqref{eq;gammaFCI} and Corollary \ref{cor;gammaGm}.
For $F\in \ulMPST$ and $X\in \Sm$, put 
\[F^X=\uHom_{\ulMPST}(\Ztr(X,\emptyset)),F).\]
Note that $F\in \CItsp_\Nis$ implies $F^X\in \CItsp_\Nis$.
We compute
\begin{multline*}
\ulomega_!\gamma F(X) =
\uHom_{\ulMPST}(\Lred, F)(X,\emptyset) \\
\simeq \Hom_{\ulMPST}(\Lred, F^X) =\gamma F^X (k),
\end{multline*}
\begin{multline*}
\uHom_{\PST}(\G_m,\ulomega_! F)(X) =
\Hom_{\PST}(\G_m,\uHom_{\PST}(X, \ulomega_! F))\\
\simeq \uHom_{\PST}(\G_m,\ulomega_!F^X)(k),
\end{multline*}
where the last isomorphism comes from Lemma \ref{lem;omegauHom}.
Hence it suffices to show that there exists a natural isomorphism for any $F\in \CItspNis$:
\[\gamma F (k) \simeq \Hom_{\PST}(\G_m,\ulomega_! F).\]
We have isomorphisms
\begin{multline*}\begin{aligned}
\Hom_{\PST}(\G_m,\ulomega_! F) &\overset{(*1)}{\simeq}
\Hom_{\ulMPST}(\ulomega^*\G_m,\ulomega^*\ulomega_! F)\\
& \overset{(*2)}{\simeq} \Hom_{\ulMPST}(\ulomega^*\G_m,\ulomegaCI\ulomega_! F)\\
& \overset{(*3)}{\simeq} \Hom_{\ulMPST}(\Lred,\ulomegaCI\ulomega_! F)\\
&\overset{(*4)}{\simeq} \gamma\ulomegaCI\ulomega_!F(k)\overset{(*5)}{\simeq} \gamma F(k),
\end{aligned}
\end{multline*}
where $(*1)$ follows from the fact that $\ulomega^*$ is fully faithful (cf. \S1\eqref{ulomega}), and
$(*2)$ from the adjunction from \S1\eqref{CI} (see also \eqref{omegaCIadjoint}) in view of the fact $\ulomega^*\G_m\in \CIt$ by Lemma \ref{lem;Kay}, $(*3)$ from Lemma \ref{lem;Kay}, $(*4)$ by \eqref{eq;gammaFCI} and $(*5)$ by Corollary \ref{cor;gammadoesntseeconductors},

\end{proof}

\section{Weak cancellation theorem}\label{sec3}

For $F, G \in \MPST$ we write (cf. \S1\eqref{CItspNis}, \eqref{tensor} and Lemma \ref{lem2;hcubesp})
\[F \otCI G = h_0^\bcube(F \otimes_{\MPST} G)\in \CI,\]
\[F \otCIsp G =\hcubesp(F \otimes_{\MPST} G)  \in \CIsp,\]
\[F \otCIspNis G = \hcubespNis(F \otimes_{\MPST} G)  \in \CIsp_\Nis.\]

\begin{prop}\label{prop;otCI}
The product $\otCI$ (resp. $\otCIsp$, resp. $\otCINissp$) defines a symmetric monoidal structure on $\CI$ (resp. $\CIsp$, resp. $\CIspNis$).
\end{prop}
\begin{proof}
The assertion follows immediately from the fact that $\otMPST$ defines a symmetric monoidal structure on $\MPST$ except the associativity.
We prove it only for $\otCINissp$ (other cases are similar).
We need to show a natural isomorphism for $F,G,H\in \CIspNis$:
\[(F\otCINissp G)\otCINissp H \simeq F\otCINissp (G\otCINissp H).\]
For simplicity we write $\lambda=\hcubespNis$.
For $F,G,H\in \CIspNis$, we have isomorphisms
\[\begin{aligned}
 \lambda(\lambda(F\otMPST G) \otMPST H)
& \overset{(*1)}{\simeq}  \lambda(\lambda^2(F\otMPST G) \otMPST \lambda H)\\
& \overset{(*2)}{\simeq}  \lambda(\lambda(F\otMPST G) \otMPST \lambda H)\\
& \overset{(*3)}{\simeq}  \lambda((F\otMPST G) \otMPST H)\\
\end{aligned}\]
where $(*1)$ (resp. $(*2)$, resp. $(*3)$) follows from Lemma \ref{lem2;hcubesp}
$(iii)$ (resp. $(i)$, resp. $(iii)$).
The lemma follows from this and the associativity of $\otMPST$. 
\end{proof}

For $F, G \in \CIt$ we write
\[F \otCI G = \tau_! \hcube(\tau^* F \otimes_{\MPST} \tau^*G)\in \CIt,\]
\[F \otCIsp G =  \tau_! \hcubesp (\tau^* F \otimes_{\MPST} \tau^*G) \in \CItsp,\]
\[F \otCIspNis G =  \tau_! \hcubespNis(\tau^* F \otimes_{\MPST} \tau^*G) \in \CItsp_\Nis.\]
By \S1\eqref{tauaNis}, we have a natural isomorphism
\eq{ulaNisotCI}{ \ulaNis(F \otCIsp G ) \simeq F \otCIspNis G .}
In view of the equivalences \eqref{CIspCItsp} and \eqref{CIspCItspNis}, 
Proposition \ref{prop;otCI} implies

\begin{prop}\label{prop2;otCI}
The product $\otCI$ (resp. $\otCIsp$, resp. $\otCINissp$) defines a symmetric monoidal structure on $\CIt$ (resp. $\CItsp$, resp. $\CItsp_\Nis$).
There is a natural isomorphism for $F,G,H\in \CItspNis$
\eq{eq;otCINissp}{(F\otCINissp G)\otCINissp H \simeq F\otCINissp (G\otCINissp H).}
\end{prop}
\medbreak

For $F\in \CIt_\Nis$ and an integer $d\geq 0$, we put
\eq{twist(1)}{ F(d) =  (\cubegmred)^{\otCINissp d}\otCINissp F .}
Note $F(d) = F(m)(n)$ with $d=m+n$ by \eqref{eq;otCINissp}.
\medbreak

For $F\in \CIt$ and $f\in F(\sX)$ with $\sX\in \ulMCor$, consider the composite map
\[ \cubegmred\otimes_{\ulMPST} \Ztr(\sX) \rmapo{id_{\cubegmred}\otimes f}
\cubegmred\otimes_{\ulMPST} F \to \cubegmred\otCI  F.\]
By the adjunction $(\Lred\otimes_{\ulMPST}  -)\dashv \uHom_{\ulMPST}(\Lred,-)$
this gives rise to a natural map 
\begin{equation}\label{eq;iotaF}
\iota_F :  F \to \gamma (\cubegmred\otCI  F),
\end{equation}
which induces  
\begin{equation}\label{eq;iotapreF}
\iota_F^{sp}: F^{sp}\to \gamma(\Lred \otCIsp F),
\end{equation}
noting the adjunction from \S1\eqref{semipure} and the fact that $\gamma:\ulMPST\to \ulMPST$ preserves semipure objects.

If $F\in \CIt_\Nis$, this induces a natural map
\begin{equation}\label{eq;iota}
\iota_F : F^{sp} \to \gamma (F(1)).
\end{equation}
which generalizes to a natural map for $n\in \Z_{\geq 1}$
(cf. \eqref{twist(1)} and \eqref{eq;gammaFCIn2})
\begin{equation}\label{eq;iotan}
\iota^n_F : F^{sp} \to \gamma^n (F(n)),
\end{equation}
noting
\begin{equation*}\label{eq;gammaFCIn}
\gamma^n F= \uHom_{\ulMPST}((\cubegmred)^{\otCI n},F) \qfor F\in \CIt
\end{equation*}
thanks to  the adjunction from \eqref{CI}. 

\begin{qn}\label{ques1}
For $F\in \CItsp_\Nis$, is the map \eqref{eq;iota} an isomorphism?
\end{qn}

We will prove the following variant.

\begin{thm}\label{thm;weakcancel}
For $F\in \CIt$, the map \eqref{eq;iotapreF} is an isomorphism.
\end{thm}

Before going into its proof, we give some consequences.

\begin{cor}\label{cor;weakcancel}
For $F\in \CIt$ the map \eqref{eq;iotapreF} gives an isomorphism\[
\ulomega_!\iota_F: \ulomega_!\ulaNis F\xrightarrow{\sim} 
\ulomega_!\gamma \ulaNis( \Lred \otCIsp F).
\]
For $F\in \CIt_\Nis$, 
the map \eqref{eq;iotan} induces an isomorphism\[
\ulomega_!\iota^n_F: \ulomega_! F\xrightarrow{\sim} \ulomega_!\gamma^n F(n). \]
\end{cor}
\begin{proof} 
The functors $\ulomega_!$ and $\ulaNis$ are exact and $\ulomega_!\ulaNis G\cong \ulomega_!\ulaNis G^{sp}$ for all $G\in \ulMPST$.
Hence Theorem \ref{thm;weakcancel} gives a natural isomorphism 
\[
\ulomega_!\ulaNis\iota_F: \omega_!\ulaNis F\isom 
\ulomega_!\ulaNis \gamma (\Lred \otCIsp F).
\]
This proves the first assertion since Lemma \ref{lem;gammasheafication}(3) implies
\[\ulomega_!\ulaNis \gamma (\Lred\otCIsp F)\simeq 
\ulomega_!\gamma \ulaNis(\Lred\otCIsp F).\]
The second assertion for the case $n=1$ follows directly from the first.
For $n>1$, we proceed by the induction on $n$ to assume 
\begin{equation}\label{eq0;cor;weakcancel}
\ulomega_!\iota^{n-1}_F: \ulomega_! F\xrightarrow{\sim} \ulomega_!\gamma^{n-1} F(n-1).
\end{equation}
Then we have isomorphisms
\begin{multline*}
\ulomega_!\gamma^n F(n) \overset{(*1)}{\simeq} 
\ulomega_!\gamma\gamma^{n-1}F(n)
\overset{(*2)}{\simeq} \uHom_{\PST}(\G_m,\ulomega_!\gamma^{n-1} F(n))=\\
\uHom_{\PST}(\G_m,\ulomega_!\gamma^{n-1} F(1)(n-1))
\overset{(*3)}{\simeq} \uHom_{\PST}(\G_m,\ulomega_!F(1))\\
\overset{(*4)}{\simeq}\ulomega_! \gamma F(1)\overset{(*5)}{\simeq} \ulomega_!F,
\end{multline*}
where $(*1)$ follows from
\eqref{eq;gammaFCIn2}, $(*2)$ follows from Proposition \ref{prop;omegagamma} noting 
$\gamma^{n-1}F(n)\in \CItspNis$, $(*3)$ follows from  \eqref{eq0;cor;weakcancel}, $(*4)$ follows from  Proposition \ref{prop;omegagamma} and $(*5)$ follows from the case $n=1$.
This completes the proof.\end{proof}

\begin{cor}\label{cor3;weakcancel}
For $F\in \RSC_\Nis$ and $\tF=\ulomega^{\CI} F\in \CIt_\Nis$ (cf. \eqref{omegaCIadjointNis}),
the map \eqref{eq;iotan} $\iota^n_{\tF}: \tF\to \gamma^n \tF(n)$ is an isomorphism.
\end{cor}
\begin{proof} 
We have a commutative diagram
\[\xymatrix{
\tF \ar[r]^{\iota^n_{\tF}}\ar[d]^{\cong} & \gamma^n\tF(n) \ar[d]^{\hookrightarrow} \\
\ulomega^{\CI}\ulomega_! \tF \ar[r]^{\ulomega^{\CI}\ulomega_!\iota_{\tF}}
& \hskip 5pt\ulomega^{\CI}\ulomega_! \gamma^n\tF(n) \\}\]
where the vertical arrows come from the adjunction \eqref{omegaCIadjointNis}.
The left (resp. right) vertical arrow is an isomorphism (resp. injective) since 
$\ulomega_!\ulomega^{\CI}\simeq id$ (resp. by the semipurity of $\gamma^n \tF(n)$).
Since $\ulomega^{\CI}\ulomega_!\iota^n_{\tF}$ is an isomorphism by
 Corollary \ref{cor;weakcancel}, this implies $\iota^n_{\tF}$ is an isomorphism by Snake Lemma.
\end{proof}

\def\trho{\tilde{\rho}}
\def\hF{\hat{F}}

\begin{cor}\label{cor2;weakcancel}
For $F\in \CItsp_\Nis$, there is a natural injective map
\[ \trho_F: \gamma^n F(n) \to \ulomega^{\CI}\ulomega_! F\]
whose composite with the map \eqref{eq;iotan} $\iota^n_F:F \to \gamma^n F(n)$ coincides with the unit map $u :F\to \ulomega^{\CI}\omega_! F$ for the adjunction \eqref{omegaCIadjointNis}.
In particular \eqref{eq;iotan} is injective.
\end{cor}
\begin{proof} 
Define $\trho_F$ as the composite
\[ \gamma^n F(n) \rmapo{u} \gamma^n \ulomega^{\CI}\ulomega_! F(n) \rmapo{(\iota^n_{\ulomega^{\CI}\omega_! F})^{-1}} \ulomega^{\CI}\omega_! F,\]
where the second map is the inverse of the isomorphism
$\iota^n_{\ulomega^{\CI}\omega_! F} : \ulomega^{\CI}\omega_! F \cong \gamma^n\ulomega^{\CI}\omega_! F(n)$ from Corollary \ref{cor3;weakcancel}.
Clearly we have $\trho_F\circ\iota^n_F=u$. We easily see that  
$\trho_F$ coincides with the composite
\[ \gamma^n F(n) \rmapo{u} \ulomega^{\CI}\ulomega_! \gamma^n F(n) \rmapo{\ulomega^{\CI} (\ulomega_! \iota^n_F )^{-1}} \ulomega^{\CI}\ulomega_! F,\]
where the first map is injective by the semipurity of $\gamma^n F(n)$ and 
the second map is induced by the inverse of the isomorphism
$\ulomega_! \iota^n_F: \ulomega_! F \cong \ulomega_! \gamma^n F(n)$ from
Corollary \ref{cor;weakcancel}.
This completes the proof.
\end{proof}

\bigskip

In the rest of this section we prove the following.

\def\cubegm{\bcube^{(1)}}

\begin{prop}\label{prop;conj2splitting}
For $F\in \CIt$, the map \eqref{eq;iotapreF} $\iota_F^{sp}$ is split injective. 
\end{prop}

\medbreak

For the proof of Proposition \ref{prop;conj2splitting} we first recall the construction of \cite{voecancel}.
Take $X,Y\in \Sm$.
For an integer $n>0$ consider the rational function on $\A^1_{x_1} \times \A^1_{x_2}$:
\[ g_n =\frac{x_1^{n+1}-1}{x_1^{n+1}-x_2}.\]
Let $D_{XY}(g_n)$ be the divisor of the pullback of $g_n$ to 
$(\A^1_{x_1}-0)\times X \times (\A^1_{x_2}-0)\times Y$.
Take a prime correspondence
\begin{equation}\label{eq;Z}
Z\in \Cor((\A^1_{x_1}-0)\times X,(\A^1_{x_2}-0)\times Y).
\end{equation}
Let $\Zb\subset \P^1_{x_1}\times X \times \P^1_{x_2}\times Y$ be the closure of $Z$ and $\Zb^N$ be its normalization.
\def\CorN{\Cor^{(N)}}

\begin{lemma}\label{lem;voeLem4.1}
\begin{itemize}
\item[(1)]
Let $N>0$ be an integer such that 
\begin{equation}\label{eq;NZ}
 N(0_1+\infty_1)_{|\Zb^N}\geq (0_2+\infty_2)_{|\Zb^N}.
\end{equation}
Then, for any integer $n\geq N$, 
$Z$ intersects properly with $|D_{XY}(g_n)|$ and any component of 
the intersection $Z\cdot D_{XY}(g_n)$ is finite and surjective over $X$. 
Thus we get  
\[ \rho_n(Z) \in \Cor(X,Y)\]
as the pushforward of $Z\cdot D_{XY}(g_n)$ in $X\times Y$.
\item[(2)]
If $Z=Id_{(\A^1-0)}\otimes W$ for $W\in \Cor(X,Y)$, then one can take $N=1$ in (1) and $\rho_n(Z)=W$.
\item[(3)]
For any $Z$ as in \eqref{eq;Z} such that $\rho_n(Z)$ is defined and for any $f\in \Cor(X',Y')$
with $X',Y'\in \Sm$, $\rho_n(Z\otimes f)$ for 
\[ Z\otimes f \in \Cor((\A^1_{x_1}-0)\times (X\times X'),(\A^1_{x_2}-0)\times (Y\times Y'))\]
is defined and we have 
\[\rho_n(Z\otimes f) =\rho_n(Z)\otimes f \in \Cor(X\times X',Y\times Y').\]
\item[(4)]
For an integer $N>0$ let
\[ \CorN((\A^1_{x_1}-0)\times X,(\A^1_{x_2}-0)\times Y)\]
be the subgroup of $\Cor((\A^1_{x_1}-0)\times X,(\A^1_{x_2}-0)\times Y))$
generated by prime correspondences satisfying 
the condition \eqref{eq;NZ}. Then the presheaf on $\Sm$ given by
\[ X\to \CorN((\A^1_{x_1}-0)\times X,(\A^1_{x_2}-0)\times Y)\]
is a Nisnevich sheaf.
\end{itemize}
\end{lemma}
\begin{proof}
The assertions are proved in \cite[Lem. 4.1, 4.3 and 4.5]{voecancel} except that 
(4) follows from the fact that the condition \eqref{eq;NZ} is Nisnevich local on $X$.
\end{proof}

For an integer $a\geq 1$ put $\cubegma=(\P^1,a(0+\infty))\in \MCor$.
Take $\sX=(\Xb,\Xinf), \sY=(\Yb,\Yinf)\in \MCor$ with 
$X=\Xb-|\Xinf|$ and $Y=\Yb-|\Yinf|$.
For $a\geq 1$ take a prime correspondence
\begin{equation*}
Z\in \MCor(\cubegma\otimes \sX,\cubegm\otimes\sY).
\end{equation*}
By definition $Z\in \Cor(X,Y)$ 
and it satisfies 
\begin{equation}\label{eq;MCZ}
(0_2+\infty_2)_{|\Zb^N} + (\Yinf)_{|\Zb^N} \leq a(0_1+\infty_1)_{|\Zb^N}+ (\Xinf)_{|\Zb^N} ,
\end{equation}
where $\Zb^N$ is the normalization of the closure $\Zb$ of $Z$ in 
$\P^1_{x_1}\times X \times \P^1_{x_2}\times \Yb$. 

For integers $n,m\geq a$, we consider the rational function on $\A^1_{x_1}\times \A^1_t \times \A^1_{x_2}$:
\[ h= t g_n + (1-t) g_m .\]
Let $D_{X\A^1Y}(h)$ be the divisor of the pullback of $h$ to 
$(\A^1_{x_1}-0)\times X\times \A^1_t \times (\A^1_{x_2}-0)\times Y$.
By \cite[Rem. 4.2]{voecancel}, 
$Z\times \A^1_t$ intersects properly with $|D_{X\A^1Y}(h)|$ and any component of the intersection $(Z\times\A^1_t)\cdot D_{X\A^1Y}(h)$ is finite and surjective over $X\times\A^1_t$. 
Thus we get  
\[ \rho_h(Z\times\A^1_t) \in \Cor(X\times\A^1_t,Y).\]
It is easy to see
\begin{equation}\label{eq;rhohhomotopy}
i_0^*\rho_h(Z\times\A^1_t)=\rho_m(Z) \qaq i_1^*\rho_h(Z\times\A^1_t)=\rho_n(Z) .
\end{equation}

\begin{lemma}\label{lemma;rhohomotopy}
For $n,m\geq a$,  
$\rho_h(Z\times\A^1_t) \in \MCor(\sX\otimes\bcube,\sY)$.
\end{lemma}
\begin{proof}
Let $V$ be any component of $(Z\times\A^1_t)\cdot D_{X\A^1Y}(h)$ and $\Vb$ be its closure in 
\[\P^1_{x_1}\times \Xb \times \P^1_t \times \P^1_{x_2}\times \Yb.\]
Let $W\subset X\times \A^1_t\times Y$ be the image of $V$ and 
$\Wb$ be its closure in $\Xb\times\P^1_t\times \Yb$. Then we have
$\Wb=\pi(\Vb)$, where 
\[\pi: \P^1_{x_1}\times \Xb \times \P^1_t \times \P^1_{x_2}\times \Yb\to
\Xb\times\P^1_t\times \Yb\]
is the projection. We want to show
\[ (\Yinf)_{|\Wb^N} \leq (\Xb\times \infty)_{|\Wb^N} + (\Xinf\times \P^1_t)_{|\Wb^N}.\]
Since $\pi:\Vb^N\to \Wb^N$ is proper and surjective, this is reduced to showing 
\[  (\Yinf)_{|\Vb^N} \leq (\Xb\times \infty)_{|\Vb^N} + (\Xinf\times \P^1_t)_{|\Vb^N}\]
by \cite[Lem. 2.2]{kp}. By \eqref{eq;MCZ} and the containment lemma \cite[Pr. 2.4]{kp} (see also \cite[Lem. 2.1]{bs}), we have
\[  (\Yinf)_{|\Vb^N} + (0_2+\infty_2)_{|\Vb^N} \leq a(0_1+\infty_1)_{|\Vb^N}
+ (\Xinf\times\P^1_t)_{|\Vb^N}.\]
Thus it suffices to show 
\begin{equation*}\label{eq;MC1}
a(0_1+\infty_1)_{|\Vb^N} \leq (0_2+\infty_2)_{|\Vb^N} + \infty_{|\Vb^N}.
\end{equation*}
Using \cite[Pr. 2.4]{kp} again, this follows from 
\begin{equation}\label{eq;MC2}
a (0_1+\infty_1)_{|T} \leq (0_2+\infty_2)_{|T} + \infty_{|T},
\end{equation}
where $T\subset\P^1_{x_1}\times \P^1_t \times \P^1_{x_2}$ is any component of the closure of the divisor of $h$ on 
$(\A^1_{x_1}-0)\times \A^1_t \times (\A^1_{x_2}-0)$.
By an easy computation, $T$ is contained in one of the closures 
$\ol{D(H)}$, $\ol{D(J_n)}$, $\ol{D(J_m)}$ of the divisors of
 \[ H = t  (x_1^{n+1}-x_1^{m+1})(1-x_2) + (x_1^{m+1}-1)(x_1^{n+1}- x_2),\]  
\[ J_n = x_1^{n+1} -x_2,\quad J_m = x_1^{m+1} -x_2\]
respectively. 
Letting $\P^1_{x_i}-0=\Spec k[\tau_i]$ with $\tau_i=x_i^{-1}$ for $i=1,2$, 
$\ol{D(H)}$, $\ol{D(J_n)}$, $\ol{D(J_m)}$ are defined in 
$(\P^1_{x_1}-0)\times \A^1_t \times (\P^1_{x_2}-0)$ by the ideals generated by\[ H'= t (\tau_1^{m+1}-\tau_1^{n+1})(\tau_2-1) + (1-\tau_1^{m+1})(\tau_2-\tau_1^{n+1}),\] 
\[ J_n' = \tau_2-\tau_1^{n+1} ,\quad J_m' = \tau_2-\tau_1^{m+1} .\]
Hence, $\ol{D(H)}$, $\ol{D(J_n)}$, $\ol{D(J_m)}$
do not intersect with $\infty_1\times\P^1_t\times \A^1_{x_2}$.

By the assumption $n,m\geq a$, the ideals $(J_n,x_1^a), (J_m,x_1^a)\subset k[x_1,x_2]$ contain $x_2$ and the ideals $(J_n',\tau_1^a), (J_m',\tau_1^a)\subset k[\tau_1,\tau_2]$ contain $\tau_2$, which implies \eqref{eq;MC2} (without the last term) if $T$ is contained in $\ol{D(J_m)}$ or $\ol{D(J_n)}$.

On the other hand,  the ideal $(H,x_1^a)\subset k[x_1,x_2,t]$ contains
$x_2$ and the ideal $(H',\tau_1^a)\subset k[\tau_1,\tau_2,t]$ contains
$\tau_2$. Over $\P^1_t-0=\Spec k[u]$ with $u=t^{-1}$, 
$\ol{D(H)}\cap\big(\A^1_{x_1}\times (\P^1_t-0) \times \A^1_{x_2}\big)$ is the zero divisor of 
 \[ \tilde{H} =  (x_1^{n+1}-x_1^{m+1})(1-x_2) +u (x_1^{m+1}-1)(x_1^{n+1}- x_2),\]  
and $\ol{D(H)}\cap\big((\P^1_{x_1}-0)\times (\P^1_t-0) \times (\P^1_{x_2}-0)\big)$ is the zero divisor of 
\[ \tilde{H}'= (\tau_1^{m+1}-\tau_1^{n+1})(\tau_2-1) + u(1-\tau_1^{m+1})(\tau_2-\tau_1^{n+1}),\] 
The ideal $(\tilde{H},x_1^a)\subset k[x_1,x_2,u]$ contains $ux_2$ and 
the ideal $(\tilde{H}',\tau_1^a)\subset k[\tau_1,\tau_2,u]$ contains $u\tau_2$. 
This shows \eqref{eq;MC2} if $T\subset \ol{D(H)}$ and completes the proof of the claim.
\end{proof}
\medbreak
\def\MCorN{\MCor^{(N)}}

\begin{lemma}\label{lem;rhoMC}
For $n\geq a$ we have $\rho_n(Z) \in \ulMCor(\sX,\sY)$.
\end{lemma}
\begin{proof}
This follows from Lemma \ref{lemma;rhohomotopy} and \eqref{eq;rhohhomotopy}.
\end{proof}
\medbreak

For an integer $N\geq a$ let
\[ \MCorN(\cubegmreda\otimes \sX,\cubegmred\otimes\sY)\subset 
\MCor(\cubegmreda\otimes \sX,\cubegmred\otimes\sY)\]
be the subgroup generated by prime correspondences lying in
\[\CorN((\A^1-0)\times X,(\A^1-0)\times Y).\]
By Lemma \ref{lem;rhoMC}, we get a map for $n\geq N\geq a$
\begin{equation}\label{eq;rhon}
\rho_{n}^{(a)} : \MCorN(\cubegmreda\otimes \sX,\cubegmred\otimes\sY)
\to \MCor(\sX,\sY).
\end{equation}
The map \eqref{eq;rhon} induces a map of cubical complexes
\begin{equation}\label{eq;rhoncomplex}
\rho_{n}^{(a)\bullet} : \MCorN(\cubegmreda\otimes \sX\otimes\bcube^\bullet,\cubegmred\otimes\sY)
\to \MCor(\sX\otimes \bcube^\bullet,\sY).
\end{equation}
By the construction the following diagram is commutative if $n\geq N\geq b\geq a$:
\begin{equation}\label{eq2;rhoncomplex}
\xymatrix{
\MCorN(\cubegmreda\otimes \sX\otimes\bcube^\bullet,\cubegmred\otimes\sY)
	\ar[r]^-{\rho_{n}^{(a)\bullet }} \ar[d]^{\beta^*} & \MCor(\sX\otimes \bcube^\bullet,\sY)\\
\MCorN(\cubegmredb\otimes \sX\otimes\bcube^\bullet,\cubegmred\otimes\sY)
\ar[ru]_-{\rho_{n}^{(b)\bullet }} \\}
\end{equation}
where $\beta^*$ is induced by the natural map $\beta:\cubegmredb \to \cubegmreda$.

\begin{cor}\label{cor;rhohomotopy}
For $m,n\geq N\geq a$, $\rho_{n}^{(a)\bullet}$ and $\rho_{m}^{(a)\bullet}$ are homotopic.
\end{cor}
\begin{proof}
By Lemma \ref{lemma;rhohomotopy}, we get a map
\begin{equation}\label{eq;smn}
s_{m,n}=\rho_h(- \times \A^1_t) : \MCorN(\cubegmreda\otimes \sX,\cubegmred\otimes\sY)
\to \MCor(\sX\otimes\bcube,\sY)
\end{equation}
such that $\partial \circ s_{m,n} = \rho_{m}^{(a)}-\rho_{n}^{(a)}$, where
\[ \partial =i_0^*-i_1^* : \MCor(\sX\otimes\bcube,\sY) \to \MCor(\sX,\sY).\]
Let
\begin{equation*}
s^i_{m,n} : \MCorN(\cubegmreda\otimes \sX\otimes\bcube^i,
\cubegmred\otimes\sY) \to \MCor(\sX\otimes\bcube^{i+1},\sY)
\end{equation*}
be the map \eqref{eq;smn} defined replacing $\sX$ by $\sX\otimes\bcube^i$.
Then we have that\[
\partial\circ((-1)^is^i_{m,n}) + (-1)^{i-1}s^{i-1}_{m,n}\circ \partial = \rho_{n}^{(a),i} - \rho_{m}^{(a),i}
\]
hence $\{(-1)^is^{i}_{m,n}\}_i$ gives the desired homotopy.
\end{proof}

Let $Z\in \MCorN(\cubegmreda\otimes \sX,\cubegmred\otimes \sY)$, then for all $W\in \MCor(\sX',\sX)$, by \cite[Lemma 4.4]{voecancel}
\[Z\circ(Id_{\A^1-\{0\}}\otimes W)\in \CorN((\A^1-0)\times X,(\A^1-0)\times Y).\]
Moreover, by \cite[Prop 1.2.4(i)]{kmsyI} we have\[
Z\circ(Id_{\A^1-\{0\}}\otimes W)\in \MCor(\cubegmreda\otimes \sX,\cubegmred\otimes \sY),
\] 
which implies that
\begin{multline*}
    L_a(\sY)^{(N)}=\uHom^{(N)}_{\MPST}(\cubegmreda,\cubegmred\otimes\Ztr(\sY)) \\
=\MCorN(\cubegmreda\otimes (-),\cubegmred\otimes \sY)
\end{multline*}
is an object of $\MPST$, which is a subobject of 
\[L_a(\sY)=\uHom_{\MPST}(\cubegmreda,\cubegmred\otimes\Ztr(\sY)) \in \MPST,\]
and we have
\eq{LaYLaYN}{L_a(\sY)=\colim_{N>0} L_a(\sY)^{(N)}.}
The above construction gives a map of complexes in 
$\MPST$:
\[ \rho^{(a)\bullet}_N: C_\bullet L_a(\sY)^{(N)} \to C_\bullet(\sY),\]
where $C_\bullet(-)$ is the cubical Suslin complex. Let
\[ \rho^{(a)}_N: H_i(C_\bullet L_a(\sY)^{(N)}) \to H_i(C_\bullet(\sY))\]
be the map in $\MPST$ induced on cohomology presheaves.
Thanks to Corollary \ref{cor;rhohomotopy}, the diagram
\[\xymatrix{
H_i(C_\bullet L_a(\sY)^{(N)}) \ar[r]^-{\rho^{(a)}_N} \ar[d] & h_i^{\bcube}(\sY)\\
H_i(C_\bullet L_a(\sY)^{(N')}) \ar[ru]_-{\rho^{(a)}_{N'}} \\}\]
commutes for integers $N'\geq N$. Hence, by \eqref{LaYLaYN} we get maps
\[ \rho^{(a)} : H_i(C_\bullet L_a(\sY)) \to h_i^{\bcube}(\sY).\]
Putting $\Phi=\cubegmred\otimes\sY$, we have
\[ C_\bullet(L_a(\sY))=\uHom_{\MPST}(\cubegmreda,\uHom_{\MPST}(\bcube^\bullet,\Phi) ).\]
Recall that for $F\in \MPST$ and $\sX\in \MCor$, we have by the Hom-tensor adjunction an isomorphism:
\[
h_0^{\bcube}\uHom_{\MPST}(\Ztr(\sX),F)\cong \uHom_{\MPST}(\Ztr(\sX),h_0^{\bcube}(F)).
\]
Hence, we get an isomorphism
\[ H_0(C_\bullet L_a(\sY))\simeq \uHom_{\MPST}(\cubegmreda,h_0^{\bcube}(\Phi)),\]
where $h_i^{\bcube}(\Phi)=H_i(C_\bullet(\Phi))$ and we have an isomorphism
\[ h_0^{\bcube}(\Phi)= h_0^{\bcube}(\cubegmred\otimes \sY) = \cubegmred\otCI \sY\in \CI.\]
Hence we get a natural map
\begin{equation}\label{rhoa}
 \rho^{(a)}_\sY : \gamma_a(\cubegmred\otCI \sY) \to h_0^{\bcube}(\sY) . 
\end{equation}
where
\[  \gamma_a(F):= \uHom_{\ulMPST}(\cubegmreda,F)\qfor F\in \ulMPST ,\]
and by abuse of notation, for $C\in \CI$, we let $C$ denote also $\tau_!C\in \CIt$
(cf. \S1\eqref{CItspNis}). 
In view of \eqref{eq2;rhoncomplex}, the following diagram is commutative (recall that we assume $b\geq a$), :
\[\xymatrix{ 
\uHom_{\ulMPST}(\cubegmreda,h_0^{\bcube}(\Phi)) \ar[r]^-{\rho^{(a)}_\sY} \ar[d]^{\beta^*} & h_0^{\bcube}(\sY)\\
\uHom_{\ulMPST}(\cubegmredb,h_0^{\bcube}(\Phi)) \ar[ru]_-{\rho^{(b)}_\sY} \\}\]

Now take any $F\in \CIt$ and consider a presentation in $\ulMPST$:
\[ A\to B\to F\to 0,\]
where $A,B$ are the direct sum of $h_0^{\bcube}(\sY)$ for varying $\sY\in \MCor$.
We then get a commutative diagram
\[\xymatrix@C=.4cm{
\gamma_a(\cubegmred \otCI A) \ar[r]\ar[d]^{\rho^{(a)}_A} & 
\gamma_a(\cubegmred\otCI B) \ar[r]\ar[d]^{\rho^{(a)}_B} & 
\gamma_a(\cubegmred \otCI F) \ar[r] & 0\\
 A \ar[r] & B \ar[r] & F \ar[r] & 0,\\}\]
where the vertical maps are induced by \eqref{rhoa}.
The upper sequence is exact by the right-exactness of $\otCI$ and the fact that
$\cubegmreda$ is a projective object of $\ulMPST$.
Thus we get the induced map in $\ulMPST$:
\begin{equation}\label{rhoaF}
 \rho^{(a)}_F : \gamma_a(\cubegmred \otCI F) \to F.
\end{equation}

Write $\rho_F =\rho^{(1)}_F $. 

\begin{claim}\label{claim;rhoiota}
The map $\rho_{F}$ splits $\iota_F$.
\end{claim}
\begin{proof}
By the construction of $\rho_F$, this is reduced to the case $F=h_0^{\bcube}(\sY)$ for $\sY\in \MCor$, which follows from Lemma \ref{lem;voeLem4.1}(2).
\end{proof}

The following result concludes the proof of Proposition \ref{prop;conj2splitting}:

\begin{lemma}\label{lm;trick}
For $F\in \CIt$, $\rho_F$ factors through
\[ \rho_{F}^{sp}:\gamma(\Lred\otCIsp F)\to F^{sp}.\]
Moreover it splits the map $\iota_F^{sp}$ from \eqref{eq;iotapreF}.
\end{lemma}
\begin{proof}
Take $\sX\in \MCor$ and let $\phi$ be in the kernel of 
\[\Hom_{\ulMPST}(\Lred\otimes \sX , \Lred\otCI F)\to \Hom_{\ulMPST}(\Lred\otimes \sX , \Lred\otCIsp F).\]
Note that the map is surjective since $\Lnredd{a} \otimes \sX$ is a projective 
object of $\ulMPST$ by Yoneda's lemma. 
By the definition of semi-purification (cf. \S1\eqref{semipure}),
there exists an integer $m>0$ such that
\[\beta_m^*\phi=0\text{ in }\Hom_{\ulMPST}(\Lnredd{m}\otimes \sX^{(m)}, \Lred\otCI F),\]
where $\beta_m: \Lnredd{m}\otimes \sX^{(m)}\to \Lnredd{1}\otimes \sX$ (cf. \S1\eqref{sXn}).
Then the maps from \eqref{rhoaF} induce a commutative diagram
\[\xymatrix{
\Hom_{\ulMPST}(\Lred\otimes \sX ,\Lred\otCI F )   \ar@/^-9.0pc/[dd]_{\beta^*_m}\ar[d]\ar[r]^{\hskip 70pt\rho_{F}} & 
F(\sX) \ar[d]^{\theta_m^*} \\
\Hom_{\ulMPST}(\Lnredd{1}\otimes \sX^{(m)}, \Lred\otCI F)  
\ar[r]^{\hskip 70pt \rho_{F}}\ar[d]
& F(\sX^{(m)})\\ 
\Hom_{\ulMPST}(\Lnredd{m}\otimes \sX^{(m)}, \Lred\otCI F)\ar[ur]_{\rho^{(m)}_F}
}\]
where $\theta_m^*$ is induced by $\theta_m:\sX^{(m)}\to \sX$ and the triangle commutes by \eqref{eq2;rhoncomplex}. We have 
\[\theta_m^*\rho_F(\phi)=\rho^{(m)}_F\beta^*_m(\phi)=0.\]
Hence $\rho_F(\phi)$ lies in the kernel of $\theta^*_m$, which is contained in the kernel of the map
\[sp_\sX:F(\sX)\to F^{sp}(\sX)\]
by the definition of semi-purification. Hence the composite map
\[sp_\sX\circ \rho_F: \Hom_{\ulMPST}(\Lred\otimes \sX ,\Lred\otCI F )\to F^{sp}(\sX) \]
factors through $\Hom_{\ulMPST}(\Lred\otimes \sX ,\Lred\otCIsp F)$ inducing the desired map $\rho_F^{sp}$.
Finally, to show the last assertion, consider the commutative diagram
\[\xymatrix{
F \ar[r]^{\hskip -20pt \iota_F} \ar[d]&  \gamma( \Lnredd{1}\otCI F)\ar[d]  \ar[r]^{\hskip 30pt \rho_F} \ar[d]& F \ar[d] \\
F^{sp}  \ar[r]^{\hskip -20pt \iota_F^{sp}}
& \gamma( \Lnredd{1}\otCIsp F)\ar[r]^{\hskip 30pt \rho_F^{sp}}
& F^{sp}\\ }\]
where $\rho_F\iota_F=id_F$ by Claim \ref{claim;rhoiota}.
This implies $\rho_F^{sp}\iota_F^{sp}=id_{F^{sp}}$ since $F\to F^{sp}$ is surjective.
This completes the proof of Lemma \ref{lm;trick}.
\end{proof}

\medbreak

\section{Completion of the proof of the main theorem}\label{sec4}

In this section we prove the following result:

\begin{prop}\label{prop;almostsurj}
For $\phi \in \Hom_{\ulMPST}(\Lred\otimes \sX ,\Lred\otimes \sY)$ with $\sX,\sY\in \MCor$, there exists $f \in \MCor(\sX,\sY)$ such that $\phi$ and $id_{\Lred}\otimes f$ have the same image in 
$\Hom_{\ulMPST}(\Lred\otimes \sX , \Lred\otCIsp \sY)$.
\end{prop}

First we deduce Theorem \ref{thm;weakcancel} follows from Proposition \ref{prop;almostsurj}. By Proposition \ref{prop;conj2splitting} it suffices to show the surjectivity of the map \eqref{eq;iotapreF} $\iota_F^{sp}$.
Proposition \ref{prop;almostsurj} implies that the following composition 
\[h_0^\bcube(\sY)\to \gamma(\Lred \otCI \sY)\to 
\gamma(\Lred\otCIsp \sY)\simeq \gamma(\Lred\otCIsp h_0^\bcube(\sY))\]
is surjective. Since the last object is semi-pure,
it factors through $h_0^\bcube(\sY)^{sp}$, proving the desired surjectivity for 
$F=h_0^\bcube(\sY)$. 

For a general $F\in \CIt$ consider a surjection \[
q:\bigoplus_{\sY\to F} h_0^\bcube (\sY)\to F
\]
which gives a commutative diagram\[\xymatrix{
\bigoplus h_0^\bcube(\sY)^{sp}  \ar[r]^{\hskip -20pt\oplus\iota_{\sY}^{sp}} \ar[d]^{q^{sp}}
& \bigoplus \gamma( \Lnredd{1}\otCIsp  \sY) \ar[d]   \\
F^{sp}  \ar[r]^{\hskip -20pt\iota_F^{sp}}& \gamma( \Lnredd{1}\otCIsp  F)\\ }\]
where the top arrow is surjective 
and the vertical arrows are surjective since representable presheaves are projective objects of $\ulMPST$ by Yoneda's lemma and the functors $(\_)^{sp}$ and $ \Lred\otCI\_$ commute with direct sums and preserves surjective maps. 
This proves the desired surjectivity of $\iota_F$. 

\medbreak

\def\End{\mathrm{End}}
 
The proof of Proposition \ref{prop;almostsurj} requires a construction analogous to the one in \cite{gramot}. For a variable $T$ over $k$ and for $i\geq 1$, we put 
\[ \bcube^{(i)}_{T}=(\P^1_T,i(0+\infty))\]
where $\P^1_T$ is the compactification of $\G_{m,T}=\Spec k[T,T^{-1}]$.
We also put (cf. \eqref{eq;cuberedgm}):
\[ \bcube^{(i)}_{T,red} =\Ker\big(\Ztr(\bcube^{(i)}_{T}) \rmapo{pr} \Z=\Ztr(\Spec k,\emptyset)\big)
\in \MPST,\]
where $pr:\P^1_T \to\Spec k$ is the projection. 
Let $e$ be the composite of $pr$ and $i_1: \Z\to \Ztr(\LT)$ induced by $1\in \P^1_T$. Then $e$ is an idempotent of $\End_{\MPST}(\LT)$ and 
$id-e\in \End_{\MPST}(\LT)$, with $id$ denoting the identity on $\bcube^{(i)}_T$, is a splitting of $\bcube^{(i)}_{T,red}\to \bcube^{(i)}_{T}$. Thus, we get a direct sum decomposition in $\MPST$ (cf. \eqref{eq;cuberedgm}):
\[ \bcube^{(i)}_T= \bcube^{(i)}_{T,red} \oplus \Z \qwith \;\;
\bcube^{(i)}_{T,red}=(id-e)\bcube^{(i)}_T.\]
For $F\in \ulMPST$ and integers $i_1,\dots,i_n\geq 1$, let
\[\pi: \Hom_{\ulMPST}(\bcube^{(i_1)}_T\otimes\cdots\otimes \bcube^{(i_n)}_T,F)
\to
\Hom_{\ulMPST}(\bcube^{(i_1)}_{T,red}\otimes\cdots\otimes \bcube^{(i_n)}_{T,red},F) 
\]
be the projection induced by the above decomposition.

\medbreak
For $X\in \Sm$ and $a\in \Gamma(X,\sO^\times)$, let 
$[a] \in \Cor(X,\A^1-\{0\})$ be the map given by $z\to a$, where
$\A^1=\Spec k[z]$.

\begin{lemma}\label{TUTplusU}
\begin{itemize}
\item[(1)]
The correspondences 
\[[T], [U], [TU], [1] \in \Cor((\A^1_T-\{0\})\times (\A^1_U-\{0\}),(\A^1-\{0\}))\]
lie in $\MCor(\LT\otimes \LU, \cubegm)$.
Moreover we have
\[[T]+[U] -[TU] - [1]=0 \in  \tMCor(\LT\otimes \LU, h_0^{\bcube}(\cubegm)).\]
\item[(2)] The correspondences 
\[[-T], [-U], [-TU], [-1] \in \Cor((\A^1_T-\{0\})\times (\A^1_U-\{0\}),(\A^1-\{0\}))\]
lie in $\MCor(\LT\otimes \LU, \cubegm)$.
Moreover we have
\[[-T]+[-U] -[-TU] - [-1]=0 \in  \tMCor(\LT\otimes \LU, h_0^{\bcube}(\cubegm)).\]
\end{itemize}
\end{lemma}
  \begin{proof}
  The first assertion of (1) follows from the fact 
\[   [T]=\mu(id\otimes [1]),\qquad [U]=\mu(id\otimes [1]),\qquad [TU]=\mu
  \]
  where $\mu:(\A_T^1-\{0\})\times(\A_U^1-\{0\})\to (\A_W^1-\{0\})$ is the multiplication $W=TU$, which lies in $\MCor(\LT\otimes \LU, \LW)$ by \cite[Claim 1.21]{shuji}.

  To show the second assertion of (1), consider as in \cite[p.142]{blochkato} the finite correspondence $Z$ given by the following algebraic subset:
  \begin{multline}\label{homoptyZ}
  \{V^2 -\big (W(T+U) + (1-W)(TU + 1)\big)V + TU=0\} \\
  \in \mathbf{Cor}((\AA_T^1-\{0\})\times (\AA_U^1-\{0\})\times \AA^1_W,
\AA_V^1-\{0\})    
  \end{multline}
  Let \[
  i_0,i_1: (\AA_T^1-0)\times (\AA_U^1-0)\times (\AA_V^1-0)\to (\AA_T^1-0)\times(\AA_U^1-0)\times \AA_W^1\times(\AA_V^1-0)
  \]
  be the maps induced by the inclusion of $0_W$ and $1_W$ in $\AA^1_W$. It is clear that $(i_0^*-i_1^*)(Z)=([TU]+[1])-([T]+[U])$ since
  \begin{align*}
  &V^2 - (TU + 1)V + TU=(V-TU)(V-1),\\ 
  &V^2 - (T+U)V + TU=(V-T)(V-U)
  \end{align*}
We need to check that $Z$ lies in 
$\MCor(\LT\otimes \LU\otimes\bcube_W, \LV)$.
Consider the compactification $(\P^1)^{\times 4}$ of $\A^1_T\times\A^1_U\times\A^1_W\times\A^1_V$ given coordinates with the usual convention $[0:1]=\infty$ and $[1:0]=0$:
\[  ([T_0:T_\infty],[U_0:U_\infty],[W_0:W_\infty],[V_0:V_\infty]).\]
Then the closure of $Z$ is the hypersurface given by the following polyhomogeneous polynomial:
  \begin{multline*}
      T_0U_0W_0V_\infty^2 - \big(W_\infty(T_0U_\infty+T_\infty U_0) + (W_0-W_\infty)(T_\infty U_\infty + T_0U_0)\big)V_\infty V_0\\
      + T_\infty U_\infty W_0V_0^2.
  \end{multline*}
  We have to check that it satisfies the modulus condition: letting
 \[  \phi: \overline{Z}\to (\PP^1)^{\times 4}   \]
  be the inclusion and letting
\[   D_1=(\{0\} + \{\infty\})\times\PP^1_U\times\PP^1_W\times \PP^1_V + 
\PP^1_T\times(\{0\} + \{\infty\})\times\PP^1_W\times\PP^1_V + 
\PP^1_T\times \PP^1_U\times \{\infty\}\times\PP^1_V,\]
\[   D_2=\PP^1_T\times\PP^1_U\times\PP^1_W\times(\{0\} + \{\infty\}),\] 
  we have to check the following inequality:
\begin{equation}\label{MCD}
  \phi^*(D_1) \geq \phi^*(D_2).  
\end{equation}
  Consider the Zariski cover of $(\PP^1)^{\times 4}$ given by:
\[  \Bigl\{\mathcal{U}_{\alpha,\beta,\gamma,\delta} = (\PP^1-\alpha)\times(\PP^1-\beta)\times(\PP^1-\gamma)\times(\PP^1-\delta),\ \alpha,\beta,\gamma,\delta\in \{0,\infty\}\Bigr\}.
  \]
Define $t_\alpha=T_\infty/T_0$ if $\alpha=\infty$ and $t_\alpha=T_0/T_\infty$ if $\alpha=0$ and $u_\beta$, $w_\gamma$, $v_\delta$ similarly. Then
  \[ \mathcal{U}_{\alpha,\beta,\gamma,\delta}=
\Spec(k[t_\alpha,u_\beta,w_\gamma,v_\delta]).  \]
  On this cover, the Cartier divisors $D_1$ and $D_2$ are given by the following system of local equations:
\[ D_1=\Bigl\{(\mathcal{U}_{\alpha,\beta,0,\delta},t_\alpha u_\beta w_0),
(\mathcal{U}_{\alpha,\beta,\infty, \delta},t_\alpha u_\beta)\Bigr\}
\qquad D_2=\Bigl\{(\mathcal{U}_{\alpha,\beta,\gamma,\delta},v_\delta)\Bigr\}
  \]
  The equation of $\ol{Z}$ on
$(\P^1)^{\times 4} - \{0\}$ is of the form\[
T_0U_0W_0 - v_0F,\qquad\textrm{for some }F\in k[v_0][U_0,U_{\infty},\ldots]\]
Hence \eqref{MCD} is satisfied on $\mathcal{U}_{\alpha,\beta,\gamma, 0}$ if $\alpha = 0$ or $\beta = 0$ or $\gamma = 0$. Furthermore $\ol{Z}\cap \mathcal{U}_{\infty,\infty,\infty, 0}\cap D_2=\emptyset$. Similarly, the equation of $\ol{Z}$ on
$(\P^1)^{\times 4} - \{\infty\}$ is of the form \[
T_{\infty}U_{\infty}W_0 - v_{\infty}G,\qquad\textrm{for some }G\in k[v_{\infty}][U_0,U_{\infty},\ldots]
\] 
Hence \eqref{MCD} is satisfied on $\mathcal{U}_{\alpha,\beta,\gamma, \infty}$ if $\alpha = \infty$ or $\beta = \infty$ or $\gamma = 0$. Furthermore $\ol{Z}\cap \mathcal{U}_{0,0,\infty,\infty}\cap D_2=\emptyset$.
\medbreak
(2) is proved by the same argument using the following correspondence instead of 
\eqref{homoptyZ}: 
  \begin{multline*}
  \{V^2 +\big (W(T+U) + (1-W)(TU + 1)\big)V + TU=0\} \\
  \in \mathbf{Cor}((\AA_T^1-\{0\})\times (\AA_U^1-\{0\})\times \AA^1_W,
\AA_V^1-\{0\})   . 
  \end{multline*}
\end{proof}

\begin{cor}\label{muzero}
  $\pi([TU])= 0\in \tMCor( \LTred \otimes \LUred ,h_0^{\bcube}(\cubegm))$.
 \end{cor}
 \begin{proof}
This follows from Lemma \ref{TUTplusU} since
\begin{multline*}
 [TU]\circ \big((id-e)\otimes(id-e)\big) = [TU] -[TU]\circ(1\otimes e) -
[TU]\circ(e\otimes 1) + [TU]\circ(e\otimes e)\\
=[TU]-[T]-[U] + [1] \text{ in }
\Hom_{\MPST}(\LT\otimes \LU,\cubegm) .\end{multline*}
 \end{proof}
\medbreak

For $X\in \Sm$ and $a,b\in \Gamma(X,\sO^\times)$, let 
\[[a,b] \in \Cor(X,(\A^1-\{0\})\otimes (\A^1-\{0\}))\]
 be the map given by $z\to a,\; w\to b$, where $z$ (resp. $w$) is the standard coordinate of the first (resp. second) $\A^1$.

  \begin{cor}  \label{highdim}
In $\tMCor(\LT\otimes \LU \otimes \LV, h_0^{\bcube}(\cubegm\otimes \cubegm))$,
we have:
\begin{multline*}
[T,V] +[U,V]-[TU,V]- [1,V]=[-T,V] +[-U,V]-[-TU,V]- [-1,V]=0.
\end{multline*}
 \end{cor}
  \begin{proof} 
This follows from Lemma \ref{TUTplusU} noting the endofunctor $\_\otimes \cubegm$ on
$\MPST$ is additive and 
$h_0^{\bcube}(\cubegm \otimes \cubegm)$ is a quotient of $h_0^{\bcube}(\cubegm)\otimes \cubegm$.
  \end{proof}


\begin{prop}  \label{inversion}
The correspondences
\[[U,T],\; [T^{-1},U] \in \Cor((\A^1_T-\{0\})\times (\A^1_U-\{0\}), (\A^1-\{0\})\times (\A^1-\{0\}))\]
lie in $\MCor(\LT\otimes \LU, \cubegm\otimes\cubegm)$.
Moreover, the element 
\[\pi([U,T])-\pi([T^{-1},U]) \in \tMCor(\LTred\otimes\LUred, h_0^\bcube(\cubegm\otimes\cubegm))\]
lies in the kernel of the map
\begin{multline*}
    \tMCor(\LTred\otimes\LUred, h_0^\bcube(\cubegm\otimes\cubegm))\to \\
    \tMCor(\LLTred\otimes\LLUred, h_0^\bcube(\cubegm\otimes\cubegm))
\end{multline*}
\begin{proof}(see \cite[Corollary 9]{gramot})
The first assertion is easily checked. To show the second, consider the map in $\MCor$:
\[  \LLS  \to \LT \otimes \LU \;;\; T \to S,\; U \to S^{-1}.  \]
  Composing this with the correspondences of Lemma \ref{TUTplusU}(1), we get
  
\[ [S]+[S^{-1}]- 2[1] =0 \in  \tMCor (\LLSred, h_0^{\bcube}(\cubegm)).  \]
Noting $\pi([1])=(id-e)\circ[1]=0$, we get
\[  \pi([S]+[S^{-1}]) =0 \in  \tMCor (\LLSred, h_0^{\bcube}(\cubegm)).\]

This implies
\begin{equation}\label{SVSinvV}
 \pi([S,V] +[S^{-1},V]) =0 \in 
\tMCor (\LLSred\otimes \LVred, h_0^{\bcube}(\cubegm \otimes \cubegm)).
\end{equation}
again noting that the endofunctor $\_\otimes \LV$ on $\MCor$ is additive and 
$h_0^{\bcube}(\cubegm \otimes \cubegm)$ is a quotient of $h_0^{\bcube}(\cubegm)\otimes \cubegm$.

  On the other hand, by tensoring the correspondence of Corollary \ref{muzero} with another copy of itself we get
\begin{multline}
    \label{eq;TUVW}
  \pi([TU, VW])=0\\ 
  \text{ in }\tMCor( \LTred\otimes\LUred\otimes \LVred\otimes \LWred,h_0^{\bcube}(\cubegm\otimes\cubegm) ).  
\end{multline}
There is a map in $\MCor$:
\begin{multline*}
\LL_{S_1}\otimes \LL_{S_2} \to \LT\otimes \LU \otimes \LV\otimes \cubegm_W\;;\;\\
T\to S_1,\; U\to S_2,\; V\to -S_1,\, W\to S_2,
\end{multline*}
which induces an element of  
\[\tMCor(\LLred{S_1}\otimes \LLred{S_2}, \LTred\otimes \LUred\otimes \LVred\otimes \LWred).\]
Composing this with \eqref{eq;TUVW} and changing variables $(S_1,S_2)$ to $(T,U)$,
we get 
\begin{equation}\label{TUTU0}
\pi([TU,-TU]) =0\in  \tMCor(\LLTred\otimes\LLUred, h_0^{\bcube}(\cubegm\otimes\cubegm)). 
\end{equation}
We make the following claim:
\begin{claim}
In $\tMCor(\LTred\otimes\LUred, h_0^{\bcube}(\cubegm\otimes\cubegm))$,
we have 
\begin{equation}\label{claim:cor1}  \pi([TU,-TU])= \pi([T,-TU])+\pi([U,-TU]),
\end{equation}
\begin{equation}\label{claim:cor2}
	  \pi([T,-TU])= \pi([T,U]),
\end{equation}
\begin{equation}\label{claim:cor3} \pi([U,-TU])= \pi([U,T]).
\end{equation}
\end{claim}

Indeed, composing the first correspondence of Corollary \ref{highdim} with the map in $\MCor$:
\eq{TUTUV}{\LT\otimes \LU \to \LT\otimes \LU \otimes \LV}
given by $V\to -TU$ which is admissible by \cite[Claim 1.21]{shuji}, we get
\begin{multline*}
    [TU,-TU] + [1,-TU] - [T,-TU]-[U,-TU]=0\\
    \text{ in }\tMCor(\LT\otimes \LU, h_0^{\bcube}(\cubegm\otimes\cubegm)).
\end{multline*} 
Then \eqref{claim:cor1} follows from the equality:
\[ \pi([1,-TU])= 0\in \tMCor( \LTred \otimes \LUred ,h_0^{\bcube}(\cubegm\otimes\cubegm)).\]

Indeed, we have
\begin{multline*}
[1,-TU]\circ ((id-e)\otimes (id-e))= \\
 [1,-TU]- [1,-TU]\circ (id\otimes e) -
[1,-TU]\circ (e\otimes id) + [1,-TU]\circ (e\otimes e)\\
= [1,-TU]-[1,-T]-[1,-U]+[1,-1]\overset{(*)}{=}0\\\text{ in }\tMCor(\LT\otimes\LU, \cubegm\otimes\cubegm),
\end{multline*}
where the equality $(*)$ above follows from Corollary \ref{highdim}. 
Then \eqref{claim:cor2} and \eqref{claim:cor3} follow from Corollary \ref{highdim} by an analogous argument considering the maps \eqref{TUTUV} given by 
$V\to T, T\to -T$ and $V\to U, U\to -U$ respectively, and noticing that:
\begin{multline*}
[T,-T]\circ ((id-e)\otimes (id-e))= \\
 [T,-T]- [T,-T]\circ (id\otimes e) -
[T,-T]\circ (e\otimes id) + [T,-T]\circ (e\otimes e)\\
= [T,-T]-[T,-T]-[1,-1]+[1,-1]=0,\\
\end{multline*}
and similarly for $[U,-U]$. This completes the proof of the claim.

\medbreak

By the above claim, \eqref{TUTU0} implies 
\begin{equation}\label{TUUT}
\pi[T,U]+\pi[U,T] =0\text{ in } \tMCor(\LLTred\otimes\LLUred, h_0^{\bcube}(\cubegm\otimes\cubegm)).
\end{equation}
Putting \eqref{SVSinvV} and \eqref{TUUT} together we conclude that
\[   \pi[T,U]-\pi[U^{-1},T] =0 \text{ in } \tMCor(\LLTred\otimes\LLUred, h_0^{\bcube}(\cubegm\otimes\cubegm)).   \]
This completes the proof of Proposition \ref{inversion}.
\end{proof}
\end{prop}

\bigskip
\def\talpha{\tilde{\alpha}}
\def\tbeta{\tilde{\beta}}
\def\phicube{\phi_{\bcube}}
\def\sigmacube{\sigma_{\bcube}}
\medbreak

Take $\sX,\sY\in \MCor$ and 
\[\phi\in \Hom_{\ulMPST}(\Lred\otimes \sX ,\Lred \otimes \sY)\]
It induces  
\[\phicube\in \Hom_{\ulMPST}(\Lred\otimes \sX ,\Lred \otCI \sY).\]
Let 
\[ \phi^* \in\Hom_{\ulMPST}(\sX\otimes \Lred ,\sY \otimes \Lred)\]
be obtained from $\phi$ by the obvious permutation.
 It induces
\[ \phicube^* \in\Hom_{\ulMPST}(\sX\otimes \Lred,\sY\otCI \Lred).\]
We then put
\begin{multline*}
\phi\otimes Id_{\Lred} \in 
\Hom_{\ulMPST}(\Lred\otimes \sX\otimes \Lred, \Lred \otimes\sY\otimes \Lred),
\end{multline*}
\begin{multline*}
Id_{\Lred} \otimes \phi^* \in 
\Hom_{\ulMPST}(\Lred\otimes \sX\otimes \Lred, \Lred  \otimes \sY\otimes \Lred),
\end{multline*}
which induce
\begin{multline*}
\phicube \otimes Id_{\Lred} \in \Hom_{\ulMPST}(\Lred\otimes \sX\otimes  \Lred, \Lred \otCI \sY\otCI \Lred ),
\end{multline*}
\begin{multline*}
Id_{\Lred} \otimes \phicube^*\in \Hom_{\ulMPST}(\Lred\otimes \sX\otimes \Lred,\Lred \otCI \sY \otCI \Lred).
\end{multline*}
For $\mathcal{M}\in \MCor$, let\[
\sigma_{\mathcal{M}}\colon \Lred\otimes \mathcal{M}\otimes \Lred \to \Lred\otimes\mathcal{M}\otimes \Lred
\]
be the permutation of the two copies of $\Lred$. We have
\[ \phi\otimes Id_{\Lred} =(\sigma_\sY) \circ (Id_{\Lred} \otimes \phi^* )\circ (\sigma_\sX).\]
Let $T$ be the standard coordinate on $\A^1$ and let
\begin{equation}\label{iotadef}
\iota:\Lred\to \Lred\end{equation}
be the map given by $T\to T^{-1}$. For all $\mathcal{M}\in \MCor$ let
\[\sigma'_{\mathcal{M}}=  \sigma_{\mathcal{M}}- Id_{\Lred\otimes \mathcal{M}}\otimes \iota\colon \Lred\otimes \mathcal{M}\otimes \Lred \to \Lred\otimes\mathcal{M}\otimes \Lred.\]
We can write  
\[\phi\otimes Id_{\Lred} =Id_{\Lred} \otimes \phi^* +  (\sigma'_\sY) \circ p + q\circ (\sigma'_\sX),\]
for some
\[p,q\in \Hom_{\ulMPST}(\Lred\otimes \sX\otimes \Lred,\Lred \otimes\sY \otimes \Lred).\]
Put
\[ \Gamma_\sX= \Lred \otCI \sX\otCI \Lred\qquad 
\Gamma_\sY= \Lred \otCI \sY \otCI \Lred .\]
Hence we can write
\begin{equation}\label{eq;phiphi*}
\phicube\otimes Id_{\Lred} =Id_{\Lred} \otimes \phicube^* + 
 \sigma'_{\bcube,\sY} \circ p + q_{\bcube} \circ \sigma'_{\bcube,X},
\end{equation}
where 
\[ \sigma'_{\bcube,\sY}: \Lred\otimes \sY \otimes \Lred\to \Gamma_\sY\]
\[ \sigma'_{\bcube,\sX}: \Lred \otimes \sX \otimes \Lred \to 
\Gamma_\sX\]
\[q_{\bcube}: \Gamma_\sX
\to \Gamma_\sY \]
are induced by $\sigma'_\sY$, $\sigma'_\sX$ and $q$ respectively. For an integer $n>0$ let $\sX^{(n)}:=(X,nD)$ if $\sX=(X,D)$.
Then we consider the map\[
\Hom_{\ulMPST}(\Lred\otimes \sX \otimes \Lred,\Gamma_\sY)\xrightarrow{\beta^*_n}
\Hom_{\ulMPST}(\Lnredd{n}\otimes \sX^{(n)}\otimes \Lnredd{n},\Gamma_\sY) \]
induced by the natural map
$\beta_n: \Lnredd{n}\otimes\sX^{(n)}\otimes \Lnredd{n}\to\Lred\otimes\sX\otimes \Lred$.

\begin{claim}\label{claim1;proofmainthm}
There is $N\geq 2$ such that for all $n\geq N$ the maps $\sigma'_{\bcube,\sY} \circ p$ and $q_{\bcube} \circ \sigma'_{\bcube,\sX}$ lie in the kernel of 
\[\Hom_{\ulMPST}(\Lred\otimes \sX \otimes \Lred,\Gamma_\sY)\xrightarrow{\beta^*_n}
\Hom_{\ulMPST}(\Lnredd{n}\otimes \sX^{(n)}\otimes \Lnredd{n},\Gamma_\sY) \]
\end{claim}

\begin{proof}

By Proposition \ref{inversion}, the composite map
\[ \LLredd\otimes \LLredd \rmapo{\beta_2} \Lred\otimes \Lred\rmapo{\sigma'} \Lred\otimes\Lred \to  h_0^\bcube(\Lred)  \otCI h_0^\bcube(\Lred)\]
is zero, where
$\sigma'=  \sigma- Id_{\Lred}\otimes \iota$
with $\sigma$ the permutation of the two copies of $\Lred$ and $\iota$ from \eqref{iotadef}. This immediately implies the claim for $q_{\bcube} \circ \sigma'_{\bcube,\sX}$. We now show the claim for $\sigma'_{\bcube,\sY} \circ p$. Choose such an integer $N$ that for all $n\geq N$ there is map 
\[p^{(n)}\in \Hom_{\ulMPST}(\Lnredd{n}\otimes \sX^{(n)}\otimes \Lnredd{n},\LLredd\otimes \sY^{(2)}\otimes \LLredd)\]
induced by $p$. For $\mathcal{M},\mathcal{N}\in \MCor$, write
\[ \Lambda_{\mathcal{M},\mathcal{N}}= \Hom_{\ulMPST}(\Lred\otimes \mathcal{M}\otimes \Lred, \Lred \otCI \mathcal{N} \otCI  \Lred),\]
\[ \Lambda^{(n)}_{\mathcal{M},\mathcal{N}}= \Hom_{\ulMPST}(\Lnredd{n}\otimes \mathcal{M}^{(n)}\otimes \Lnredd{n},\Lred \otCI \mathcal{N} \otCI \Lred).\] Then we have a commutative diagram for $n\geq N$
\begin{equation}\label{CDgammabeta}
	\xymatrix{
		\Lambda_{\sY,\sY} \ar[r]^{p^*}\ar[d]^{\beta_2^*} & \Lambda_{\sX,\sY} \ar[d]^{\beta_n^*}\\
		\Lambda^{(2)}_{\sY,\sY} \ar[r]^{(p^{(n)})^*} & \Lambda^{(n)}_{\sX,\sY} ,\\}
\end{equation}
The claim for $\sigma'_{\bcube,\sY} \circ p$ follows from this.

\end{proof}

We now complete the proof of Proposition \ref{prop;almostsurj}.
We consider the commutative diagram
\begin{small}
\[\xymatrix@C=1em{
\Hom_{\ulMPST}(\Lred\otimes \sX\otimes \Lred,\Lred\otimes \sY \otimes \Lred)  \ar[r]^-{\rho_1} \ar[d]^{\beta_n^*} & 
\Hom_{\ulMPST}(\sX\otimes\Lred,\sY\otCI \Lred)
\ar[d]^{\beta_n^*} \\
\Hom_{\ulMPST}(\Lnredd{n}\otimes \sX^{(n)}\otimes \Lnredd{n},\Lred\otimes \sY \otimes \Lred) \ar[r]^-{\rho_n}  & 
\Hom_{\ulMPST}(\sX^{(n)}\otimes\Lnredd{n},\sY\otCI \Lred)\\}\]
\end{small}

where the horizontal maps come from \eqref{rhoa} replacing $\sY$ with $\sY\otimes \Lred$.
By Lemma \ref{lem;voeLem4.1}(3) and (2) we have
\[ \rho_1(\phi\otimes id_{\Lred}) = \rho(\phi)\otimes Id_{\Lred}\qaq
\rho_1(Id_{\Lred} \otimes \phi^*) = \phicube^*,\;\text{ where}\]
\begin{equation}\label{eq;rhoX}
\rho: \Hom_{\ulMPST}(\Lred\otimes \sX ,\Lred\otimes \sY)
\to \Hom_{\ulMPST}(\sX,  h_0^\bcube(\sY))
\end{equation}
is the map from \eqref{rhoa}. 
In view of the diagram, \eqref{eq;phiphi*} and Claim \ref{claim1;proofmainthm} imply that there is $n\gg0$ such that
$\beta_n^*(\phicube^*-\rho(\phicube)\otimes Id_{\Lred})=0$ so that
\begin{equation}\label{eq;rhophiId}
\beta_n^*(\phicube-Id_{\Lred}\otimes \rho(\phicube))=0\in
\Hom_{\ulMPST}(\Lnredd{n}\otimes \sX^{(n)} , \Lred\otCI \sY).
\end{equation}

Consider the commutative diagram \[\xymatrix{
\Hom_{\ulMPST}(\Lred\otimes \sX ,\Lred\otCI \sY)  \ar[r] \ar[d]^{\beta^*_n} & 
\Hom_{\ulMPST}(\Lred\otimes \sX , \Lred\otCIsp \sY) \ar[d]^{\beta^*_n} \\
\Hom_{\ulMPST}(\Lnredd{n}\otimes \sX^{(n)}, \Lred\otCI \sY)  \ar[r]
& \Hom_{\ulMPST}(\Lnredd{n}\otimes \sX^{(n)} , \Lred\otCIsp \sY)\\ }\]
The two horizontal maps are surjective since representable presheaves are projective objects of $\ulMPST$ and 
$\Lred\otCI \sY\to \Lred\otCIsp \sY$ is surjective.
The map $\beta^*_n$ on the right hand side is injective since $\Lred\otCIsp \sY$ is semi-pure. Hence Proposition \ref{prop;almostsurj} follows from \eqref{eq;rhophiId}.

\section{Implications on reciprocity sheaves}\label{implicationRSC}

Let $\RSC_\Nis$ be the category of reciprocity sheaves (see \S\ref{recollection}
\eqref{RSC}). 
Recall that for simplicity, we denote for all $F\in \RSC_\Nis$
(cf. \S\ref{recollection} \eqref{omegaCI})
\[
\tF:=\ulomega^{\CI} F\in \CItsp_\Nis.
\]
By \cite{rsy} there is a \emph{lax} monoidal structure on $\RSC_\Nis$ given by
(cf. Proposition \ref{prop;otCI})
\[\bigl(F, G\bigr)_{\RSC_{\Nis}}:=\ulomega_!(\tF\otCINissp \tG).\]
Following \cite[5.21]{rsy}, we define 
\eq{otRSCdef}{F\langle 0\rangle:=F,\qquad F\langle n\rangle:= \bigl(F\langle n-1\rangle,\G_m\bigr)_{\RSC_{\Nis}}\text{ for }n\geq 1.}
By Corollary \ref{cor;gammaGm}(1), we have (cf. \eqref{twist(1)})
\eq{twistulomega}{F\langle n\rangle\cong \ulomega_!(\widetilde{F\langle n-1\rangle}(1)).}
By recursiveness of the definition we have 
\begin{equation}\label{twistrecursive}
\bigl(F\langle n\rangle\bigr)\langle m\rangle\cong F\langle n+m\rangle.
\end{equation}
By \cite[Prop. 5.6 and Cor. 5.22]{rsy}, we have isomorphisms
\begin{equation}\label{KMGatensor}
\ulomega_!((\ulomega^*\G_m)^{\otCINissp n})\cong \Z\langle n\rangle \cong \mathcal{K}_n^M,\quad
\G_a\langle n\rangle\cong \Omega^n \;\text{  if }\ ch(k)=0,
\end{equation}
where the second isomorphism is defined as follows: for an affine $X=\Spec A\in \Sm$, the composite map
\begin{equation}\label{GaGmOmega} 
\G_a(A)\otimes_\bZ \G_m(A)^{\otimes_\Z n} \to 
(\G_a\otimes_{\NST} \G_m^{\otimes_{\NST} n} )(A) \to 
\G_a\langle n\rangle(A) \rmapo{\eqref{KMGatensor}}  \Omega^n_A
\end{equation}
sends $a\otimes f_1\otimes\cdots\otimes f_n$ with $a\in A$ and $f_i\in A^\times$ to $a\dlog f_1\wedge\cdots \wedge\dlog f_n$.

By \cite[5.21 (4)]{rsy}, there is a natural surjective map for $F\in \RSC_\Nis$
\begin{equation}\label{FGm}
F\otimes_{\NST} \sK_n^M \to F\langle n\rangle.
\end{equation}

\begin{lemma}\label{lem;lambdadef}
The map \eqref{FGm} factors through a natural surjective map 
\eq{tFtensor}{ \ulomega_!(\tF\otCINissp (\ulomega^*\G_m)^{\otCINissp n})
\to F\langle n \rangle.}
\end{lemma}
\begin{proof}
By \cite[(5.21.1)]{rsy}, there is a natural surjective map
\begin{equation}\label{FGm2}
\ulomega_!\ulaNis h_0^{\bcube}(\tF\otimes_{\ulMPST} (\ulomega^*\G_m)^{\otimes_{\ulMPST}n}) \to F\langle n\rangle.
\end{equation}
By Lemma \ref{lem2;hcubesp} (ii) and (iii), we have a natural isomorphism
\[\ulomega_!\ulaNis h_0^{\bcube}(\tF\otimes_{\ulMPST} (\ulomega^*\G_m)^{\otimes_{\ulMPST}n})  \simeq \ulomega_!(\tF\otCINissp (\ulomega^*\G_m)^{\otCINissp n}).\]
Hence \eqref{FGm2} induces \eqref{tFtensor}.
We have a surjective map
\begin{multline*}
F\otimes_{\PST} \sK_n^M\overset{\eqref{KMGatensor}}\simeq 
\ulomega_!\tF \otimes_{\PST} \ulomega_!((\ulomega^*\G_m)^{\otCINissp n})\simeq\\
\ulomega_!(\tF \otimes_{\ulMPST} ((\ulomega^*\G_m)^{\otCINissp n}))
\to \ulomega_!(\tF\otCINissp (\ulomega^*\G_m)^{\otCINissp n}).
\end{multline*}
where the second isomorphism comes from the monoidality of $\ulomega_!$ (cf. \S1\eqref{tensor}). By the adjunction from \eqref{aVNis}, this induces a surjective map
\eq{FGm3}{F\otimes_{\NST} \sK_n^M = \aVNis(F\otimes_{\PST} \sK_n^M) \to 
\ulomega_!(\tF\otCINissp (\ulomega^*\G_m)^{\otCINissp n}) .}
By the construction of \eqref{FGm2}, it is straightforward to check 
that \eqref{FGm} is the composite \eqref{tFtensor} and \eqref{FGm3}. 
This completes the proof of the lemma.
\end{proof}

We have a map natural in $X\in \Sm$:
\begin{multline}\label{iotaFmX}
F(X)=\Hom_{\PST}(\Ztr(X),F) \xrightarrow{\_\otimes id_{\sK_n^M}} 
\Hom_{\PST}(\Ztr(X)\otimes_{\NST}\sK_n^M,F\otimes_{\NST}\sK_n^M)\\
\to \Hom_{\PST}(\Ztr(X)\otimes_{\NST}\sK_n^M,F\langle n\rangle),
\end{multline}
where the last map is induced by \eqref{FGm}.
Thus we get a map
\begin{equation}\label{iotaFGm}
\lambda_F^n : F \to \uHom_{\PST}(\sK_n^M,F\langle n\rangle).
\end{equation}

\begin{thm}\label{cor;cancelmilnorK}
		For $F\in \RSC_\Nis$, the map $\lambda_F^n$  is an isomorphism.
\end{thm}

The proof will be given later. First we prove the following.

\begin{prop}\label{prop;cancelRSC}
The map $\lambda_F^n$ is an isomorphism for $n=1$.
\end{prop}
\begin{proof}
Note $\sK_1^M=\G_m$ and that for $F_1,G_1,F_2,G_2\in \ulMPST$ and  maps $f: F_1\to F_2$, $g:G_1\to G_2$, the diagram
\[\xymatrix{
\ulomega_!F_1\otimes_{\PST} \ulomega_!G_1 \ar[r]^-{\ulomega_! f\otimes \ulomega_! g}\ar[d]^{\simeq}&\ulomega_!F_2\otimes_{\PST} \ulomega_!G_2\ar[d]^{\simeq}\\
\ulomega_!(F_1\otimes_{\ulMPST} G_1) \ar[r]^-{\ulomega_! (f\otimes g)}&\ulomega_!(F_2\otimes_{\ulMPST} G_2),}\]
commutes, where the vertical isomorphisms follow from the monoidality of $\ulomega_!$.
Thus, by Lemma \ref{lem;lambdadef}, \eqref{iotaFmX} with $n=1$ coincides with the composite map:
\eq{iotaFmX2}{
\begin{aligned}
F(X)&=\ulomega_!\tF(X)\xrightarrow{\ulomega_!(\_\otimes id_{\ulomega^*\G_m})(X)}\ulomega_!\uHom_{\ulMPST}(\ulomega^*\G_m, \tF\otCINissp\ulomega^*\G_m)(X)\\
&\simeq \Hom_{\ulMPST}(\ulomega^*\G_m,\uHom_{\ulMPST}(\Ztr(X,\emptyset), \tF\otCINissp\ulomega^*\G_m))\\
&\overset{(*1)}{\simeq} \Hom_{\ulMPST}(\ulomega^*\G_m,\ulomegaCI\ulomega_!\uHom_{\ulMPST}(\Ztr(X,\emptyset), \tF\otCINissp\ulomega^*\G_m))\\
&\overset{(*2)}{\simeq} \Hom_{\PST}(\G_m,\ulomega_!\uHom_{\ulMPST}(\Ztr(X,\emptyset), \tF\otCINissp\ulomega^*\G_m))\\
&\overset{(*3)}{\simeq} \Hom_{\PST}(\G_m,\uHom(\Ztr(X), \ulomega_!(\tF\otCINissp\ulomega^*\G_m))\\
&\overset{(*4)}{\simeq} \uHom_{\PST}(\G_m,F\langle 1\rangle)(X)
\end{aligned}}
where $(*1)$ is induced by the injective unit map 
$G\to \ulomegaCI\ulomega_! G$ ($G\in \CItspNis$) for the adjunction \eqref{omegaCIadjointNis} and it is an isomorphism by Corollary \ref{cor;gammadoesntseeconductors} and the fact that $\uHom_{\ulMPST}(\Ztr(X,\emptyset),\tF\otCINissp\omega^*\G_m)\in \CItspNis$,
$(*2)$ is given by the fully faithfulness of $\ulomegaCI$ and 
$\ulomegaCI \G_m =\ulomega^* \G_m$ by \cite[Lem .2.3.1]{ksyII},
 $(*3)$ follows from Lemma \ref{lem;omegauHom}, and $(*4)$ holds by the definition \eqref{otRSCdef}.
 
This gives a commutative diagram
\eq{iotalambda}{\xymatrix{
F \ar[r]^-{\lambda_F^1} \ar[d]^{\simeq} 
& \uHom_{\PST}(\G_m,F\langle 1\rangle)\\
\ulomega_!\tF \ar[r]^-{\ulomega_!\iota^1_{\tF}} 
& \ulomega_!\uHom_{\ulMPST}(\ulomega^*\G_m,\tF\otCIspNis \Lred)\ar[u]^{\simeq},}}
where $\iota_{\tF}^1=(\_\otimes id_{\ulomega^*\G_m})$ is an isomorphism from Corollary \ref{cor3;weakcancel}
(using Corollary \ref{cor;gammaGm}).
This proves the proposition. 

\end{proof}

For $F,G\in \RSC_\Nis$ let
\begin{equation}\label{eq;cancelRSC}
\iota_{F,G}: \Hom_{\PST}(F,G) \to \Hom_{\PST}(F\langle 1 \rangle,G\langle 1 \rangle)
\end{equation}
be the composite map 
\begin{multline*}
 \Hom_{\PST}(F,G) \rmapo{\ulomega^{\CI}} \Hom_{\ulMPST}(\tF,\tG) 
\rmapo{-\otCINis\ulomega^*\G_m} \\
\Hom_{\ulMPST}(\tF\otCINis \ulomega^*\G_m,\tG\otCINis \ulomega^*\G_m)\rmapo{\ulomega_!} \Hom_{\PST}(F\langle 1 \rangle,G\langle 1 \rangle).
\end{multline*}

\begin{thm}\label{thm;closeenoughtocancel}
For $F,G\in \RSC_\Nis$, $\iota_{F,G}$ is an isomorphism.
\end{thm}
\begin{proof}
  We have isomorphisms (cf. \S1 \eqref{omegaCI})
  \begin{multline}\label{HomFaG1}
    \Hom_{\PST}(F\langle 1 \rangle,G\langle 1 \rangle)\\ 
=\Hom_{\PST}(\ulomega_!(\tF\otCINissp\ulomega^*\G_m), \ulomega_!(\tG\otCINissp\Lred))\\
  \cong  \Hom_{\ulMPST}(\tF\otCINissp \ulomega^*\G_m,\ulomega^{\CI}\ulomega_!(\tG\otCINissp\Lred))\\
 \cong  \Hom_{\ulMPST}(\tF\otimes_{\ulMPST}\ulomega^*\G_m,\ulomega^{\CI}\ulomega_!(\tG\otCINissp\Lred))\\
  \cong \Hom_{\ulMPST}(\tF,\uHom_{\ulMPST}(\ulomega^*\G_m,\ulomega^{\CI}\ulomega_!(\tG\otCINissp\Lred))),
  \end{multline}
where the first (resp. second) isomorphism follows from \eqref{omegaCIadjoint} 
(resp. the fact $\ulomega^{\CI}\ulomega_!\tau_!(\tG\otCINissp\Lred)\in\CItsp_\Nis$).
Note that for $H\in \CItsp$, the natural map
$H \to  \ulomega^{\CI}\ulomega_! H$ is injective. 

Hence we get injective maps
\begin{multline}\label{uHominjmaps}
  \Hom_{\ulMPST}(\tF,\uHom_{\ulMPST}(\ulomega^*\G_m,\tG\otCINissp\Lred))\\
\hookrightarrow\Hom_{\ulMPST}(\tF,\uHom_{\ulMPST}(\ulomega^*\G_m,\ulomega^{\CI}\ulomega_!(\tG\otCINissp\Lred)))\\
\hookrightarrow\Hom_{\ulMPST}(\tF,\ulomega^{\CI}\ulomega_! \uHom_{\ulMPST}(\ulomega^*\G_m,\ulomega^{\CI}\ulomega_!(\tG\otCINissp\Lred)))\\
\overset{(*1)}{\simeq} 
\Hom_{\ulMPST}(\tF,\ulomega^{\CI}\uHom_{\PST}(\G_m,\ulomega_!(\tG\otCINissp\Lred)))\\
\overset{(*2)}{\simeq} 
\Hom_{\ulMPST}(\tF,\ulomega^{\CI}\uHom_{\PST}(\G_m,G\langle 1\rangle)),\\
\end{multline}
where the isomorphism $(*1)$ comes from Proposition \ref{prop;omegagamma} and
$\ulomega_!\ulomega^{\CI}\simeq id$ (cf. \S1 \eqref{omegaCI}) and 
$(*2)$ follows from \eqref{twistulomega}.
These maps fit into a commutative diagram
  \[\xymatrix{
&  \ar[ld]^{\simeq}_{\alpha}\Hom_{\ulMPST}(\tF,\tG)  \\
\Hom_{\ulMPST}(\tF,\uHom_{\ulMPST}(\ulomega^*\G_m,\tG\otCINissp\Lred))\ar[d]^{\hookrightarrow}&\Hom_{\PST}(F,G) \ar[d]^{\iota_{F,G}}
\ar[u]^-{\simeq}_-{\ulomega^{\CI}} \ar@/^6.0pc/[dd]^{\ulomegaCI}_{\simeq} &   \\
\Hom_{\ulMPST}(\tF,\uHom_{\ulMPST}(\ulomega^*\G_m,\ulomega^{\CI}\ulomega_!(\tG\otCINissp\Lred)))\ar[d]^{\hookrightarrow}&\ar[l]_-{\simeq}^-{\eqref{HomFaG1}}
 \Hom_{\PST}(F\langle 1 \rangle,G\langle 1 \rangle)\\
\Hom_{\ulMPST}(\tF,\ulomega^{\CI}\uHom_{\PST}(\G_m,G\langle1\rangle))& 
\ar[l]^-{\simeq}_-{\beta}\Hom_{\ulMPST}(\tF,\tG) 
}
\]
The two right vertical isomorphisms follow from the full faithfulness of $\ulomega^{\CI}$. The isomorphism $\alpha$ (resp. $\beta$) comes from 
$\iota_{\tG}^1$ from Corollaries \ref{cor3;weakcancel} and \ref{cor;gammaGm} (resp. $\lambda_G^1$ from Proposition \ref{prop;cancelRSC}).
The squares are commutative by \eqref{iotalambda}
noting that the left vertical maps are viewed as inclusions under the identifications
	\begin{multline*}
		\ulomega_!\uHom_{\ulMPST}(\ulomega^*\G_m,\tG\otCINissp\Lred) \simeq
		\uHom_{\PST}(\G_m,G\langle1\rangle))\\
		\simeq 
		\ulomega_! \uHom_{\ulMPST}(\ulomega^*\G_m,\ulomega^{\CI}\ulomega_!(\tG\otCINissp\Lred)))
	\end{multline*}
	coming from Proposition \ref{prop;omegagamma}.
This proves that the map $\iota_{F,G}$ is an isomorphism as desired.

\end{proof}

\begin{cor}
\label{cor;cancelinjRSC}
For $F,G\in \RSC_\Nis$, there exists a natural injective map in $\NST$ for internal hom:
\begin{equation}\label{eq2;cancelRSCintro}
\uHom_{\PST}(F\langle 1 \rangle,G\langle 1 \rangle)\hookrightarrow 
\uHom_{\PST}(F,G),
\end{equation}
which coincides with the inverse of \eqref{eq;cancelRSC} on the $k$-valued points.
\end{cor}
\begin{proof}
The surjective map $F\otimes_{\NST}\G_m \to F\langle 1\rangle$ in $\NST$ from \eqref{FGm} induces an injective map
\begin{multline*}
\uHom_{\PST}(F\langle 1 \rangle,G\langle 1 \rangle)\hookrightarrow 
\uHom_{\PST}(F\otimes_{\NST}\G_m,G\langle 1 \rangle)\\
\simeq \uHom_{\PST}(F,\uHom_{\PST}(\G_m,G\langle 1 \rangle)\end{multline*}
and the latter is isomorphic to $\uHom_{\PST}(F,G)$ by Proposition \ref{prop;cancelRSC}. This completes the proof.
\end{proof}
\medbreak

\begin{proof}[Proof of Theorem \ref{cor;cancelmilnorK}]
Consider the map induced by \eqref{FGm}:\[
q: \uHom_{\PST}(\sK_n^M,F\otimes_{\NST} \sK_n^M) \to \uHom_{\PST}(\sK_n^M,F\langle n \rangle).
\]
The map \eqref{iotaFGm} is then the composition of $q$ and the map 
\begin{equation}\label{tensoridMilnor}
    F\to \uHom_{\PST}(\sK_n^M,F\otimes_{\NST} \sK_n^M);\;\; s\mapsto s\otimes id_{\sK_n^M}.
\end{equation} 
On the other hand, we have isomorphisms $ \sK_{i-1}^M\langle 1 \rangle \cong \sK_i^M$ for all $i\geq 1$ by \eqref{KMGatensor}.
Hence the map \eqref{eq2;cancelRSCintro} for $F=\sK_{i-1}^M$ gives an injective map
\begin{equation}\label{eq2;cancelRSCintroMilnorK}
    \uHom_{\PST}(\sK_i^M,F\langle i \rangle)\to \uHom_{\PST}(\sK_{i-1}^M,F\langle i-1 \rangle).
\end{equation}
Composing \eqref{eq2;cancelRSCintroMilnorK} for all $i\leq n$, we get an injective map
\begin{equation}\label{eq2;cancelRSCMilnorKn}
\uHom_{\PST}(\sK_n^M,F\langle n \rangle)\hookrightarrow F
\end{equation}
which by definition sends $q(s\otimes id_{\sK_n^M})$ to $s$ for a section $s$ of $F$. Hence the composition\[
F\xrightarrow{\eqref{iotaFGm}} \uHom_{\PST}(\sK_n^M,F\langle n \rangle)\overset{\eqref{eq2;cancelRSCMilnorKn}}{\hookrightarrow} F
\]
is the identity, so \eqref{iotaFGm} is an isomorphism, which completes the proof of Theorem \ref{cor;cancelmilnorK}.
\end{proof}

Let $G\in \RSC_{\Nis}$ and $X\in \Sm$. By Lemma \ref{lem;omegauHom} we have a natural isomorphism 
\[ \ulomega_!\uHom_{\ulMPST}((X,\emptyset),\ulomegaCI G) \simeq
\uHom_{\PST}(X, G). \]
Hence, the unit map $id\to \ulomegaCI\ulomega_!$ from \eqref{omegaCIadjointNis} induces a natural map
\begin{equation}\label{theta}
 \uHom_{\ulMPST}((X,\emptyset),\ulomegaCI G) \to 
\ulomegaCI\uHom_{\PST}(X,G).
\end{equation}
It is injective by the semipurity of 
$\uHom_{\ulMPST}(\Ztr(X,\emptyset),\ulomegaCI G) $,
and becomes an isomorphism after taking $\ulomega_!$. 
Moreover the following diagram is commutative:
\begin{equation}\label{thetaCD}
\xymatrix{
\uHom_{\ulMPST}((X,\emptyset),\ulomegaCI G) \ar[r]^-{\eqref{theta}}\ar[d]^{\hookrightarrow} &  
\ulomegaCI\uHom_{\PST}(X,G)\ar[d]^{\hookrightarrow} \\
 \uHom_{\ulMPST}((X,\emptyset),\ulomega^* G) \ar[r]^-{\simeq} &  
\ulomega^*\uHom_{\PST}(X,G)\\}
\end{equation}
where the isomorphism comes from Lemma \ref{lem0;omegauHom}.

For $G\in \RSC_{\Nis}$ and $X\in \Sm$, we define the following condition:

\begin{itemize}
    \item[$(\clubsuit)_X$] The maps \eqref{theta} is an isomorphism.
\end{itemize}

\begin{thm}
Let $F,G\in \RSC_\Nis$. Assume one of the following:
\begin{enumerate}
    \item[(a)] $G$ satisfies $(\clubsuit)_X$ for any $X\in \Sm$. 
    \item[(b)] $G$ satisfies $(\clubsuit)_{\Spec(K)}$ 
for any function field $K$ over $k$ and $F$ is the quotient of a direct sum of representable objects. 
\end{enumerate}
Then \eqref{eq2;cancelRSCintro} is an isomorphism.
\end{thm}
\begin{proof}
Assume the condition $(a)$. Letting $\tG=\ulomegaCI G$, we have isomorphisms for $X\in \Sm$
\begin{multline}\label{eq;trickclub}
\uHom_{\PST}(F,G)(X)=\Hom_{\PST}(F,\uHom_{\PST}(X,G)) \\
\underset{(*1)}{\cong}\Hom_{\ulMPST}(\tF,\ulomega^{\CI}\uHom_{\PST}(X,G))\underset{(*2)}{\cong}\Hom_{\ulMPST}(\tF,\uHom_{\ulMPST}((X,\emptyset),\tG)),
\end{multline}
where the isomorphism $(*1)$ (resp. $(*2)$) comes from the full faithfullness of $\ulomega^{\CI}$
(resp. $(\clubsuit)_X$).
Moreover, we have isomorphisms
\begin{multline}\label{eq;alphainternal}
    \uHom_{\ulMPST}((X,\emptyset),\tG)\underset{(*3)}{\cong} 
\uHom_{\ulMPST}((X,\emptyset),\uHom_{\ulMPST}(\ulomega^*\G_m,\tG(1)))\\
    \cong \uHom_{\ulMPST}(\ulomega^*\G_m,\uHom_{\ulMPST}((X,\emptyset),\tG(1))),
\end{multline}
where the isomorphism $(*3)$ comes from Corollaries \ref{cor3;weakcancel} and \ref{cor;gammaGm}.
We also have isomorphisms
\begin{multline}\label{HomFaG1internal}
    \uHom_{\PST}(F\langle 1 \rangle,G\langle 1 \rangle)(X) = \Hom_{\PST}(F\langle 1 \rangle,\uHom_{\PST}(X,G\langle 1 \rangle))\\
\underset{(*4)}{\cong} \Hom_{\PST}(\ulomega_!(\tF\otCINis \ulomega^*\G_m),\ulomega_!\uHom_{\ulMPST}((X,\emptyset),\tG(1)))\\
 \underset{(*5)}{\cong}  \Hom_{\ulMPST}(\tF\otimes_{\ulMPST} \ulomega^*\G_m,\ulomega^{\CI}\ulomega_!\uHom_{\ulMPST}((X,\emptyset),\tG(1)))\\
\cong \Hom_{\ulMPST}(\tF,\uHom_{\ulMPST}(\ulomega^*\G_m,\ulomega^{\CI}\ulomega_!\uHom_{\ulMPST}((X,\emptyset),\tG(1))) ,
\end{multline}
where $(*4)$ (resp. $(*5)$) comes from Lemma \ref{lem;omegauHom}
(resp. the adjunction \eqref{omegaCIadjoint}). These maps fit into a commutative diagram
\begin{small}
  \[\hskip -10pt
\xymatrix{
\ar[d]^{\simeq}_{\eqref{eq;alphainternal}}\Hom_{\ulMPST}(\tF,\uHom_{\ulMPST}((X,\emptyset),\tG))  \\
\Hom_{\ulMPST}(\tF,\uHom_{\ulMPST}(\ulomega^*\G_m,\uHom_{\ulMPST}((X,\emptyset),\tG(1))))\ar[d]^{\hookrightarrow}_{(\dagger)}&\uHom_{\PST}(F,G)(X)
\ar[ul]^-{\simeq}_-{\eqref{eq;trickclub}} \\
\Hom_{\ulMPST}(\tF,\uHom_{\ulMPST}(\ulomega^*\G_m,\ulomega^{\CI}\ulomega_!\uHom_{\ulMPST}((X,\emptyset),\tG(1)))))&\ar[l]_-{\simeq}^-{\eqref{HomFaG1internal}}
 \uHom_{\PST}(F\langle 1 \rangle,G\langle 1 \rangle)(X)\ar[u]^{\hookrightarrow}_{\eqref{eq2;cancelRSCintro}} \\
}
\]
\end{small}
where the injective map $(\dagger)$ comes from the counit map $id\to \ulomegaCI\ulomega_!$ from the adjunction \eqref{omegaCIadjoint}.
The diagram commutes since the map \eqref{eq;alphainternal} is induced by the map
\begin{multline*}
\uHom_{\ulMPST}((X,\emptyset),\tG)\to
\uHom_{\ulMPST}(\ulomega^*\G_m,\uHom_{\ulMPST}((X,\emptyset),\tG(1)))\\
\simeq\uHom_{\ulMPST}((X,\emptyset)\otimes \ulomega^*\G_m,\tG\otCINissp \ulomega^*\G_m)
\end{multline*}
given by $f\mapsto f\otimes id_{\ulomega^*\G_m}$, and the map \eqref{eq2;cancelRSCintro} is induced by the surjection
$F\otimes_{\NST}\G_m \to F\langle 1\rangle$ from \eqref{FGm} and the isomorphism inverse of \eqref{iotaFGm}:
\[	\uHom_{\PST}(F\otimes \G_m, G\langle 1 \rangle) \isom \uHom_{\PST}(F,G)	\] 
given by $f\otimes id_{\G_m} \mapsto f$, 
and the maps \eqref{eq;trickclub} and $(\dagger)$ are inclusions under the identifications
	\begin{multline*}
		\ulomega_!\uHom_{\ulMPST}(\ulomega^*\G_m,\uHom_{\ulMPST}(X,\emptyset),\tG(1)) \simeq
		\uHom_{\PST}(\G_m\otimes X,G\langle1\rangle))\\
		\simeq 
		\ulomega_! \uHom_{\ulMPST}(\ulomega^*\G_m,\ulomega^{\CI}\ulomega_!\uHom_{\ulMPST}((X,\emptyset),\tG\otCINissp\Lred)) 
	\end{multline*}
	coming from Lemma \ref{lem;omegauHom} and Proposition \ref{prop;omegagamma}.
This proves that \eqref{eq2;cancelRSCintro} is an isomorphism.
\medbreak

Next assume the condition $(b)$. In view of Lemma \ref{lem;homrec}, we have 
$\uHom_{\PST}(F,G)$ and $\uHom_{\PST}(F\langle 1 \rangle,G\langle 1 \rangle)$ are in $\RSC_{\Nis}$. Hence, by Lemma \ref{lem;RSCexactness}, it is enough to prove that \eqref{eq2;cancelRSCintro} induces an isomorphism
\[\uHom_{\PST}(F\langle 1 \rangle,G\langle 1 \rangle)(K) \cong \uHom_{\PST}(F,G)(K)\]
for any function field $K$ over $k$. This follows from the same computations as above.
\end{proof}

\begin{lemma}
$F\in \HI_{\Nis}$ satisfies $(\clubsuit)_X$ for all $X\in \Sm$.
\end{lemma}
\begin{proof} We have
\begin{multline*}
\uHom_{\ulMPST}((X,\emptyset),\ulomega^{\CI}F)=\uHom_{\ulMPST}((X,\emptyset),\ulomega^*F)\underset{(*1)}{\cong} \ulomega^*\uHom_{\PST}(X,F)\\
\underset{(*2)}{\cong}\ulomega^{\CI}\uHom_{\PST}(X,F),
\end{multline*}
where the isomorphism $(*1)$ follows from Lemma \ref{lem0;omegauHom} 
and $(*2)$ from the fact that $\uHom_{\PST}(X,F)\in\HI$ so that
$\ulomega^*\uHom_{\PST}(X,F))\in \CIt$ by \cite[Lem. 2.3.1]{ksyII}.
This completes the proof.
\end{proof}

\begin{lemma}\label{uHomOmega}
If $\ch(k)=0$, $\Omega^i$ satisfies $(\clubsuit)_X$ for all $X\in \Sm$.
\end{lemma}
\begin{proof}
Put $\Gamma=\uHom_{\PST}(\Ztr(X),\Omega^i)$ and
\[ G=\uHom_{\ulMPST}(\Ztr(X,\emptyset),\ulomegaCI \Omega^i),\quad
G^*=\ulomegaCI\uHom_{\PST}(\Ztr(X),\Omega^i).\]
Note that $\Gamma\in \RSC_\Nis$ by Lemma \ref{lem;homrec}. 
By \cite[Cor. 6.8]{rs}, for $\sY=(Y,D)\in \ulMCor$ where $Y\in \Sm$ and $D_{\red}$ is a simple normal crossing divisor, we have
\begin{equation}\label{eq1;uHomOmega}
G(\sY)= \Gamma(Y\times X,\Omega^i(\log D_{red}\times X)((D-D_{\red}) \times X)).
\end{equation}
Hence the conductor $c^G$ associated to $G$ in the sense of \cite[Def. 4.14]{rs}
is given as follows (note that Lemma \ref{lem;omegauHom} implies $G\in \CI(\Gamma)$ under the notation of loc. cite.): Let $\Phi$ be as \cite[Def. 4.1]{rs}. For 
\[a\in G(L)=H^0(X\otimes_k L,\Omega^i)\;\text{ with } L\in \Phi,\]
put $c_L^{G}(a)=0$ if $a\in H^0(X\otimes_k\sO_L,\Omega^i)$. Otherwise, put
\[c_L^{G}(a)=
\min\left\{n\ge 1\mid a\in H^0(X\otimes_k \sO_L, \frac{1}{t^{n-1}}\cdot
\Omega^i_{X\otimes_k \sO_L}(\log)) \right\},\]
where $t$ is a local paramter of $\sO_L$ and $\Omega^\bullet_{X\otimes_k \sO_L}(\log)$ is the differential graded subalgebra of $\Omega^\bullet_{X\otimes_k L}$ generated by $\Omega^\bullet_{X\otimes_k \sO_L}$ and $\dlog\; t$ (cf. \cite[\S6.1 6.3]{rs}). 
Moreover, one easily sees that for $\sY=(Y,D)\in \ulMCor$ as \eqref{eq1;uHomOmega}, 
\[G(\sY)= \left\{a\in G(Y-D) \mid c_L^G(a)\leq v_L(D)\;\text{ for any } L\in \Phi \right\}\]
(see \cite[Notation 4.2]{rs} for $v_L(D)$).
Hence, by \cite[Th. 4.15(4)]{rs}, it suffices to show $c^{G^*}=c^G$.
We know $c^{G^*}\leq c^G$ by loc. cite so that it suffices to show the following: Let $L\in \Phi$ and $a\in G(L)$. For $r\in \Z_{\geq 0}$, we have
\begin{equation*}\label{eq2-1;uHomOmega}
c_L^{G^*}(a)\leq r \Rightarrow  c_L^G(a)\leq r.
\end{equation*}
We prove it by the descending induction on $r$.
By \cite[Cor. 4.44]{rs} this is reduced to showing the following:
Choose a ring homomorphism $K\inj \sO_L$ such that $K\to \sO_L\to \sO_L/(t)$ is an identity and extend it in the canonical way to $\sigma: K(x)\inj \sO_{L_x}$,
where $x$ is a variable and $L_x=\Frac(\sO_L[x]_{(t)}^h)$. Assume $c_L^G(a)\leq r+1$.
Then the following implication holds
\begin{equation}\label{eq2;uHomOmega}
(a, 1-xt^{r})_{L_x,\sigma}=0\in G(K(x)) \,\Rightarrow\, c_L^G(a)\leq r,
\end{equation}
where $(-,-)_{L_x,\sigma}$ is the local symbol for $\Gamma^i=\uHom_{\PST}(\Ztr(X),\Omega^i)$ from \cite[\S4.3 4.37]{rs}.
Since the local symbol is uniquely determined by 
the properties (LS1) - (LS4) from \cite[\S4.3 4.38]{rs}, we see that it is given by 
\[(a, 1-xt^{r})_{L_x,\sigma}= \Res_t(a\; \dlog (1-xt^{r})),\]
where 
\[ \Res_t: \Gamma^{i+1}(L_x)= H^0(X\otimes_k L_x,\Omega^{i+1}) \to 
\Gamma^{i}(K(x))  =H^0(X\otimes_k K(x),\Omega^{i})\]
is induced by the residue map
$\Omega^{i+1}_{L_x}\to \Omega^i_{K(x)}$, which is defined using 
the isomorphism $L_x \simeq K(x)((t))$ induced by $\sigma: K(x)\inj \sO_{L_x}$.
To prove the implication \eqref{eq2;uHomOmega}, we may assume after
replacing $a$ by $a-b$ for some $b\in \Gamma(L)$ with $c_L^G(b)\leq r$,
\[a = \frac{1}{t^r}\alpha + \beta\frac{dt}{t^{r+1}}\;\;\text{for }
\alpha\in H^0(X\otimes_k K, \Omega^i),\;\beta\in H^0(X\otimes_k K, \Omega^{i-1}).\]
Then we compute in $H^0(X\otimes_k K(x), \Omega^i)$
\[\Res_t(a \;\dlog (1-xt^{r}))=-rx\alpha+\beta dx.\]
This shows \eqref{eq2;uHomOmega} and completes the proof.

\end{proof}

\section{Internal hom's for \texorpdfstring{$\Omega^n$}{K\"ahler differentials}}\label{intHomOmega}

\def\uPhi{\underline{\Phi}}
In this section, we assume $\ch(k)=0$. Note that a section of 
$\uHom_{\PST}(\Omega^n,\Omega^m)$ over $X\in \Sm$ is given by a collection of 
maps 
\[\phi_Y: H^0(Y,\Omega^n) \to H^0(X\times Y,\Omega^m)\;\text{ for } Y\in \Sm,\]
which are natural in $Y\in \Cor$. For 
\[(\alpha,\beta)\in H^0(X,\Omega^{m-n})\oplus H^0(X,\Omega^{m-n-1}),\]
we define 
\[ \phi_{Y,\alpha,\beta}^{n,m} :  H^0(Y,\Omega^n) \to H^0(X\times Y,\Omega^m)\;;\;
\omega \to p_X^*\alpha\wedge p_Y^*\omega + p_X^*\beta\wedge p_Y^*d\omega,\]
where $p_X: X\times Y \to X$ and $p_Y: X\times Y \to Y$ are the projections.
The naturalness of $\phi_{Y,\alpha,\beta}^{n,m} $ in $Y\in \Cor$ follows from \cite{cr}.
Thus we get a natural map in $\NST$:
\begin{equation}\label{uPhi}
\Omega^{m-n}\oplus \Omega^{m-n-1} \to
\uHom_{\PST}(\Omega^n,\Omega^m)\;;\;
(\alpha,\beta)\to \{\phi_{Y,\alpha,\beta}^{n,m}\}_{Y\in \Sm},
\end{equation}
where $\Omega^i=0$ for $i<0$ by convention.
Taking the sections over $\Spec k$, we get a natural map
\begin{equation}\label{Phi}
\Phi^{n,m} : 
\Omega^{m-n}_k\oplus \Omega^{m-n-1}_k \to
\Hom_{\PST}(\Omega^n,\Omega^m) .
\end{equation}
We also consider the composite map in $\NST$:
\begin{equation}\label{uPsi}
\Omega^{m-n} \rmapo{\eqref{uPhi}} \uHom_{\PST}(\Omega^n,\Omega^m)
\rmapo{\dlog^*} \uHom_{\PST}(\sK_n^M,\Omega^m),
\end{equation}
where the second map is induced by the map $\dlog: \sK_n^M \to \Omega^n$.
Taking the sections over $\Spec k$, we get a natural map
\begin{equation}\label{Psi}
\Psi^{n,m} : 
\Omega^{m-n}_k \to \Hom_{\PST}(\sK_n^M,\Omega^m) .
\end{equation}

The main result of this subsection is the following.

\begin{thm}\label{thm;cancelOmega}
The maps \eqref{uPhi} and \eqref{uPsi} are isomorphisms.
\end{thm}

First we prove the following.

\begin{prop}\label{pr;cancelOmega}
The maps \eqref{Phi} and \eqref{Psi} are isomorphisms.
\end{prop}

This follows from Lemmas
\ref{lem;Phicompatibiity}, \ref{lem;homgaomegavanishing} and \ref{lem;homgaomega} below, in light of Theorem \ref{thm;closeenoughtocancel}.
For $i\geq 0$, let us fix the isomorphisms
\begin{equation}\label{eq;sigma}
\sigma^i: \Omega^{i-1}\langle 1\rangle\isom \Omega^i ,\;\; \varsigma^i:\sK_{i-1}^M\langle 1\rangle\isom \sK_i^M
\end{equation}
coming from \eqref{twistrecursive} and \eqref{KMGatensor}

\begin{lemma}\label{lem;Phicompatibiity}

\begin{itemize}
\item[(1)]
The following diagram is commutative:
\[\xymatrix{
\Omega^{m-n}_k\oplus \Omega^{m-n-1}_k \ar[r]^-{\Phi^{n,m}}\ar[d]^{\Phi^{n-1,m-1}}
& \Hom_{\PST}(\Omega^n,\Omega^m) \\
\Hom_{\PST}(\Omega^{n-1},\Omega^{m-1}) \ar[r]^-{\eqref{eq;cancelRSC}} & 
\Hom_{\PST}(\Omega^{n-1}\langle 1\rangle,\Omega^{m-1}\langle 1\rangle) 
\ar[u]
}\]
where the right vertical map is induced by $\sigma^m$ and $(\sigma^n)^{-1}$ of \eqref{eq;sigma}.
\item[(2)]
The following diagram is commutative:
\[\xymatrix{
\Omega^{m-n}_k \ar[r]^-{\Psi^{n,m}}\ar[d]^{\Psi^{n-1,m-1}}
& \Hom_{\PST}(\sK_n^M,\Omega^m) \\
\Hom_{\PST}(\sK_{n-1}^M,\Omega^{m-1}) \ar[r]^-{\eqref{eq;cancelRSC}} & 
\Hom_{\PST}(\sK_{n-1}^M\langle 1\rangle,\Omega^{m-1}\langle 1\rangle) 
\ar[u]
}\]
where the right vertical map is induced by $\sigma^m$ and $(\varsigma^n)^{-1}$ of \eqref{eq;sigma}.
\end{itemize}
\end{lemma}
\begin{proof}
By \cite[Cor. 5.22]{rsy}, for an affine $X=\Spec A\in \Sm$ and $i\geq 0$, the composite map
\[ \theta^i : \Omega^{i-1}_A\otimes_\bZ A^\times \to 
(\Omega^{i-1}\otimes_{\NST} \G_m)(A) \rmapo{\eqref{FGm}} 
\Omega^{i-1}\langle 1\rangle(A) \rmapo{\sigma^i} \Omega^i_A\]
sends $\omega\otimes f$ with $\omega\in \Omega^{i-1}_A$ and $f\in A^\times$ to
$\omega \wedge \dlog f$. Moreover, for $\phi\in \Hom_{\PST}(\Omega^{n-1},\Omega^{m-1})$ and $\phi'=\sigma^m\circ\phi\langle 1\rangle\circ(\sigma^n)^{-1}$, the diagram
\[\xymatrix{
\Omega^{n-1}_A\otimes_\bZ A^\times \ar[r]^{\hskip 20pt \theta^n}  \ar[d]^{\phi\otimes id_{A^\times}} 
& \Omega^n_A\ar[d]^{\phi'} \\
\Omega^{m-1}_A\otimes_\bZ A^\times \ar[r]^{\hskip 20pt\theta^m} & \Omega^m_A\\ }\]
is commutative. Hence (1) follows from the equation
\[ \alpha\wedge (\omega \wedge \dlog f) + \beta\wedge d(\omega \wedge \dlog f)
= (\alpha\wedge \omega + \beta\wedge d\omega) \wedge \dlog f,\]
where $\alpha\in \Omega^{m-n}_k$ and $\beta\in \Omega^{m-n-1}_k$.

(2) follows from (1) and the commutativity of the diagram
\[\xymatrix{
\sK_{n-1}^M\langle 1\rangle\ar[r]^-{\dlog\langle 1\rangle} \ar[d]^{\varsigma^n}
& \Omega^{n-1}\langle 1\rangle\ar[d]^{\sigma^n}\\
\sK_n^M \ar[r]^-{\dlog} & \Omega^n \\}\]
which can be verified using \eqref{GaGmOmega}.
\end{proof}

\begin{lemma}\label{lem;homgaomegavanishing}
For an integer $n\geq 1$, we have 
\begin{equation}\label{OmegaKMGatensor}
 \Hom_{\PST}(\Omega^n,\G_a)= \Hom_{\PST}(\mathcal{K}^M_n,\G_a)=0.
\end{equation}
\end{lemma}
\begin{proof}
 We have isomorphisms
\begin{multline*}
\Hom_{\PST}(\Omega^n,\G_a)\simeq 
\Hom_{\PST}(\ulomega_!(\widetilde{\Omega^{n-1}}\otCI\ulomega^*\G_m),\G_a)\\
\simeq \Hom_{\ulMPST}(\widetilde{\Omega^{n-1}}\otCI\ulomega^*\G_m,\ulomega^{\CI}\G_a)\\
\simeq \Hom_{\ulMPST}(\widetilde{\Omega^{n-1}}\otimes_{\ulMPST}\ulomega^*\G_m,\ulomega^{\CI}\G_a)\\
\simeq \Hom_{\ulMPST}(\widetilde{\Omega^{n-1}},\uHom_{\ulMPST}(\ulomega^*\G_m,\ulomega^{\CI}\G_a)).\\
\end{multline*}
where the first isomorphism is induced by $(\sigma^n)^{-1}$, inverse of \eqref{eq;sigma}, and the second from 
\eqref{omegaCIadjoint}.
Similarly we have an isomorphism using $(\varsigma^n)^{-1}$ instead of $(\sigma^n)^{-1}$ 
\[\Hom_{\PST}(\mathcal{K}^M_n,\G_a)\simeq 
\Hom_{\ulMPST}(\ulomega^*\mathcal{K}^M_{n-1},\uHom_{\ulMPST}(\ulomega^*\G_m,\ulomega^{\CI}\G_a)).\]
We compute 
\begin{multline*} \uHom_{\ulMPST}(\ulomega^*\G_m,\ulomega^{\CI}\G_a)(X)\hookrightarrow \uHom_{\ulMPST}(\ulomega^*\G_m,\ulomega^{\CI}\G_a)(K(X))\\
\simeq\Coker\big(\ulomega^{\CI}\G_a(K(X)) \to \ulomega^{\CI} \G_a(\P^1_{K(X)},0+\infty)\big)\\
\simeq \Coker\big(K(X) \to  H^0(\P^1_{K(X)},\sO)\big)=0\\
\end{multline*}
where the first map is injective by \cite[Corollary 0.3]{shuji}, and the first (resp. last) isomorphism follows from Corollary \ref{cor;gammaGm}(1) (resp. \cite[Cor. 6.8]{rs}).
This completes the proof of Lemma \ref{lem;homgaomegavanishing}.
\end{proof}

\begin{lemma}\label{lem;homgaomega}
The maps \eqref{Phi} and \eqref{Psi} are isomorphisms for $n=0$.
\end{lemma}
\begin{proof}
The assertion for \eqref{Psi} is obvious since $\sK_n^M=\Z$ for $n=0$.
We prove it for \eqref{Phi}.
We have isomorphisms
\begin{multline}\label{eq2}
 \Hom_{\PST}(\G_a,\Omega^i) \simeq 
\Hom_{\PST}(\aVNis\omega_! h_0^\bcube(\cubega),  \Omega^i)\\
\simeq \Hom_{\MPST}(h_0^\bcube(\cubega), \omegaCI \Omega^i)\\
\simeq \Hom_{\MPST}(\cubega, \omegaCI \Omega^i)\\
\simeq \Ker\big(H^0(\P^1,\Omega^i_{\P^1}(\log \infty)(\infty))\rmapo{i_0^*} 
\Omega^i_k\big),\\
\end{multline}
where the first (resp. last) isomorphism follows from \eqref{Gah0cube}
(resp. \cite[Cor. 6.8]{rs}).
Since $\Omega^i_{\P^1/k}(\log\infty)=0$ for $i>1$, 
$\sO_{\P^1}(\log\infty)=\sO_{\P^1}$ and $\Omega^1_{\P^1/k}(\log\infty)=\Omega^1_{\P^1/k}(\infty)$, the standard exact sequence
\[ 0\to \sO_{\P^1}\otimes_k\Omega^1_k \to \Omega^1_{\P^1}(\log\infty) \to \Omega^1_{\P^1/k}(\log\infty)\to 0\]
induces an exat sequence
\[ 0\to \sO_{\P^1}(\infty)\otimes_k\Omega^i_k \to \Omega^i_{\P^1}(\log\infty)(\infty) \to 
\Omega^1_{\P^1/k}(2\infty) \otimes_k \Omega_k^{i-1}\to 0,\]
where $\Omega^{i-1}_k=0$ if $i=0$ by convention.
Letting $t$ be the standard coordinate of $\A^1\subset \P^1$, we have 
\[ H^0(\P^1,\sO_{\P^1}(\infty)) = k\cdot 1 \oplus k\cdot t,\quad
H^0(\P^1,\Omega^1_{\P^1/k}(2\infty)) = k\cdot dt,\]
and $dt$ lifts canonically to a section $dt\in H^0(\P^1, \Omega^1_{\P^1}(\log\infty)(\infty)) $. Hence we get an isomorphism
\begin{multline}\label{eq2-1}
H^0(\P^1,\Omega^i_{\P^1}(\log \infty)(\infty))\simeq
(k\cdot 1\oplus k\cdot t)\otimes_k \Omega^{i}_{k} \oplus (k\cdot dt)\otimes_k \Omega^{i-1}_{k}. 
\end{multline}
Thus the last group of \eqref{eq2} is isomorphic to 
\[k\cdot t\otimes_k \Omega^{i}_{k} \oplus k\cdot dt\otimes_k \Omega^{i-1}_{k}
\simeq \Omega^{i}_{k}\oplus \Omega^{i-1}_{k}.\]
Hence, from \eqref{eq2}, we get a natural isomorphism
\begin{equation}\label{HomGHaOmega}
 \Omega^{i-1}_k\oplus \Omega^i_k \isom \Hom_{\PST}(\G_a,\Omega^i).
\end{equation}
\medbreak

Next we claim that the map \eqref{HomGHaOmega} coincides with \eqref{Phi} for $n=0$. By Lemma \ref{lem;GmGa}(2), we have a commutative diagram
\begin{equation}\label{eq1;lem;explicit}
\xymatrix{
\Ztr(\A^1_t) \ar[r]^{\lamga} \ar[d]^\simeq & \G_a \\
\omega_!\Ztr(\P^1,2\infty)\ar[r]&\omega_!h_0^\bcube(\cubega) 
\ar[u]_{\eqref{Gah0cube}}\\}\end{equation}
where $\lamga$ is given by $t\in \G_a(\A^1_t)=k[t]$. The standard isomorphism 
\[\Omega^i(\A_t^1) \simeq (\Omega^i_k\otimes_k k[t])  
\oplus (\Omega^{i-1}_k\otimes_k k[t] dt)\]
induces a natural  isomorphism
\begin{equation}\label{HomGaOmega}
\Hom_{\PST}(\Ztr(\A^1_t),\Omega^i)=\Omega^i(\A^1_t)\simeq \Omega^i_k[t]\oplus \Omega^{i-1}_k[t]\wedge dt,\end{equation}
where 
\[  \Omega^i_k[t] =\underset{m\in \Z_{\geq 0}}{\bigoplus} \Omega^i_k \cdot t^m,
\quad \Omega^{i-1}_k[t]\wedge dt = 
\underset{m\in \Z_{\geq 0}}{\bigoplus} \Omega^{i-1}_k\wedge t^m dt.\]
The map $\lamga$ induces the inclusion 
\[\lamga^*: \Hom_{\PST}(\G_a,\Omega^i)\hookrightarrow
 \Hom_{\PST}(\Ztr(\A^1_t),\Omega^i)=\Omega^i(\A_t^1)\]
such that
\begin{equation}\label{eq;lamgaformula}
\lamga^*(\phi) = \phi_{\A_t^1}(t)\qfor \phi\in \Hom_{\PST}(\G_a,\Omega^i),
\end{equation}
where $\phi_{\A_t^1}: \G_a(\A^1_t) =k[t] \to \Omega^i(\A_t^1)$ is induced by $\phi$. The following claim follows from \eqref{eq2}, \eqref{eq2-1} and \eqref{eq1;lem;explicit}.

\begin{claim}\label{claim;explicit}
The image of $\lamga^*$ is identified under \eqref{HomGaOmega} with
\[\Omega^i_k\cdot t \oplus \Omega^{i-1}_k \wedge dt\subset
\Omega^i_k[t]\oplus \Omega^{i-1}_k[t]\wedge dt,\]
and the composite map
\[ \Omega^i_k\oplus\Omega^{i-1}_k \rmapo{\eqref{HomGHaOmega}} 
\Hom_{\PST}(\G_a,\Omega^i)
\rmapo{\lamga^*} \Omega^i_k\cdot t \oplus \Omega^{i-1}_k \wedge dt\]
is given by the obvious identifications 
$\Omega^i_k =\Omega^i_k\cdot t$ and $\Omega^{i-1}_k =\Omega^{i-1}_k\wedge dt$.
\end{claim}

Let
\begin{equation}\label{eq;lem;explicit}
\Hom_{\G_a}(\G_a, \Omega^i)\subset \Hom_{\PST}(\G_a, \Omega^i)
\end{equation} 
be the subgroup of $\G_a$-linear morphisms.
There is a natural isomorphism\[ 
\xi: \Omega^i_k\cong \Hom_{\G_a}(\G_a, \Omega^i)\;;\quad \omega \mapsto \{\lambda\mapsto \lambda \omega\} \; (\lambda\in \G_a).\]
\eqref{eq;lem;explicit} is a direct summand since we have a splitting given by
\[\Hom_{\PST}(\G_a, \Omega^i) \to \Hom_{\G_a}(\G_a, \Omega^i)\;;\;
\phi \mapsto \{\lambda\mapsto \lambda \phi(1)\}.\]
The other summand is
\[\Hom_{\PST}(\G_a,\Omega^i)^{0}:=\{\phi|\;  \phi(1)=0\}.\]
There is a natural map
\[\xi':\Omega^{i-1}_k\to \Hom_{\PST}(\G_a,\Omega^i)^{0}\; ;\quad \omega \mapsto \{\alpha \mapsto \omega \wedge d\alpha\}.\]
By \eqref{eq;lamgaformula}, under the identification \eqref{HomGaOmega}, we have
\[\lamga^*(\xi(\omega))=\omega\cdot t,\;\; \lamga^*(\xi'(\eta))=\eta\wedge dt
\quad (\omega\in \Omega^i,\; \eta\in \Omega^{i-1}).\]
Hence the composite map
\[ \Omega^i_k\oplus\Omega^{i-1}_k \rmapo{\xi\oplus\xi'} 
\Hom_{\PST}(\G_a,\Omega^i)
\rmapo{\lamga^*} \Omega^i_k\cdot t \oplus \Omega^{i-1}_k \wedge dt\]
is given by the obvious identifications 
$\Omega^i_k =\Omega^i_k\cdot t$ and $\Omega^{i-1}_k =\Omega^{i-1}_k\wedge dt$.
By Claim \ref{claim;explicit} this proves the desired claim and completes the proof of Lemma \ref{lem;homgaomega}.
\end{proof}
\bigskip

To deduce Theorem \ref{thm;cancelOmega} from Proposition \ref{pr;cancelOmega},
we need some preliminaries. 

Let $K$ be the function field of $S\in \Sm$ and 
define $\Cor_K$, $\PST_K$, $\ulMCor_K$, $\ulMPST_K$, etc. defined as 
$\Cor$, $\PST$, $\ulMCor$, $\ulMPST$, etc. where the base field $k$ is replaced by $K$. We have then a map
\begin{equation}\label{rKHom}
 r_K : \Hom_{\PST_K}(\Omega^n,\Omega^m) \to 
\uHom_{\PST}(\Omega^n,\Omega^m)(K)\;;\; \phi \to \{\psi_Y\}_{Y\in \Sm},
\end{equation}
where $\psi_Y$ for $Y\in \Sm$ is the composite map
\[ H^0(Y,\Omega^n) \to H^0(Y\times_k K, \Omega^n)  \to 
H^0(Y\times_k K, \Omega^m),\]
where the second map is $\phi_{Y\times_k K}$ (note $Y\times_k K\in \Sm_K$)
and the first is the pullback by the projection $p_Y: Y\times_k K \to Y$.
Similarly we can define a map
\begin{equation}\label{rKHom2}
 r_K : \Hom_{\PST_K}(\sK_n^M,\Omega^m) \to \uHom_{\PST}(\sK_n^M,\Omega^m)(K).
\end{equation}
By definitions, the following diagrams are commutative.
\[\xymatrix{
\Omega^{m-n}_K\oplus \Omega^{m-n-1}_K \ar[r]^-{\eqref{Phi}} 
\ar[rd]_-{\eqref{uPhi}} 
& \Hom_{\PST_K}(\Omega^n,\Omega^m)  \ar[d]^{r_K}\\ 
&\uHom_{\PST}(\Omega^n,\Omega^m)(K)\\}\]
\[\xymatrix{
\Omega^{m-n}_K \ar[r]^-{\eqref{Psi}} \ar[rd]_-{\eqref{uPsi}} 
& \Hom_{\PST_K}(\sK_n^M,\Omega^m)  \ar[d]^{r_K}\\ 
&\uHom_{\PST}(\sK_n^M,\Omega^m)(K)\\}\]
In view of Lemma \ref{lem;RSCexactness}, Theorem \ref{thm;cancelOmega} follows from Proposition \ref{pr;cancelOmega} and the following.
\def\cubeGa{\bcube_{\G_a}}
\def\cubeGm{\bcube_{\G_m}}
\def\cubeomegan{\bcube_{\Omega^n}}
\def\cubeomeganK{\bcube_{\Omega^n,K}}

\begin{lemma}\label{rK}
The maps \eqref{rKHom} and \eqref{rKHom2} are isomorphisms.
\end{lemma}

For the proof we need the following. Recall from Conventions that for $U=\lim_i U_i\in \widetilde{\Sm}$ and $F\in \PST$ we let\[
F(U):= \colim_i F(U_i).
\]
In general for $(Y,D_Y)\in \MCor$ and $F\in \ulMPST$, we let \[\uHom_{\ulMPST}(U,F)(Y,D_Y) := \colim_i F(U_i\times_k Y,U_i\times_k D_Y),\]
and for $G\in \ulMPST$ we have\[
\Hom(G,\uHom(U,F)) = \colim_i \Hom(G,\uHom(U_i,F)).
\]
\begin{lemma}\label{lem0;rK}
For $\sX=(X,D)\in \MCor$ and $\sX_K=(X_K,D_K)\in \MCor(K)$ with $X_K=X\times_k K$ and $D_K=D\times_k K$, we have a natural isomorphism
\begin{multline*}
    \Hom_{\ulMPST_K}(\Ztr(\sX_K), \ulomega^{\CI_K}\Omega^n)\cong  
\Hom_{\ulMPST}(\Ztr(\sX),\uHom_{\ulMPST}(K, \ulomega^{\CI}\Omega^n)).
\end{multline*}
\end{lemma}
\begin{proof}
By \cite[Pr. 1.9.2 c)]{kmsyI} and resolutions of singularities (recall that we are assuming $\ch(k)=0$) we may assume $X\in \Sm$ and $D_{\red}$ is a simple normal crossing divisor.
From the explicit computation of $\ulomegaCI\Omega^n$ in \cite[Cor. 6.8]{rs}, 
\begin{multline*}\label{eq;lem0;rK}
(\ulomega^{\CI_K}\Omega^n)(X_K,D_K)=
H^0(X_K,\Omega^n_{X_K}(\log(D_K))(D_K-D_{K,red}))\\
=(\ulomega^{\CI}\Omega^n)(X_K,D_K):
= \colim_{U\subset S} (\ulomega^{\CI}\Omega^n)(X\times_k U,D\times_k U),
\end{multline*}
where $U$ ranges over the open subsets of $S$.
This proves the lemma.
\end{proof}
\medbreak

We now prove Lemma \ref{rK}. 
We only prove the assertion for \eqref{rKHom}. The proof for \eqref{rKHom2} is similar. Put
\[ \cubeomegan=\cubeGa\otimes_{\MPST} \cubegmm^{\otimes n},\]
where $\cubeGa$ and $\cubegmm$ are from Lemma \ref{lem;GmGa}.
By \eqref{Gmh0cube} and \eqref{Gah0cube} and \cite[Theorem 5.20]{rsy},
we have an isomorphism in $\PST$:
\begin{equation}\label{Omegah0cube}
\aVNis \omega_!  h_0^\bcube(\cubeomegan) \isom \Omega^n.
\end{equation}
Let $\bcube_K=(\P^1_K,\infty)\in \MCor_K$ and $\cubeomeganK\in \MPST_K$ be
defined as $\cubeomegan$. 
We have isomorphisms
\begin{multline}\label{eq1;rK}
 \Hom_{\PST_K}(\Omega^n,\Omega^m) 
\simeq  \Hom_{\PST_K}(\omega_!  h_0^{\bcube_K}(\cubeomeganK),\Omega^m) \simeq \\
\Hom_{\ulMPST_K}(\cubeomeganK,\ulomega^{\CI_K}\Omega^m) 
\simeq\Hom_{\ulMPST}(\cubeomegan,\uHom_{\ulMPST}(K,\ulomega^{\CI}\Omega^m) ),
\end{multline}
where the last isomorphism comes from Lemma \ref{lem0;rK}.
On the other hand, by \eqref{Omegah0cube} and Lemma \ref{lem;homrec}, we have
$\uHom(\Z_{\tr}(U),\Omega^m)\in\RSC_{\Nis}$ for $U\in \Sm$. Hence, writing $\Spec(K) = \lim_i U_i$ with $U_i\in \Sm$, we have isomorphisms (see Conventions)
\begin{multline}\label{eq2;rK}
 \uHom_{\PST}(\Omega^n,\Omega^m)(K)= 
\colim_i\Hom_{\PST}(U_i,\uHom_{\PST}(\Omega^n,\Omega^m))\simeq\\
\colim_i\Hom_{\PST}(\Omega^n,\uHom_{\PST}(U_i,\Omega^m))\simeq \colim_i \Hom_{\PST}(\omega_!  h_0^{\bcube}(\cubeomegan),\uHom_{\PST}(U_i,\Omega^m))\\\simeq 
\Hom_{\ulMPST}(\cubeomegan,\ulomegaCI\uHom_{\PST}(K,\Omega^m)) .
\end{multline}
Hence Lemma \ref{rK} follows from Lemma \ref{uHomOmega} and the following.

\begin{claim}\label{claim;rK}
The following diagram is commutative.
\begin{equation}\label{eq;claim;rK}
\xymatrix{
\Hom_{\PST_K}(\Omega^n,\Omega^m) \ar[r]^-{\eqref{eq1;rK}}\ar[d]^{r_K}
& \Hom_{\ulMPST}(\cubeomegan,\uHom_{\ulMPST}(K,\ulomega^{\CI}\Omega^m) )
\ar[d]\\
\uHom_{\PST}(\Omega^n,\Omega^m)(K) \ar[r]^-{\eqref{eq2;rK}}
& \Hom_{\ulMPST}(\cubeomegan,\ulomegaCI\uHom_{\PST}(K,\Omega^m)) \\}
\end{equation}
where the right vertical map is induced by the map \eqref{theta}.
\end{claim}
\def\AOmegan{\A_{\Omega^n}}
\def\AOmeganK{\A_{\Omega^n,K}}
To show the above claim, write $\AOmegan=\A^1\times (\A^1-\{0\})^n$ and 
$\AOmeganK=\AOmegan\otimes_k K$.
Take the standard coordinates $y$ on $\A^1$ and $(x_1,\dots,x_n)$ on 
$(\A^1-\{0\})^n$ so that 
\[\AOmegan=\Spec k[y,x_1,\dots,x_n][x_1^{-1},\dots.x_n^{-1}].\]
By the definition of $\cubeomegan$, we have natural maps in $\ulMPST$
\begin{equation}\label{eq1;claim;rK}
\Ztr(\AOmegan,\emptyset) \to
(\P^1,2\infty) \otimes (\P^1,0+\infty)^{\otimes n} \to \cubeomegan, 
\end{equation}
which induces a map in $\PST$:
\begin{equation}\label{eq1-1;claim;rK}
\lambda_{\Omega^n} : \Ztr(\AOmegan) \to \omega_! \cubeomegan \to \Omega^n,
\end{equation}
where the last map is induced by \eqref{Omegah0cube}. Let
\begin{equation}\label{eq1-2;claim;rK}
\lambda_{\Omega^n,K} : \Ztr(\AOmeganK) \to \Omega^n
\end{equation}
be defined as \eqref{eq1-1;claim;rK} replacing $k$ by $K$. 
By the definition of $\lamgm$ and $\lamga$ (cf. Lemma \ref{lem;GmGa})
and \eqref{GaGmOmega}, $\lambda_{\Omega^n}$ corresponds to 
\begin{equation}\label{omega0}
\omega_0:= y \frac{dx_1}{x_1}\wedge\cdots \wedge \frac{dx_n}{x_n}\in
\Omega^n(\AOmegan).
\end{equation}
The map \eqref{eq1;claim;rK} induces injective maps
\begin{equation}\label{eq2;claim;rK}
\Hom_{\ulMPST}(\cubeomegan,\uHom_{\ulMPST}(K,\ulomega^{\CI}\Omega^m) )
\hookrightarrow H^0(\AOmeganK,\Omega^m),
\end{equation}
\begin{equation}\label{eq2b;claim;rK}
 \Hom_{\ulMPST}(\cubeomegan,\ulomegaCI\uHom_{\PST}(K,\Omega^m)) 
\hookrightarrow H^0(\AOmeganK,\Omega^m),
\end{equation}
which are compatible with the right vertical map in \eqref{eq;claim;rK}
since applying $\ulomega_!$, the map \eqref{theta} is identified with the identity on 
$\uHom_{\PST}(K,\Omega^m)$ via the isomorphism in Lemma \ref{lem;omegauHom}.
Hence it suffices to show the commutativity of the diagram
\begin{equation}\label{eq3;claim;rK}
\xymatrix{
\Hom_{\PST_K}(\Omega^n,\Omega^m) \ar[r]^{\alpha} \ar[d]^{r_K}
& H^0(\AOmeganK,\Omega^m) \\
\uHom_{\PST}(\Omega^n,\Omega^m)(K) \ar[ur]^{\beta}\\}
\end{equation}
where $\alpha$ (resp. $\beta$) is the composite of \eqref{eq1;rK} and \eqref{eq2;claim;rK} (resp. \eqref{eq2;rK} and \eqref{eq2b;claim;rK}).
By the definition, $\alpha$ is induced by 
the map $\lambda_{\Omega^n,K}$ from \eqref{eq1-2;claim;rK}.
As $\lambda_{\Omega^n,K}$ is given by the image $\omega_{0,K}$ of $\omega_0$ from \eqref{omega0} under the pullback map
$p^*:  \Omega^n(\AOmegan) \to \Omega^n(\AOmeganK)$, we have
\begin{equation*}\label{eq4-1;claim;rK} 
\alpha(\phi) = \phi_{\AOmeganK}(\omega_{0,K})\qfor \phi\in\Hom_{\PST_K}(\Omega^n,\Omega^m)  ,
\end{equation*}
where 
$ \phi_{\AOmeganK}: \Omega^n(\AOmeganK) \to \Omega^m(\AOmeganK) $ is induced by $\phi$.
On the other hand, by the definition of $\beta$, we have a commutative diagram
\begin{equation*}
\xymatrix{
H^0(\AOmeganK,\Omega^m)\ar[r]^-{\simeq} &  
\Hom_{\PST}(\AOmegan,\uHom_{\PST}(K,\Omega^m))\\
\uHom_{\PST}(\Omega^n,\Omega^m)(K) \ar[u]^{\beta}\ar[r]^-{\simeq}&
\Hom_{\PST}(\Omega^n,\uHom_{\PST}(K,\Omega^m))\ar[u]_{\lambda_{\Omega^n}^*}\\}
\end{equation*}
where $\lambda_{\Omega^n}^*$ is induced by $\lambda_{\Omega^n}$ from \eqref{eq1-1;claim;rK}. Hence we have
\begin{equation*}\label{eq4-2;claim;rK} 
 \beta(\psi) = \psi_{\AOmegan}(\omega_{0})\qfor \psi\in\uHom_{\PST}(\Omega^n,\Omega^m)(K),
\end{equation*}
where 
$ \psi_{\AOmegan}: \Omega^n(\AOmegan) \to \uHom_{\PST}(K,\Omega^m)(\AOmegan)=\Omega^m(\AOmeganK) $ is induced by $\psi$. Then, for
$\phi\in \Hom_{\PST_K}(\Omega^n,\Omega^m)$, we get
\[ \beta(r_K(\phi))= r_K(\phi)_{\AOmegan}(\omega_0)=
\phi_{\AOmeganK}(p^*\omega_0)= \phi_{\AOmeganK}(\omega_{0,K})=\alpha(\phi),\]
which proves the commutativity of \eqref{eq3;claim;rK}. This concludes the proof of Claim \ref{claim;rK}, and hence of Lemma \ref{lem0;rK} and therefore also the proof of Theorem \ref{thm;cancelOmega}.

\end{document}